\documentclass[reqno, 11pt]{amsart}
\usepackage[left=2cm, right=2cm, top=2cm, bottom=2cm]{geometry}
\usepackage[T1]{fontenc}
\usepackage[utf8]{inputenc}	
\usepackage[english]{babel}
\usepackage[babel]{csquotes}
\usepackage{amssymb}
\usepackage{amsmath}
\usepackage{amsthm}
\usepackage{latexsym}
\usepackage{xcolor}
\usepackage{mathrsfs}
\usepackage{bbm}
\usepackage{verbatim}
\usepackage[colorlinks=true, pdfstartview=FitV, linkcolor=blue, citecolor=blue, urlcolor=blue]{hyperref}
\usepackage[style=numeric, sorting=nty, sortcites=true, backend=bibtex, maxbibnames=10, giveninits=true]{biblatex}
\usepackage{faktor}
\usepackage[shortlabels]{enumitem}
\usepackage{macros}

\numberwithin{equation}{section}
\newtheorem{thm}{Theorem}[section]
\newtheorem{cor}[thm]{Corollary}
\newtheorem{lem}[thm]{Lemma}
\newtheorem{prop}[thm]{Proposition}
\theoremstyle{definition}
\newtheorem{defin}[thm]{Definition}
\newtheorem{remark}[thm]{Remark}

\newtheorem{condition}[thm]{Condition} 

\makeatletter
\DeclareFontFamily{OMX}{MnSymbolE}{}
\DeclareSymbolFont{MnLargeSymbols}{OMX}{MnSymbolE}{m}{n}
\SetSymbolFont{MnLargeSymbols}{bold}{OMX}{MnSymbolE}{b}{n}
\DeclareFontShape{OMX}{MnSymbolE}{m}{n}{
	<-6>  MnSymbolE5
	<6-7>  MnSymbolE6
	<7-8>  MnSymbolE7
	<8-9>  MnSymbolE8
	<9-10> MnSymbolE9
	<10-12> MnSymbolE10
	<12->   MnSymbolE12
}{}
\DeclareFontShape{OMX}{MnSymbolE}{b}{n}{
	<-6>  MnSymbolE-Bold5
	<6-7>  MnSymbolE-Bold6
	<7-8>  MnSymbolE-Bold7
	<8-9>  MnSymbolE-Bold8
	<9-10> MnSymbolE-Bold9
	<10-12> MnSymbolE-Bold10
	<12->   MnSymbolE-Bold12
}{}

\let\llangle\@undefined
\let\rrangle\@undefined
\DeclareMathDelimiter{\llangle}{\mathopen}%
{MnLargeSymbols}{'164}{MnLargeSymbols}{'164}
\DeclareMathDelimiter{\rrangle}{\mathclose}%
{MnLargeSymbols}{'171}{MnLargeSymbols}{'171}
\makeatother
\AtEveryBibitem{%
\clearfield{doi}%
\clearfield{url}%
\clearfield{issn}%
\clearfield{isbn}
} 
\renewbibmacro*{publisher+location+date}{%
\ifentrytype{book}
{\printlist{publisher}%
\iflistundef{location}
{\setunit*{\addcomma\space}}
{\setunit*{\addcomma\space}}%
\printlist{location}}
{\printlist{location}%
\iflistundef{publisher}
{\setunit*{\addcomma\space}}
{\setunit*{\addcomma\space}}%
\printlist{publisher}}
\setunit*{\addcomma\space}%
\usebibmacro{date}%
\newunit} 

\usepackage[]{amsaddr}
\newcommand{\review}[1]{{\color{black}#1}}
\begin{document}
    \title[Invariant measures for the
    stochastic Allen--Cahn--Navier--Stokes system]
    {Existence, uniqueness and asymptotic stability of invariant measures
    for the stochastic Allen--Cahn--Navier--Stokes system with singular potential}
    
    \author{Andrea Di Primio, Luca Scarpa \and Margherita Zanella}
    \address{Dipartimento di Matematica,
        Politecnico di Milano, Via E.~Bonardi 9, 20133 Milano, Italy}
    \email{andrea.diprimio@polimi.it}
    \email{luca.scarpa@polimi.it}
    \email{margherita.zanella@polimi.it}
    
    \begin{abstract}
    We study the long-time behaviour of a stochastic Allen-Cahn-Navier-Stokes system modelling the dynamics of binary mixtures of immiscible fluids. The model features two stochastic forcings, one on the velocity in the Navier-Stokes equation and one on the phase variable in the Allen-Cahn equation, and includes the thermodynamically-relevant Flory-Huggins logarithhmic potential. 
    We first show existence of ergodic invariant measures and characterise their support by exploiting ad-hoc regularity estimates and suitable Feller-type and Markov properties. Secondly, we prove that if the noise acting in the Navier-Stokes 
    equation is non-degenerate along a sufficiently large number of low modes, and the Allen-Cahn equation is highly dissipative, then the stochastic flow admits a unique invariant measure and is asymptotically stable with respect to a suitable Wasserstein metric.
    \end{abstract}
    \maketitle


    \section{Introduction}
	The dynamics of complex immiscible (or partially miscible) fluids attracts significant interest in many research areas. This is mainly due to its importance in several practical applications, ranging from biology to materials science. A peculiar feature of these physical processes is that the interface between two different fluid species can exhibit a variety of topological transitions, possibly giving rise to interesting phenomena such as, for instance, pattern formation. 
    
    The mathematical understanding of models for multicomponent systems has remarkably advanced during the last decades. In particular, several contributions have been devoted to investigating the so-called diffuse interface approach, that is, studying models allowing a degree of mixing between two different species in a layer of a certain positive thickness (i.e., the diffuse interface). The main advantage of this method with respect to its sharp counterpart lies in the fact that it prevents the presence of jumps, allowing a continuous transition when moving across different regions that are rich in a single phase. The diffuse interface approach revolves around an order parameter (usually a relative rescaled concentration difference) to keep track of the interface. Focusing on binary systems for the sake of simplicity, and letting $\varphi$ denote the order parameter, the domains that are rich in a single species are usually represented by the sets $\{\varphi = 1\}$ and $\{\varphi = -1\}$, respectively. Whenever the order parameter attains intermediate values, i.e., when $\varphi \in (-1,1)$, some mixing takes place: this region is exactly the diffuse interface.
    Of course, in order for such a mathematical description of binary fluids to be relevant, it must be physically consistent:
    it is then necessary to give the order parameter a precise physical meaning, ensuring e.g.~that
    $\varphi \in [-1,1]$.  
    
    Let $d \in \{2,3\}$ and let $\OO \subset \mathbb R^d$ be a bounded and sufficiently smooth domain. As the interaction between the two components is the result of the interplay between the mixing entropy and demixing effects leading to phase separation, the energy of the system is given by the functional
	\[
	\mathcal E(\varphi) = \int_\OO
    \left[
    \dfrac{\nup}{2}|\nabla \varphi|^2 + F(\varphi)
    \right]\,\d x,
	\]
	also usually called Helmholtz free energy. Here, the parameter $\nup$ is related to the interface thickness, while the function $F$ represents a double-well potential energy density. Thermodynamical consistency requires the derivative of $F$ to be singular at $\pm 1$: indeed, the relevant choice for $F$ is the so-called Flory--Huggins logarithmic potential (see \cite{Flory42, Huggins41}) given by
	\begin{equation}
		\label{F_log}
		F_{\text{log}}(r):=\frac\theta2\left[(1+r)\ln(1+r) + (1-r)\ln(1-r)\right] - \frac{\theta_0}2r^2, 
		\quad r\in(-1,1),
	\end{equation}
	where $0<\theta<\theta_0$ are fixed constants related to the absolute and to the critical temperature of the material under consideration. The evolution of $\varphi$ over time can be described by various systems of partial differential equations. The most prominent examples are the Cahn--Hilliard equation (i.e., the $H^{-1}$-gradient flow of the energy functional, see \cite{CahnHilliard58, CH1961}) and its second-order relaxation, the Allen--Cahn equation (i.e., the $L^2$-gradient flow of the energy functional, see \cite{AC1979}). Up to fixing suitable boundary and initial conditions, the former reads
	\[
	\begin{cases}
		\partial_t \varphi- \Delta w=0 &\quad \text{in }\OO \times (0,T),\\
		w = -\nup\Delta \varphi + F'(\varphi)&\quad \text{in }\OO \times (0,T),\\
	\end{cases}
	\]
	while the latter is
	\[
	\begin{cases}
		\partial_t \varphi + w = 0 &\quad \text{in }\OO \times (0,T),\\
		w = -\nup\Delta \varphi + F'(\varphi)&\quad \text{in }\OO \times (0,T),\\
	\end{cases}
	\]
	where in both systems $T > 0$ is some arbitrary time horizon, while $w$ is the variational derivative of the energy functional and denotes the chemical potential.
    Despite their original application related to materials science, in the last decades the model has proven to be remarkably flexible in describing numerous segregation-driven problems related to cell biology, see for instance \cite{Hym, Dol, Ber} and the more recent work \cite{Hes}. For example, 
    the modelling of liquid-liquid phase separation in cell biology \cite{Dol,Hes,A19,H14,M22,S17}, as well as the study of its respective equilibrium configurations, can help for studying the impact of radiation on biological tissue in the context of FLASH radiotherapy, a pioneering treatment modality poised to revolutionise the field of radiotherapy.
    Both models are now quite well understood from the mathematical standpoint, and even some of their variants have been analyzed (see, for instance, \cite{mir-CH, GGW20} and the references therein). 
    
    A refinement of these diffuse interface models for multicomponent fluids accounts for hydrodynamic effects introducing a coupling with a forced Navier--Stokes system. The source is called the Korteweg force, and models capillarity effects. Denoting by $\b u$ a suitably averaged velocity field for the fluid mixture, and by $\pi$ the corresponding pressure field, the refined versions of the models read
	\[
	\begin{cases}
		\partial_t \b u - \nuu \Delta \b u + (\b u \cdot \nabla) \b u + \nabla \pi = w \nabla \varphi &\quad \text{in }\OO \times (0,T),\\
		\div \b u = 0  &\quad \text{in }\OO \times (0,T),\\
		\partial_t \varphi + \b u \cdot \nabla \varphi - \Delta w =0 &\quad \text{in }\OO \times (0,T),\\
		w = -\nup\Delta \varphi + F'(\varphi)&\quad \text{in }\OO \times (0,T),\\
	\end{cases}
	\]
	when employing the (convective) Cahn--Hilliard equation, and
	\[
	\begin{cases}
		\partial_t \b u - \nuu \Delta \b u + (\b u \cdot \nabla) \b u + \nabla \pi = w \nabla \varphi &\quad \text{in }\OO \times (0,T),\\
		\div \b u = 0  &\quad \text{in }\OO \times (0,T),\\
		\partial_t \varphi + \b u \cdot \nabla \varphi + w= 0 &\quad \text{in }\OO \times (0,T),\\
		w = - \nup\Delta \varphi + F'(\varphi)&\quad \text{in }\OO \times (0,T),\\
	\end{cases}
	\]
	when employing the (convective) Allen--Cahn equation. In both systems, the constant $\nuu > 0$ denotes the kinematic viscosity. After the original classification provided in \cite{HH}, the mathematical literature on these systems has quite developed (we refer, for instance, to \cite{GPV, Blesgen,GG2010,GG2010-2} and to \cite{AGG2012,GKL2018,LT1998} for further refinements). 
    
    Although these classes of systems are widely studied under several aspects, deterministic diffuse interface models fail to render unpredictable oscillations happening at microscopic spatial scales, typically of thermal nature or due to configurational effects. This has lead to the investigation of their stochastic counterparts, starting from the pioneering contribution \cite{cook} proposing the first version of the model now known as stochastic Cahn--Hilliard equation. On the topic, we refer to \cite{daprato-deb,deb-zamb, deb-goud, EM1991,goud} and to the more recent works \cite{scar-SCH, scar-SVCH,scarpa21},
    and to \cite{HLR22, orr-scar, Bertacco21} on the stochastic Allen--Cahn equation. 
    
    Nonetheless, stochastic diffuse interface models coupled with random partial differential equations or even to other stochastic differential equations are far less understood with respect to their deterministic counterparts. In this regard, one example is given by the asymptotic analysis of the problem. On the deterministic side, several results have been established in terms of the existence of attractors or even the convergence to a single stationary state (see, for instance, \cite{AW07, mir-zel} or the more recent contribution \cite{GMT2019}), while the only result obtained in the stochastic case, in the direction of existence and uniqueness of invariant measures, is the very recent work on the Allen--Cahn equation with singular potential and degenerate noise \cite{SZ}. It is our aim, in this work, to perform the longtime behavior analysis on the stochastic Allen--Cahn--Navier--Stokes equation, that is,
	\begin{equation}
		\label{AC_NS}
		\begin{cases}
			{\rm d} \b{u}+ \left[ -\nuu \Delta \b{u}+(\b{u} \cdot \nabla )\b{u}+ \nabla \pi-w \nabla \varphi\right]\, {\rm d}t=G_1(\b{u})\, {\rm d}W_1 & \text{in} \ (0,T) \times \mathcal{O},
			\\
			\div \b{u}=0 & \text{in} \ (0,T) \times \mathcal{O},
			\\
			{\rm d}\varphi +\left[\bu \cdot \nabla \varphi+ w\right]\, {\rm d}t =G_2(\varphi)\, {\rm d}W_2 & \text{in} \ (0,T) \times \mathcal{O},
			\\
			w=-\nup\Delta \varphi +F'(\varphi) & \text{in} \ (0,T) \times \mathcal{O},
			\\
			\b{u}=0 & \text{in} \ (0,T) \times \partial\mathcal{O},
            \\
			\alpha_d \varphi + \alpha_n\nabla \varphi \cdot \b{n}=0  & \text{in} \ (0,T) \times \partial\mathcal{O},
			\\
			\b{u}(0)=\b{u}_0 & \text{in}\ \mathcal{O} \\ \varphi(0)=\varphi_0 & \text{in}\ \mathcal{O},   \end{cases}
	\end{equation}
	where $\OO \subset \mathbb{R}^2$ is a bounded sufficiently smooth domain, $W_1$ and $W_2$ are two independent cylindrical	Wiener processes on some (possibly different) separable Hilbert spaces,
	and $G_i$ is a suitable stochastically integrable process with respect to $W_i$, for $i \in \{1,2\}$. Here $\bn$ stands for the outward normal unit vector to $\partial\OO$ and the parameters $\nuu,\,\nup > 0$ are given. Concerning the boundary conditions for the order parameter $\varphi$, we assume that either $\alpha_d = 1$ and $\alpha_n = 0$, so that Dirichlet conditions are imposed, or $\alpha_n = 1$ and $\alpha_d = 0$, so that Neumann boundary conditions are imposed. In the latter case, the existence of martingale solutions in dimension two and three, as well as pathwise uniqueness and existence of probabilistically-strong solutions in dimension two for system \eqref{AC_NS} have been proven in \cite{DPGS}. In this follow-up work, we complete the analysis in the two-dimensional case by showing the existence, uniqueness and the asymptotic stability of invariant measures, taking suitable care of the assumptions on the noise in order to claim as much generality as possible.

    The mathematical literature on the long-time behaviour and ergodicity for stochastic systems is extremely developed. 
In general, showing the existence of an invariant measure is much easier than proving uniqueness. 
A classical method for proving the existence of an invariant measure is the Krylov–Bogoliubov method: see \cite{dapratozab} for a general presentation and 
\cite{BDP2006, MS20} for applications to problems in variational form.
A modification of this method has been proposed by Maslowski and Seidler in \cite{MS2}. This latter technique, that allows to work  with weak topologies, has been successful in the study of invariant measures in a series of recent works, see \cite{BMO} and \cite{BF19} for the stochastic Navier-Stokes equations, \cite{BOS} for the stochastic nonlinear beam and wave equations, \cite{BesFer} for the stochastic damped Euler equation, \cite{BLGR} for the stochastic inviscid multilayer quasi-geostrophic equation and \cite{BFZ} for the stochastic damped nonlinear Schr\"odinger equation.

Proving the uniqueness of the invariant measure and the convergence to it is in general a more challenging problem (see e.g. \cite{FerZannew} for an overview of the different techniques that can be used).
In \cite{GHMR17}, the authors identified a conceptually and intuitive framework for proving the uniqueness of the invariant measure by a generalized coupling technique. This approach has been developed in \cite{KS} where it is shown that by similar arguments one can also gets the asymptotic stability of the invariant measure. In \cite{GHMR17} and \cite{KS} many examples of PDEs driven by an additive noise are considered, for which this framework led to streamlined proof of uniqueness of the invariant measure. The main thread between these systems is the existence of a finite number of determining modes (low modes) and a sufficiently rich stochastic forcing term to ensure that the low modes are excited. This is usually referred to as ``effectively elliptic'' setting where all the presumptively unstable directions are stochastically forced. The key idea of the methods in \cite{GHMR17} and \cite{KS} is to introduce a suitable shift in the driving Wiener process to force solutions, that start at different initial conditions, together asymptotically as time goes to infinity. For strongly dissipative systems, in the spirit of \cite{FP}, it is usually enough to control a finite number of unstable directions by introducing a finite-dimensional shift and requiring a sufficiently rich stochastic forcing to ensure that the unstable modes are excited. These methods have been successfully used to prove ergodic properties of fluid dynamics PDEs driven by an additive noise, see the many examples in \cite{GHMR17} and \cite{KS}. These techniques have been recently extended in \cite{FZ23} to deal with multiplicative-type noises as well. 

Let us comment now on the mathematical challenges arising in the study of the system \eqref{AC_NS}.
The first main issue is that in our framework the system \eqref{AC_NS} can be solved only if the initial datum of the Allen-Cahn equation satisfies a nonlinear-type condition, namely $|\varphi_0|\leq1$ a.e.~in $\OO$.
This condition cannot be avoided as it is due to the singularity of the potential 
\eqref{F_log} and must then be taken into account in the setting up of the stochatsic flow. This forces to intend the transition semigroup as operators acting on functions defined on a general topological space, and not on a Hilbert space as in more classical settings. Such topological space would be the natural one associated to the stochastic flow \eqref{AC_NS}. 
The second challenge is due to a lack of satisfactory 
continuous dependence estimate (see \cite{DPGS}) with respect to the initial datum for the SPDE \eqref{AC_NS}.
This results in the impossibility of deducing Feller-type properties via standard techniques and 
naturally calls for a tailored proof based on ad-hoc arguments. The main idea is to first show that the transition semigroup is sequential weak Feller through stochastic compactness arguments, and then to infer the 
classical Feller property via variational energy methods.
Once the Feller and Markov properties of the semigroup are obtained, the proof of existence of invariant measures is presented in extreme generality, by exploiting either the sequential weak Feller property or the classical Feller property. 
In this direction, 
as far as the potential $F$ and the stochastic forcing terms are concerned, 
we work under the same assumptions that guarantee the existence and uniqueness of a probabilistically strong solution.

The paramount challenge and the main contribution of the entire work is related instead to proving uniqueness of the invariant measure and asymptotic stability.
Here, the crucial difficulty consists in handling the degeneracy of the noise in the Allen-Cahn equation (see Assumption \ref{hyp:diffusionAC} below), which is necessary to solve the SPDE \eqref{AC_NS}, and the singular coupling between the Allen-Cahn and the Navier-Stokes equations.
In the spirit of \cite{FZ23}, for the system \eqref{AC_NS} subjected to Dirichlet boundary conditions, we prove the uniqueness of the invariant measure and the asymptotic stability under the following two main requirements. The coefficient $\beta$ in the Allen-Cahn equation has to be sufficiently large and the noise driving the Navier-Stokes equations has to be non-degenerate on the unstable directions, i.e., the image of the covariance of the operator $G_1$ contains a finite number of low (unstable) modes.
We strongly emphasise that the non-degeneracy condition that we propose is only on the noise in the Navier-Stokes equation, as the operator $G_2$
in the noise driving the Allen-Cahn equation is indeed assumed to be degenerate. 
Such degeneracy forces to impose
the condition on $\beta$ to be sufficiently large, whereas the viscosity coefficient $\nu$ can be taken arbitrarily small. More precisely, what we manage to show is the following: provided that the parameter $\beta$ is large enough, for every arbitrarily small viscosity coefficient $\nu$ there is a minimum number $N$ of modes such that, if the noise operator $G_1$
driving the Navier-Stokes equation is non-degenerate on the first $N$ modes, then the invariant measure of the entire system is unique and asymptotically stable.
Actually, what we obtain is a joint condition on the parameter $\beta$ and the number $N$ of unstable directions forced by the noise. Roughly speaking, this means that,
provided that $\beta$ is sufficiently large, for every 
$\nu>0$ arbitrarily small there exists a number $ N$ of low modess
such that uniqueness holds.
 A joint condition on the parameters is not surprising since we are dealing with a coupled system of equations. Note that in the literature a condition corresponding to large dissipation is not new, see for example  
\cite{SZ} for the Allen-Cahn equation and \cite{BFZ23} for the nonlinear Schr\"odinger equation. 

The crucial technical tool to apply the techniques of \cite{GHMR17} and \cite{KS} to our problem is given by a Foias-Prodi-type estimate. The Foias-Prodi estimates describe the following property for an infinite dimensional dynamical system:
 given any two solutions, if they synchronise in the limit when time goes to infinity on a sufficient (but finite) number of components (i.e.~the low modes), then in fact all components synchronise. 
In our setting, the Foias-Prodi estimates describe the behavior for which any two solutions of system \eqref{AC_NS}, with different initial data, converge to each other as $t\to+\infty$ if a control acts on a sufficient finite number $N$ of modes in the Navier-Stokes system and if the coefficient $\beta$ in the Allen-Cahn equation is sufficiently large. 
We highlight that in the literature the Foias-Prodi estimates appear both in a pathwise formulation and  
in expectation: when the SPDE has a strong dissipation and an additive noise then the Foias-Prodi estimate can be proved pathwise,
whereas when there is a weak dissipation like a damping term (see \cite{DO}, \cite{GHMR21}) or when the noise is multiplicative (see \cite{Oda2008}, \cite{FZ23}), then the Foias-Prodi estimates can be proved in expected value.  

The Foias-Prodi estimates that we obtain are in expectation
and they pivot around a very subtle refined continuous-dependence result. Given any two solutions to the system, we show a
continuous-dependence estimate in stronger (energy) norms, where
the implicit constant appearing depends {\em only on one} of the solutions, and not on both. This delicate technical study 
is essential in order to possibly exploit the Foias-Prodi argument described above and is categorically non-trivial in general.
Besides the Navier-Stokes equation, where such type of 
stability estimate is known, for nonlinear systems it is false in general as the stability constant might depend on both solutions.
In our framework, we provide a refined stability estimate with such feature for the entire system, by exploiting the structural coupling between the two equations.
To our knowledge, this is the first contribution that presents a nontrivial example of such a stability estimate in the context of
a general system of SPDEs, and we expect it to have important developments beyond the Foias-Prodi argument. For example, 
we plan to use it in forthcoming works also to obtain quantitative mixing results.

\textbf{Main results.} We briefly 
list here the main contributions that we obtain in this work.
\begin{itemize}
    \item We prove refined regularity results for the system \eqref{AC_NS}, that allow to obtain continuous dependence estimates in energy-type norms with respect to the initial data,
    see Proposition~\ref{prop:uniqueness}.
    \item We prove that the transition semigroup associated to the system \eqref{AC_NS} is Markov, sequential weak Feller, and Feller in the classical sense, see Theorems~\ref{th:Feller} and \ref{th:Markov}. 
    As a consequence, we 
    show existence of ergodic invariant measures and characterise their support, see Theorems~\ref{th:ex_inv}, \ref{th:support}
    and \ref{th:ergodic}.
    \item We prove uniqueness of the invariant measure in the case of Dirichlet boundary conditions for the system \eqref{AC_NS}
    by proposing a joint condition involving the viscosity coefficients and the number of active low modes in the Navier-Stokes equation, see Theorem~\ref{asinvmeas}.
    \item We prove that the unique invariant measure is also asymptotically stable with respect to a suitable Wasserstein metric, see Theorem~\ref{weak_conv_ACNS}.
\end{itemize}

\textbf{Plan of the paper.} The contents of the paper are structured as follows. Section~\ref{sec:setting} contains all the necessary preliminary material. Section~\ref{sec:acns} illustrates known results on the well-posedness of the problem and establishes a new continuous dependence estimate in energy norms. Section~\ref{sec:existence} is devoted to investigating the existence of invariant measures, Section~\ref{sec:erg}
presents the characterisation of the support of the invariant measures and the ergodic result,
while the uniqueness and asymptotic stability issues are studied in Section~\ref{sec:uniqueness}. Finally, a collection of useful technical estimates are presented in Appendix~\ref{sec:appA} and Appendix~\ref{sec:appB}.


\section{Mathematical setting} \label{sec:setting}

\subsection{Notation}
For any Banach space $E$, we use the boldface symbol $\boldsymbol E$ for the space of two-dimensional vectors or 2-by-2 matrices with components in $E$. Moreover, we denote its topological dual by $E^*$. 
The duality pairing between $E$ and $E^*$ will be indicated by $\ip{\cdot}{\cdot}_{E^*,E}$. If no confusion arises, we may drop the subscripts.
For any real Hilbert space $H$, we denote by $\|\cdot\|_H$ and $(\cdot, \cdot)_H$ 
the norm and the scalar product, respectively. 
Given any two Banach spaces 
$E$ and $F$, we use the symbol $\mathcal{L}(E,F)$ for the space of linear bounded operators from $E$ to $F$. Furthermore, we write $E \hookrightarrow F$, if $E$ is continuously embedded in $F$.
If $H$ and $K$ are separable Hilbert spaces, we employ the symbol
$\LL_{HS}(H, K)$ for the space of Hilbert--Schmidt operators from $H$ to $K$, endowed with its canonical norm $\norm{\cdot}_{\LL_{HS}(H, K)}$. 
For a fixed $T>0$ and $\alpha\in(0,1]$, we denote by 
$\C ([0,T];E)$, $\C_\text{w}([0,T]; E)$,
and $\C^{\alpha}([0,T];E)$ the spaces of
functions from $[0,T]$ to $E$ which are 
strongly continuous, weakly continuous, 
and $\alpha$-H\"older-continuous, respectively.
Moreover, for every $p\in[1,+\infty]$
we use the classical symbol $L^p(0,T; E)$ for the space
of strongly measurable
$p$-Bochner integrable functions from $(0,T)$ to $E$.
Analogously, for every $\gamma>0$, we denote by
$W^{\gamma, p}(0,T;E)$ and $W^{-\gamma,p}(0,T;E)$
the usual $E$-valued Sobolev spaces. As for the probabilistic framework, let $(\Omega,\cF,\mathbb{F}:=(\cF_t)_{t\in[0,T]},\P)$ be a filtered probability space satisfying the usual conditions (namely it is saturated and right-continuous), with $T>0$ being a prescribed final time. Given a measurable space $(E,\mathscr{E})$, the law of a random variable $\xi: \Omega \to E$ will be denoted by $\text{Law}_{\mathbb{P}}(\xi)$. We denote by $\cP$ the progressive sigma algebra on $\Omega\times[0,T]$.
For every $s\in[1,+\infty]$ and for every Banach space $E$
the symbol $L^s(\Omega; E)$ 
indicates the usual space of strongly measurable $s$-Bochner-integrable functions from $\Omega$ to $E$.
Whenever $E$ is itself a Bochner space, for all $s,r\in[1,+\infty)$ we write
$L^s_\cP(\Omega;L^r(0,T; E))$ to stress that measurability is intended with respect to $\cP$. For all $s\in(1,+\infty)$ and for every separable and reflexive Banach space $E$ we also define
	\[
	L^s_w(\Omega; L^\infty(0,T; E^*)):=
	\left\{v:\Omega\to L^\infty(0,T; E^*) \text{ weakly* measurable}\,:\,
	\norm{v}_{L^\infty(0,T; E^*)}\in L^s(\Omega)
	\right\}
	\]
 so that one has the canonical identification
 $L^s_w(\Omega; L^\infty(0,T; E^*))=
 (L^{s'}_\cP(\Omega; L^1(0,T; E)))^*$ with $s':=\frac{s}{s-1}$.
 
 Throughout the whole work, we reserve the symbol $K$ (eventually indexed) for universal constants possibly only depending on the domain $\OO$. Analogously, we reserve the symbol $C$ (eventually indexed) for constants depending on some or all the structural parameters of the problem. If the dependencies of such constants are relevant, they will be explicitly pointed out. The values of these constants may change within the same argument without relabelling.

\subsection{Functional spaces and operators}
	Consider a bounded domain $\OO\subset\erre^2$ with smooth boundary $\partial\OO$, outward
	normal unit vector $\bn$ and Lebesgue measure denoted by $|\OO|$.
    The parameters $\alpha_d,\alpha_n\in\{0,1\}$ are fixed and will identify the boundary condition of the Allen-Cahn equation, 
    namely it holds that either $(\alpha_d,\alpha_n)=(1,0)$ (Dirichlet) or $(\alpha_d,\alpha_n)=(0,1)$ (Neumann).
	The symbol $W^{s,p}(\OO)$, where $s\in\erre$ and $p\in[1,+\infty]$, denotes the usual real Sobolev space of order $(s,p)$ and we denote by $\norm{\cdot}_{W^{s,p}(\OO)}$ its canonical norm. In the Hilbert case $p = 2$, we define $H^s(\OO):=W^{s,2}(\OO)$, $s\in\erre$,
	endowed with its canonical norm $\norm{\cdot}_{H^s(\OO)}$.
	As customary, we introduce the zero-trace space
    \[
    H^1_0(\OO) := \{v \in H^1(\OO) : v = 0 \text{ a.e.~on } \partial\OO \},
    \]
    endowed with the norm
    \[
    \|u\|_{H^1_0(\OO)} := \|\nabla u\|_{\b L^2(\OO)} \qquad \forall \: u \in H^1_0(\OO),
    \]
    owing to the Poincaré inequality. According to the type of boundary conditions, for the sake of convenience, we let also
	\[
	H:=L^2(\OO), \qquad V_1:=\begin{cases}
	    H^1(\OO) & \quad \text{if }\alpha_d = 0 \text{ and } \alpha_n = 1, \\
        H^1_0(\OO) & \quad \text{if }\alpha_d = 1 \text{ and } \alpha_n = 0,
	\end{cases}
    \]
    and
    \[
    V_2:=\begin{cases}
	    \left\{\psi\in H^2(\OO):\;\partial_\bn\psi=0\text{ a.e.~on } \partial\OO \right\} & \quad \text{if }\alpha_d = 0 \text{ and } \alpha_n = 1, \\
        H^2(\OO) \cap H^1_0(\OO) & \quad \text{if }\alpha_d = 1 \text{ and } \alpha_n = 0,
	\end{cases}
	\]
	endowed with their standard norms $\norm{\cdot}_H$,
	$\norm{\cdot}_{V_1}$, and $\norm{\cdot}_{V_2}$, respectively.
	As usual, we identify the pivot Hilbert space $H$ with its dual through
	the Riesz isomorphism, so that we have the variational structure
	\[
	V_2\embed V_1\embed H \embed V_1^* \embed V_2^*\,,
	\]
	with dense and compact embeddings. We will denote by $A:V_1\to V_1^*$ the variational realization of the negative Laplace operator with homogeneous Neumann or Dirichlet boundary conditions, namely
	\[
	\ip{A\psi}{\phi}_{V_1^*,V_1}=\int_\OO\nabla\psi\cdot\nabla\phi\,,
	\qquad\forall \: \psi,\phi\in V_1.	\]
    Let us recall that, in the case of Dirichlet boundary conditions, the operator $A$ is the Riesz isomorphism between $V_1$ and its topological dual $V_1^*$. 
    
    In order to handle the Navier--Stokes velocity field, we define the following solenoidal vector-valued spaces
	\begin{align*}
		\bHs := \overline{\{ \bv \in  \b{\C}^\infty_0(\OO): \div \bv = 0 \text{ in } \OO\}}^{\b{L}^2(\OO)}, \quad
		\bVs := \overline{\{ \bv \in \b{\C}^\infty_0(\OO) : \div \bv = 0 \text{ in } \OO\}}^{\b{H}^1(\OO)}.
	\end{align*}
    The space $\bHs$ is endowed with the Hilbert structure inherited by $\b H$, and we denote its norm and scalar product by $\norm{\cdot}_{\bHs}$ and $(\cdot,\cdot)_{\bHs}$, respectively.
	By means of the Poincaré inequality, on the space $\bVs$ we can use the norm $\norm{\bv}_{\bVs}:=\norm{\nabla\bv}_{\b L^2(\OO)}$, for all $\bv\in\bVs$, and its respective
	scalar product $(\cdot,\cdot)_{\bVs}$. By identifying the Hilbert space $\bHs$ with its dual we have the variational  structure
	\[
	\bVs \embed \bHs \embed  \bVsd,
	\]
	 with dense and compact embeddings. The Stokes operator $\b A:\bVs\to\b V_\sigma^*$ is defined as
	the canonical Riesz isomorphism of $\bVs$, i.e.,
	\[
	\ip{\b A\bv}{\b w}_{\b V^*_\sigma, \b V_\sigma}:=
	(\nabla\bv, \nabla\b w)_{\b H}\,, \qquad \forall \: \bv,\b w\in\bVs.
	\]
	Recalling the spectral properties of the operator $\b{A}$, it is possible to define the family of power operators $\b{A}^s$ for any $s \in \mathbb{R}$. In particular, using $\bHs$ as pivot space, we highlight the general structure
		\[
		D(\b{A}^s) \embed D(\b{A}^t) \embed \bHs \equiv D(\b{A}^0) \embed  D(\b{A}^{-t}) \embed D(\b{A}^{-s})
		\]
		for any $s > t > 0$, with dense and compact embeddings (see also \cite{DPGS}).
We define the usual Stokes trilinear form
	$b$ on $\bVs \times \bVs \times \bVs$
	\[
	b(\bu, \b v, \b w):=\int_{\OO}(\b u\cdot\nabla)\b v\cdot\b w
	=\sum_{i,j=1}^2\int_\OO u_i\frac{\partial v_i}{\partial x_j}w_j\,,
	\qquad \forall \: \b u, \b v, \b w \in \bVs\,,
	\]
	and the associated bilinear form $\b B:\bVs\times\bVs\to \b V_\sigma^*$ as
	\[
	\ip{\b B(\bu, \bv)}{\b w}_{\b V_\sigma^*, \bVs}:=b(\bu, \b v, \b w)\,,
	\qquad\forall \: \b u, \b v, \b w \in \bVs\,.
	\]
	Let us recall that $b(\bu, \bv, \b w)=-b(\bu, \b w, \bv)$ for all
	$\b u, \b v, \b w \in \bVs$, from which it follows in particular
	that $b(\bu, \bv, \bv)=0$ for all $\b u, \b v\in \bVs$. 
    
    Eventually, let us deal with the stochastic terms appearing in equation \eqref{AC_NS}. On the normal filtered probability space $(\Omega,\cF,(\cF_t)_{t\in[0,T]},\P)$, we consider two independent
	cylindrical Wiener processes $W_1$ and $W_2$ with values in some separable Hilbert spaces $U_1$ and $U_2$,
	respectively. For convenience, we fix once and for all two complete orthonormal systems
	$\{u^1_j\}_{j\in\enne}$ on $U_1$ and $\{u^2_j\}_{j\in\enne}$ on $U_2$.
    We recall that, as a cylindrical process on $U_i$, $i\in \{1,2\}$, $W_i$ admits the following formal representation
	\begin{equation} \label{eq:representation}
		W_i = \sum_{k\in\enne} \beta_k u^i_k,
	\end{equation}
	where $\{\beta_k\}_{k \in \enne}$ is a family of real and independent Brownian motions. It is well known that \eqref{eq:representation} does not converge in $U_i$, in general. Nonetheless, there always exists a larger Hilbert space $U_0^i$, such that
	$U_i \embed U_0^i$ with Hilbert--Schmidt embedding $\iota_i$, enabling the identification of $W_i$ as a $Q^0_i$-Wiener process on $U_0^i$, with $Q^0_i=\iota_i \circ \iota_i^*$ being trace-class on $U_0^i$ (see \cite[Subsections 2.5.1]{LiuRo}). In the following, we may implicitly assume this extension by simply saying that $W_i$ is a cylindrical process on $U_i$. This holds also for stochastic integration with respect to $W_i$: 
  for every real Hilbert space $K$ and
  for every process $G \in L^2_\cP(\Omega; L^2(0,T;\cL^2(U_i,K)))$,
  the stochastic integral
	\[
	\int_0^\cdot G(s)\,\d W_i(s)
	\]
is well defined in terms of the usual stochastic integration with respect to $W_i$ as a $Q_0^i$-Wiener process. The definition is well posed and does not depend on the choice of $U_i^0$ or $\iota_i$ (see \cite[Subsection 2.5.2]{LiuRo}).

\subsection{Embeddings and interpolation inequalities}
Throughout this work, we shall make frequent use of two-dimensional embeddings and interpolation theorems. For the sake of clarity, we shall list several of them hereafter.
\begin{enumerate}[(i)] \itemsep0.3em
    \item The Ladyzhenskaya inequality
    \[
    \|u\|_{L^4(\OO)} \leq \KL\|u\|_H^\frac 12\|u\|_{V_1}^\frac 12
    \]
    holds for all $u \in V_1$.
    \item The Gagliardo--Nirenberg inequality
    \[
    \|u\|_{W^{1,4}(\OO)} \leq \KGN\|u\|_{L^\infty(\OO)}^\frac 12\|u\|_{H^2(\OO)}^\frac 12
    \]
    holds for all $u \in H^2(\OO)$.
    \item The two-dimensional embedding $V_1 \embed L^q(\OO)$ for all finite $q \geq 1$ implies that
    \[
    \|u\|_{L^q(\OO)} \leq K_q \|u\|_{V_1}
    \]
    for all $u \in L^q(\OO)$.
    \item The elliptic regularity theory for the Laplace problem
    with either Neumann or Dirichlet conditions implies that
    \[
    \|u\|_{H^2(\OO)} 
    \leq K_\Delta\left( \|\Delta u\|_{H} + \alpha_n\|u\|_H \right)
    \]
    for all $u \in V_2$.
\end{enumerate}

\subsection{Weak topologies and weak continuity} 
\label{ssec:weaktop}
Let us recall some classical results on weak topologies that will be useful later on. Of course, what follows is by no means exhaustive, and, for further detail, we refer the reader to \cite{BF19,MS} and the references therein.
For a separable Hilbert space $E$, we denote by $\tau_\text{s}^E$, $\tau_\text{w}^E$, 
and $\tau_{\text{bw}}^E$ the strong, weak, and bounded-weak topology, respectively. We recall that the bounded-weak topology $\tau_{\text{bw}}^{E}$ is defined as the finest topology on $E$ that coincides with the weak topology of $E$ on every norm-bounded subset of $E$, i.e., a set $A \subset E$ is closed with respect to the bounded-weak topology if and only if $A \cap U$ is closed with respect to the weak topology in $U$ for any bounded set $U$ (equivalently, for any ball).
We also recall that since $E$ is separable, then the strong,
weak, and bounded-weak Borel $\sigma$-algebras coincide.
Hence, we shall simply speak of the Borel $\sigma$-algebra with no ambiguity and it will be denoted by $\mathscr B(E)$. Analogously, 
we will use the symbol $\mathcal B(E)$ 
to denote the space of 
real-valued Borel-measurable functions on $E$, and the symbol $\mathcal B_b(E)$
to denote the space of real-valued Borel-measurable functions on $E$ which are also bounded.
We also recall that a function $f:E\to\mathbb R$ is continuous with respect to the bounded-weak topology $\tau_{\text{bw}}^{E}$ if and only if $f$ is sequentially continuous with respect to the weak topology $\tau_\text{w}^E$: this is a consequence of the fact that the weak topology on every closed ball of a separable Hilbert space $E$ is metrizable. Given any topological space $(E, \tau)$, we denote by $\C(E)$, $\C_b(E)$, $\mathcal S \mathcal C(E)$ and $\S \C_b(E)$ the spaces of continuous, bounded continuous, sequentially continuous, and bounded sequentially real-valued continuous functions on the topological space $(E, \tau)$, respectively. We recall that $f: E \rightarrow \mathbb{R}$ is called sequentially continuous if for every sequence $\{e_n\}_{n \in \enne}\subset E$
such that $e_n \to e \in E$ with respect to its topology, we have that $f(e_n)$ converges to $f(e)$ in $\erre$.

\subsection{Spaces of trajectories}
Let $\delta \in \left(0, \frac 12\right)$ be fixed arbitrarily.
For every $T>0$,
we define the locally convex spaces 
\begin{align}
\label{Z_phi}
\mathcal Z_T^{\varphi} &:= \C([0,T];H) \cap L^2(0,T;V_1) \cap 
\C_\text{w}([0,T];V_1),\\
\label{Z_u}
\mathcal Z_T^{\bu} &:=\C([0,T];D(\b{A}^{-\delta}))\cap L^2(0,T;\bHs) \cap \C_\text{w}([0,T];\bHs),
\end{align}
endowed with the topologies given by the supremum of the corresponding topologies in the right-hand side.
Let us point out that the topologies on the spaces of weakly continuous functions are generated by the family of seminorms
\[\left\{\varphi\mapsto\sup_{t\in[0,T]}|\langle\psi,\varphi(t)\rangle_{V_1^*, V_1}|
\right\}_{\psi\in V_1^*}\]
for $\C_\text{w}([0,T];V_1)$, and 
\[\left\{\bu\mapsto\sup_{t\in[0,T]}|\langle \b v,\bu(t)\rangle_{\bVsd, \bVs}|
\right\}_{\b v\in \b V_\sigma^*},\]
for $\C_\text{w}([0,T];\bHs)$.
When working on the unbounded time interval $[0,+\infty)$, we consider 
the locally convex topological spaces
\begin{align}
\label{Z_inf_phi}
\mathcal Z_\infty^{\varphi} &:= \C([0,+\infty);H) \cap L^2_{\text{loc}}(0,+\infty;V_1) \cap 
\C_{\text{w}}([0,+\infty);V_1)\,,\\
\label{Z_inf_u}
\mathcal Z_\infty^{\bu} &:=\C([0,+\infty);D(\b{A}^{-\delta}))\cap L^2_{\text{loc}}(0,+\infty;\bHs) \cap
\C_{\text{w}}([0,+\infty);\bHs)\,,
\end{align}
endowed again with the topologies
given by the supremum of the corresponding topologies in the right-hand side.
Recall that on an unbounded time interval, if $E$ is a separable Hilbert space, then
\begin{enumerate}[(i)]
 \item for the space of strongly continuous functions $\C([0,+\infty);E)$, we consider the topology induced by the metric
\[(v,w)\mapsto\displaystyle\sum_{k\in\enne}
     \frac
     1{2^k}\frac{\|v-w\|_{\C([0,k];E)}}{1+\|v-w\|_{\C([0,k];E)}}, \qquad \forall \: v,w \in \C([0,+\infty);E);\]
\item for the space of locally square integrable functions
$L^2_{\text{loc}}(0,+\infty;E)$, we consider the topology induced by the metric
\[(v,w)\mapsto\displaystyle\sum_{k\in\enne}
     \frac 1{2^k}\frac{\|v-w\|_{L^2(0,k;E)}}{1+\|v-w\|_{L^2(0,k;E)}}, \qquad \forall \: v,w \in L^2_{\text{loc}}([0,+\infty);E);\]    
\item for the space of weakly continuous functions
$\C_{\text{w}}([0,+\infty);E)$, we consider the topology generated by the family of seminorms
\[ 
\left\{w\mapsto\sup_{0\le t\le k}|\langle v, w(t)\rangle_{E^*,E} |\right\}_{k \in \mathbb N,\,v \in E^*}.\]
\end{enumerate}
Let now $\Gamma$ be a family of probability measures defined on the $\sigma$-algebra of Borel subsets of $\mathcal Z_T^{\varphi}$ ($\mathcal Z_T^{\bu}$, respectively). We recall that the family $\Gamma$ is said to be tight if for any $\eta >0$ there exists a compact set $K_\eta \subset \mathcal Z_T^{\varphi}$ ($\mathcal Z_T^{\bu}$, respectively) such that 
\begin{equation*}
\inf_{\nu \in \Gamma}\nu(K_\eta) \ge 1- \eta.
\end{equation*}
The tightness in $\mathcal Z_\infty^{\varphi}$ and $\mathcal Z_\infty^{\bu}$
is equivalent to the tightness
of the traces of the measures in $\Gamma$ in $\mathcal Z_N^{\varphi}$ and $\mathcal Z_N^{\bu}$,
for all $N \in\mathbb N$, respectively, see e.g., \cite{fland-gat}.

\subsection{Assumptions}
Let us state the set of assumptions that will be used throughout the paper. We work in a similar framework with respect to the one considered in \cite{DPGS}.
\begin{enumerate}[label=\textbf{(A\arabic*)}] \itemsep0.3em
\item \label{hyp:structural} The spatial domain $\OO$ is a bounded smooth domain in $\mathbb R^2$. The parameters $\nuu$ and $\nup$ are strictly positive, and we set either $\alpha_d = 0$ and $\alpha_n = 1$ (Neumann boundary conditions for the order parameter) or $\alpha_n = 0$ and $\alpha_d = 1$ (Dirichlet boundary conditions for the order parameter).
\item \label{hyp:potential} The potential $F:[-1,1]\to[0,+\infty)$ enjoys the following properties:
  \begin{enumerate}[(i)]
  \item \label{hyp:regularity} it satisfies $F\in C^0([-1,1])\cap C^2(-1,1)$ and $F'(0)=0$,
  \item \label{hyp:convex} there exists a constant $L_F>0$ such that $F''(r)\geq -L_F$ for all $r\in(-1,1)$,
  \item \label{hyp:singular} it holds that 
  \[
  \lim_{r\to-1^{+}}F'(r)=-\infty\,, 
  \qquad
  \lim_{r\to1^{-}}F'(r)=+\infty\,.
  \]
  \end{enumerate}
  As customary, the function $F$ is extended with value $+\infty$ outside $[-1,1]$ without relabeling, so that it can be regarded as a proper lower semicontinuous function $F:\erre\to[0,+\infty]$. 
\item \label{hyp:diffusionNS}
The operator
$G_1: \bHs \to \LL_{HS}(U_1, \bHs)$ is bounded in $\bHs$,
i.e. there exists $L_{G_1} > 0$ such that
\[
\|G_1(\bv)\|^2_{\LL_{HS}(U_1,\bHs)} \leq L_{G_1}
\quad\forall\,\bv \in \bHs.
\]
Moreover, we assume that $G_1: Z \to \LL_{HS}(U_1, Z)$ is $L_{G_1}$-Lipschitz-continuous where  $Z$ is either $\bHs$ or $\bVsd$.
\item \label{hyp:diffusionAC} Setting $\BB$ to be the closed unit ball in $L^\infty(\OO)$, the operator $G_2:\BB \rightarrow \mathcal{L}_{HS}(U_2,H)$ satisfies  
\begin{equation*}
	G_2(\psi)[u_k^2]=g_k(\psi), \quad \psi \in \BB, \quad k \in \mathbb{N},
\end{equation*}
for some sequence $\{g_k\}_{k \in \mathbb{N}} \subset W^{1, \infty}(-1,1)$ satisfying
\[
g_k(\pm 1)=0, \qquad F^{\prime \prime}g_k^2,\,F^\prime g_k \in L^\infty(-1,1)
\] 
for every $k \in \mathbb{N}$ such that
\begin{equation*}
\label{C_B}
L_{G_2} := \sum_{k \in \mathbb{N}}
\left(\|g_k\|^2_{W^{1, \infty}(-1,1)} +
\norm{F''g_k^2}_{L^\infty(-1,1)}
+\norm{F'g_k}^2_{L^\infty(-1,1)}\right) < +\infty.
\end{equation*}
\end{enumerate}

\begin{remark} 
Let us briefly comment the above assumptions. First, we point out once and for all that, whenever Assumption \ref{hyp:structural} holds, we are considering both Neumann or Dirichlet boundary conditions for the order parameter. Whenever only one of the two choices is admissible, it will be explicitly pointed out. As for Assumption \ref{hyp:potential}, one can easily check that the logarithmic potential \eqref{F_log} (up to some additive constant) satisfies all the properties listed therein. Moreover, the validity of Assumption \ref{hyp:potential} ensures that, for instance,
\begin{equation}
\label{F'_prop}
F'(r)r\ge L_Fr^2-2L_F \qquad \forall \: r \in (-1,1),
\end{equation} 
and that $F$ can be decomposed as the sum of a proper convex lower semicontinuous function $F_1:\mathbb [-1,1]\to[0,+\infty]$ and a regular function $F_2:\erre \to \erre$ such that $F_2\in C^2(\mathbb R) \cap C^{1,1}(\mathbb R)$.
Finally, observe that Assumption \ref{hyp:diffusionAC} implies that $G_2$ is $\sqrt{L_{G_2}}$-Lipschitz-continuous with respect to the $H$-metric on $\BB$. Indeed, since
	\[
	G_2(\psi)[h]:= \sum_{k \in \mathbb{N}}(h,u_k^2)_{U_2}g_k(\psi),
	\quad \psi \in \BB, \quad h \in U_2,
	\]
	by a direct computation (see also e.g.~\cite[Section 2]{Bertacco21}) it follows, for all $\psi,\phi \in \BB$, that 
	\begin{align}
		\label{HS_norm}
		\|G_2(\psi)\|_{\mathcal{L}_{HS}(U_2,H)}^2 &\le L_{G_2}|\OO|,\\
		\|G_2(\psi)-G_2(\phi)\|_{\mathcal{L}_{HS}(U,H)}^2 &\le L_{G_2}\|\psi-\phi\|^2_H.
	\end{align}
\end{remark}


\section{Results on the Allen--Cahn--Navier--Stokes system} \label{sec:acns}
\subsection{Existence and uniqueness of solutions}
Existence and uniqueness of a probabilistically-strong solution to the Allen--Cahn--Navier--Stokes problem \eqref{AC_NS} in the two-dimensional case is proved in \cite[Theorems 2.7 and 2.9]{DPGS} in the case $\nup = \nuu = 1$ and Neumann boundary conditions, i.e., when $\alpha_d = 0$ and $\alpha_n = 1$. However, the general case of $\nuu,\,\nup > 0$, as well as the Dirichlet case, follow from the same arguments with very minor modifications that will be omitted.
Let us recall here the concept of strong and martingale  solution and the main well-posedness result.
\begin{defin}
		\label{def:strong_sol}
		Assume that Assumption \ref{hyp:structural}--\ref{hyp:diffusionAC} hold and let $T > 0$ be fixed. Let further $p \geq 2$ and let $(\bu_0, \varphi_0)$ satisfy
		\begin{align}
			\label{eq:u0}
			\bu_0&\in L^p(\Omega,\cF_0; \bHs),\\
			\label{eq:phi0}
			\varphi_0&\in L^p(\Omega,\cF_0; V_1), \quad
			F(\varphi_0)\in L^{\frac p2}(\Omega,\cF_0; L^1(\OO)).
		\end{align}
		A probabilistically-strong solution to problem \eqref{AC_NS}
		with respect to the initial datum $(\bu_0, \varphi_0)$
		is a pair of processes $(\bu, \varphi)$ such that,
		\begin{align}
			\label{u}
			&\bu\in L^p_\cP(\Omega; L^\infty(0,T; \bHs)) \cap
			L^p_\cP(\Omega; L^2(0,T; \bVs)),\\
			\label{phi}
			&\varphi \in L^p_\cP(\Omega; C^0([0,T]; H))\cap
			L^p_w(\Omega; L^\infty(0,T; V_1)) \cap
			L^p_\cP(\Omega; L^2(0,T; V_2)),\\
			& |{\varphi}| < 1 
            \quad\text{a.e. in } \Omega \times \OO \times [0,T], \\
			\label{mu}
			&w:=-\nup\Delta\varphi+F'(\varphi) \in L^{p}_\cP(\Omega; L^2(0,T; H)),\\
			\label{initial}
			&(\bu(0), \varphi(0))=(\bu_0, \varphi_0),
		\end{align}
		and
		\begin{align}
			\nonumber
			&(\bu(t),\bv)_{\bHs} +
			\int_0^t\left[\nuu\ip{\b A \bu(s)}{\bv}_{\b V_\sigma^*, \b V_\sigma}
			+\ip{\b B( \bu(s),  \bu(s))}{\bv}_{\b V_\sigma^*, \b V_\sigma}
			-\int_\OO w(s)\nabla \varphi(s)\cdot\bv
			\right]\,\d s\\
			\label{var1}
			&\qquad= (\bu_0,\bv)_{\bHs} +
			\left(\int_0^t G_1( \bu(s))\,\d  W_1(s), \bv\right)_{\bHs}
			\qquad\forall\,\bv \in \bVs,
            \quad\forall\,t\geq0,
            \quad\P\text{-a.s.},\\
			\nonumber
			&( \varphi(t),\psi)_H +
			\int_0^t\!\int_\OO\left[ \bu(s)\cdot\nabla \varphi(s) +  w(s)\right]\psi\,\d s\\
			\label{var2}
			&\qquad= ( \varphi_0,\psi)_{H} +
			\left(\int_0^t G_2( \varphi(s))\,\d  W_2(s), \psi\right)_{H}
			\qquad\forall\,\psi\in V_1,
            \quad\forall\,t\geq0,
            \quad\P\text{-a.s.}
		\end{align} 
	\end{defin}
\begin{remark}
	Observe that the paths of $\bu$ lie also in the space of weakly continuous functions $\C_{\text w}([0,T]; \bHs)$ or even in $\C([0,T]; D(\b A^{-\delta}))$ for a sufficiently small $\delta > 0$.
\end{remark} 

\begin{defin}
		\label{def:mart_sol}
		Assume that Assumptions \ref{hyp:structural}--\ref{hyp:diffusionAC} hold and let $T > 0$ be fixed. Let further $p \geq 2$ and let $(\bu_0, \varphi_0)$ satisfy
        \eqref{eq:u0}--\eqref{eq:phi0}.
		A martingale solution to problem \eqref{AC_NS} with respect to the initial datum $(\bu_0, \varphi_0)$
		is a family
		\[
		\left(\left(\widehat\Omega, \widehat\cF, (\widehat\cF_t)_{t\in[0,T]}, \widehat\P\right),
		 \widehat W_1, \widehat W_2,
		\widehat\bu, \widehat\varphi\right),
		\]
		where $(\widehat\Omega, \widehat\cF, (\widehat\cF_t)_{t\in[0,T]}, \widehat\P)$ is a filtered
		probability space satisfying the usual conditions; $\widehat W_1,\, \widehat W_2$ are two independent cylindrical Wiener processes on $U_1$ and $U_2$, respectively; the pair of processes $(\widehat \bu, \widehat\varphi)$ satisfies,
        by setting $(\widehat\bu_0,\widehat\varphi_0):=
        (\widehat\bu(0), \widehat\varphi(0))$, that
		\begin{align}
			\label{eq:u_hat}
			&\widehat\bu\in L^p_w(\widehat\Omega; L^\infty(0,T; \bHs)) \cap
			L^p_\cP(\widehat\Omega; L^2(0,T; \bVs)),\\
			\label{eq:phi_hat}
			&\widehat\varphi \in L^p_\cP(\widehat\Omega; C^0([0,T]; H))\cap
			L^p_w(\widehat\Omega; L^\infty(0,T; V_1)) \cap
			L^p_\cP(\widehat\Omega; L^2(0,T; V_2)),\\
			& |\widehat{\varphi}| < 1 \text{ a.e. in } \hom \times \OO \times [0,T], \\
			\label{eq:mu_hat}
			&\hw:=-\nup\Delta\widehat\varphi+F'(\widehat\varphi) \in
			L^{p}_\cP(\widehat\Omega; L^2(0,T; H)),\\
			\label{eq:initial_hat}
			&\text{Law}_{\hP}(\widehat \bu_0, \widehat \varphi_0)=\text{Law}_{\P}(\bu_0, \varphi_0)
			\text{ on } \bHs\times V_1;
		\end{align}
		and
		\begin{align}
			\nonumber
			&(\widehat \bu(t),\bv)_{\bHs} +
			\int_0^t\left[\nuu\ip{\b A\widehat \bu(s)}{\bv}_{\b V_\sigma^*, \b V_\sigma}
			+\ip{\b B(\widehat \bu(s), \widehat \bu(s))}{\bv}_{\b V_\sigma^*, \b V_\sigma}
			-\int_\OO \hw(s)\nabla\widehat \varphi(s)\cdot\bv
			\right]\,\d s\\
			\label{eq:var1_hat}
			&\qquad= (\widehat \bu_0,\bv)_{\bHs} +
			\left(\int_0^t G_1(\widehat \bu(s))\,\d \widehat W_1(s), \bv\right)_{\bHs}
			\qquad\forall\,\bv \in \bVs,
            \quad\forall\,t\geq0,
            \quad\P\text{-a.s.},\\
			\nonumber
			&(\widehat \varphi(t),\psi)_H +
			\int_0^t\!\int_\OO\left[\widehat \bu(s)\cdot\nabla\widehat \varphi(s) + \hw(s)\right]\psi\,\d s\\
			\label{eq:var2_hat}
			&\qquad= (\widehat \varphi_0,\psi)_{H} +
			\left(\int_0^t G_2(\widehat \varphi(s))\,\d \widehat W_2(s), \psi\right)_{H}
			\qquad\forall\,v \in V_1,
            \quad\forall\,t\geq0,
            \quad\P\text{-a.s.}.
		\end{align}
	\end{defin} \noindent
\begin{remark}
    We point out, once and for all, that the existence of martingale solutions to \eqref{AC_NS} is shown in \cite{DPGS} in the case $\nup = \nuu = 1$ and with Neumann boundary conditions. Once again, we stress that the general case where $\nup,\,\nuu > 0$ and/or with Dirichlet boundary conditions can be obtained analogously.
\end{remark}

\noindent
The strong well-posedness result obtained in \cite{DPGS}
requires the choice $Z=\bVsd$ and reads as follows.
\begin{thm}
		\label{th:probstrongsol}
		Let Assumptions \ref{hyp:structural}-\ref{hyp:diffusionAC} hold with $Z = \bVsd$ and let $T > 0$ be fixed. Then,
		for every initial datum $(\bu_0,\varphi_0)$ satisfying \eqref{eq:u0}-\eqref{eq:phi0} with $p > 2$, there exists a unique probabilistically-strong solution
		$(\bu, \varphi)$ for problem \eqref{AC_NS} in the sense of Definition~\ref{def:strong_sol}.
  Moreover, for every pair of initial data 
  $(\bu_0^1,\varphi_0^1)$ and $(\bu_0^2,\varphi_0^2)$ satisfying \eqref{eq:u0}-\eqref{eq:phi0},
  there exist a sequence $\{C_n\}_{n \in \mathbb N}$ of positive real numbers depending on the structural parameters of the problem and a sequence $\{\tau_n\}_{n \in \mathbb N}$ of positive stopping times such that 
  $\tau_n\nearrow T$ $\P$-almost surely as $n \to \infty$ and 
  the respective solutions $(\bu_1, \varphi_1)$ and $(\bu_2,\varphi_2)$ to problem \eqref{AC_NS} satisfy, for all $n \in \mathbb N$, the stopped continuous dependence estimate
  \begin{multline*}
      \|(\bu_1-\bu_2)^{\tau_n}\|_{L^p_\cP(\Omega; C^0([0,T]; \b V_\sigma^*)
      \cap  L^2(0,T;\b H_\sigma))}+
      \|(\varphi_1-\varphi_2)^{\tau_n}\|_{L^p_\cP(\Omega; C^0([0,T]; H)
      \cap L^2(0,T; V_1))}\\
      \leq C_n\left[
      \|\bu_0^1-\bu_0^2\|_{L^p(\Omega; \b V_\sigma^*)} + 
      \|\varphi_0^1-\varphi_0^2\|_{L^p(\Omega; H)}
      \right].
  \end{multline*}
\end{thm}
\begin{remark}
    Let us point out that the analysis of the system \eqref{AC_NS} in \cite{DPGS} is carried out in much greater generality. In \cite{DPGS}, for instance, the existence of probabilistically-weak solutions is proven also when $d = 3$, and the existence of a pressure is established. Moreover, the Lipschitz-continuity 
    assumption on $G_1$ can be relaxed.
    Here we have just recalled the well-posedeness result that is needed in order to investigate the ergodic behaviour of the system.
    For further details, we refer the reader to \cite{DPGS}.
\end{remark}

\subsection{Refined well-posedness results} Before investigating the longtime behavior of the system, we aim to refine Theorem \ref{th:probstrongsol}. Differently from \cite[Proposition 4.1]{DPGS}, we provide here some stability estimates in the stonger energy norms and with the more natural choice $Z=\bHs$.
In this regard, we will work under the following slightly stronger assumptions (see also \cite{OS23}).
\begin{enumerate}[start=5,label=\textbf{(A\arabic*)}]
	\itemsep0.3em
	\item \label{hyp:additionalF} For any $s \geq 1$, consider the function
	\[
	\Psi_s : (-1,1) \to \mathbb R, \qquad \Psi_s(r) = \dfrac{1}{(1-r^2)^s}.
	\]
    The potential density $F:[-1,1] \to \mathbb R$ satisfies
	\[
	F''(r) \leq L_F(1+\Psi_{s_F}(r))
	\]
	for some $s_F \geq 1$ and all $r \in (-1,1)$. Moreover, $F'':(-1,1)\to\erre$ is a quasi-convex function, namely
    it holds that $F''(\theta r_1+(1-\theta)r_2)\leq
    \max\{F''(r_1), F''(r_2)\}$ for all $r_1,r_2\in(-1,1)$
    and $\theta\in[0,1]$.
	\item \label{hyp:additionalAC} The sequence $\{g_k\}_{k \in \mathbb N}$ defined in Assumption \ref{hyp:diffusionAC} satisfies
    \[
    g_k\Psi_{s_0+1} \in L^\infty(-1,1)
    \]
    for every $k\in\enne$, for some $s_0 > 2s_F - 1$, and
    \[
    \sum_{k\in\enne}
    \norm{g_k\Psi_{s_0+1}}_{L^\infty(-1,1)}^2
    <+\infty
    \]
    and we redefine without relabeling the constant $L_{G_2}$
    by adding also the infinite sum above.	
\end{enumerate}

\begin{remark}
   Note that the logarithmic potential satisfies Assumption \ref{hyp:additionalF} with $s_F = 1$. 
\end{remark}

Let us show a technical key lemma first, holding independently of the type of boundary conditions.
\begin{lem} \label{lem:F''}
    Let Assumptions \ref{hyp:structural}-\ref{hyp:additionalAC} hold. 
    Let $(\bu_0, \varphi_0)$ satisfy \eqref{eq:u0}--\eqref{eq:phi0}
    and $\Psi_{s_0}(\varphi_0) \in L^1(\Omega; L^1(\OO))$, and 
    let $(\hphi, \hbu)$ denote a martingale solution to \eqref{AC_NS}. Then, 
    by setting $\gamma:=\frac{s_0+1}{s_F} >2$,
    it holds that $F''(\hphi) \in L^\gamma_\cP(\hom; L^\gamma(0,T;L^\gamma(\OO)))$ and the estimate
    \[
    \|F''(\hphi)\|^\gamma_{ L^\gamma_\cP(\hom; L^\gamma(0,T;L^\gamma(\OO)))} \leq C\left( 1 + \|\Psi_{s_0}(\hphi_0)\|_{L^1(\hom; L^1(\OO))}\right) 
    \]
    holds for a suitable constant $C > 0$. Moreover, we have
    \begin{equation*}
        \int_\OO \Psi_{s_0}(\hphi(t)) + \int_0^t \int_\OO \Psi_{s_0+1}(\hphi(\tau)) \: \d \tau \leq C\left( 1+t+ \int_\OO \Psi_{s_0}(\hphi_0) + \int_0^t \left(\Psi'_{s_0}(\hphi(s)), G_2 (\hphi(s)) \d W(s) \right)_H\right)
    \end{equation*}
    for all $t \geq 0$, $\hP$-almost surely, and for a constant $C>0$ dependent on the structural parameters of the problem, except $\nuu$ and $\nup$.
\end{lem}
\begin{proof}
    The claim follows arguing as in \cite[Proposition 3.1]{OS23}. Indeed, following line by line its proof, we observe that the convective term in the Allen--Cahn equation does not impact the application to the It\^{o} lemma to the integral of a suitable approximating version of $\Psi_{s_0}$. Arguing as therein, we have
	\begin{equation*} 
		\hE \supp \int_\OO \Psi_{s_0}(\hphi(\tau)) + \hE \int_0^t \int_\OO |\Psi'_{s_0}(\hphi(\tau))| \: \d\tau \leq C(1+t) + \E \int_\OO \Psi_{s_0}(\hphi_0)
	\end{equation*}
	and since $|\Psi'_{s_0}(r)| \geq \frac1C\Psi_{s_0+1}(r)$ for some constant $C > 0$ independent of $r \in (-1,1)$, we infer that
	\begin{equation}\label{eq:gs0}
		\hE \supp \int_\OO \Psi_{s_0}(\hphi(\tau)) +
        \frac1C\hE \int_0^t \int_\OO \Psi_{s_0+1}(\hphi(\tau)) \: \d\tau \leq C(1+t) + \hE \int_\OO \Psi_{s_0}(\hphi_0). 
	\end{equation}
	Moreover, by Assumption \ref{hyp:additionalAC} we have $s_0 + 1 > 2s_F$, and therefore $\frac{s_0+1}{s_F} > 2$. In particular, Assumption \ref{hyp:additionalF} implies in turn that
    \[
    \begin{split}
        \hE \int_0^t \|F''(\hphi(\tau)\|_{L^{\frac{s_0+1}{s_F}}}^{\frac{s_0+1}{s_F}} \: \d \tau & \leq C\hE \int_0^t \int_\OO \left( 1 + \Psi_{s_F}(\hphi(\tau)) \right)^{\frac{s_0+1}{s_F}} \: \d \tau \\
        & \leq C\left( 1 + \hE \int_0^t \int_\OO  \Psi_{s_F}^{\frac{s_0+1}{s_F}}(\hphi(\tau))  \: \d \tau\right) \\
        & = C\left( 1 + \hE \int_0^t \int_\OO  \Psi_{s_0+1}(\hphi(\tau))  \: \d \tau\right)
        \leq C\left( 1 + \hE \int_\OO \Psi_{s_0}(\hphi_0)\right), 
    \end{split}
    \]
    where we also used \eqref{eq:gs0}.
    Therefore, we showed that $F''(\hphi) \in L^\gamma(\hom; L^\gamma(0,T;L^\gamma(\OO))$ for $\gamma = \frac{s_0+1}{s_F}> 2$, and the first claim is proved. Moreover, arguing a posteriori from \cite[Proposition 3.1]{OS23} we also have the limiting version of the estimate \cite[(3.3)]{OS23} pointwise in $\Omega$, that is
    \[
        \int_\OO \Psi_{s_0}(\hphi(t)) + \frac1C\int_0^t \int_\OO \Psi_{s_0+1}(\hphi(\tau)) \: \d \tau \leq C(1+t) + \int_\OO \Psi_{s_0}(\hphi_0) + \int_0^t \left(\Psi'_{s_0}(\hphi(s)), G_2 (\hphi(s)) \d W(s) \right)_H
    \]
    and the proof is complete.
\end{proof} 

We are now ready to present the refined continuous dependence estimate, providing stability with respect to the initial data 
in the stronger energy norms and ensuring uniqueness of a (probabilistically-strong) solution. In the particular case of Dirichlet conditions, we only provide a pathwise stability 
estimate that will be used later on.
 
\begin{prop} \label{prop:uniqueness}
	Let Assumptions \ref{hyp:structural}-\ref{hyp:additionalAC} hold with $Z = \bHs$. Consider two sets of initial conditions $(\bu_{0,i}, \varphi_{0,i})$ for $i \in \{1,2\}$
	satisfying \eqref{eq:u0}--\eqref{eq:phi0}
    with $p > 2$ and $\Psi_{s_0}(\varphi_{0,i}) \in L^1(\Omega; L^1(\OO))$. Let $(\widehat\varphi_i, \widehat\bu_i)$ denote two martingale solutions
	to \eqref{AC_NS}, defined
	on the same suitable filtered space $(\hom, \hF, (\hF_t)_t, \hP)$
	and with respect to a pair of Wiener processes $\widehat W_1,\, \widehat W_2$. Then, letting $\gamma := \frac{s_0+1}{s_F}$, there exists a constant $C_\Xi > 0$ only depending on the parameters of the problem such that, defining the process
    $\Xi: [0,T] \to \mathbb{R}$ as
    \begin{equation*}
    \Xi(t) := C_\Xi\int_0^t \left( 1 + \|\hbu_2(\tau)\|^2_{\bHs} \|\hbu_2(\tau)\|^2_{\bVs}+ \|\hphi_2(\tau)\|^2_{V_2} 
    +\sum_{i=1}^2 \|F''(\hphi_i(\tau))\|_{L^\gamma(\OO)}^2\right) \: \d \tau,
    \end{equation*}
    the weighted continuous dependence estimate
	\begin{multline*}
		\hE \supp e^{-\Xi (\tau)}\| \hphi_1(\tau) - \hphi_2(\tau)\|^2_{V_1} + \hE \supp e^{-\Xi (\tau)}\|\hbu_1(\tau)-\hbu_2(\tau)\|^2_{\bHs} \\
		\leq  C\left[ \hE \|\hphi_{0,1}-\hphi_{0,2}\|^2_{V_1} + \hE \|\hbu_{0,1}-\hbu_{0,2}\|^2_{\bHs} \right]
	\end{multline*}
	holds for all $t \in [0,T]$ and for a constant $C > 0$ depending on all the parameters of the problem. In particular, the martingale solution to \eqref{AC_NS} is pathwise unique and probabilistically-strong. 
    \\
    Eventually, 
    in the particular case $\alpha_d = 1$ and $\alpha_n = 0$ (Dirichlet conditions), the estimate
	\begin{multline} \label{eq:cd4_D}
		\dfrac{\nup}{2}\| \hphi_1(t)-\hphi_2(t)\|^2_{V_1} + \dfrac{1}{2}\|\hbu_1(t)-\hbu_2(t)\|^2_{\bH} + \dfrac{1}{2}\int_0^t \left[ \nup^2\|\Delta (\hphi_1 - \hphi_2)(\tau)\|^2_H +\nu\|\nabla(\hbu_1-\hbu_2)(\tau)\|^2_{\bHs} \right] \: \d \tau \\
		\leq  \dfrac{\nup}{2}\| \hphi_{0,1}-\hphi_{0,2}\|^2_{V_1} + \dfrac{1}{2}\|\hbu_{0,1}-\hbu_{0,2}\|^2_{\bHs} \\
        + 2\int_0^t \left[ L_{G_1} + \dfrac{27\KL^8}{4\nu^3}\|\hbu_2(\tau)\|_{\bHs}^2\|\hbu_2(\tau)\|_{\bVs}^2+\dfrac{128}{\nu}\KL^4\KGN^4K_\Delta^4\| \hphi_2(\tau)\|_{V_2}^2 \right]\dfrac{1}{2}\|\hbu_1(\tau)-\hbu_2(\tau)\|_{\bHs}^2 \: \d \tau \\
        +2 \int_0^t  \left[L_{G_2} + \frac{512}{\beta^3}\KL^8K_\Delta^2\|\hbu_2(\tau)\|_{\bHs}^2\|\hbu_2(\tau)\|_{\bVs}^2 + \frac2\beta K_{\frac{2\gamma}{\gamma-2}}^2 \sum_{i=1}^2\|F''(\hphi_i(\tau))\|^2_{L^\gamma(\OO)} \right. \\
        \left.+ \left( L_{G_2}C_4^4 + \dfrac{32\beta}{\nu^2}\KL^4\KGN^4 K_\Delta^2\right)\|\hphi_2(\tau)\|^2_{V_2}\right]\frac\beta2\|\hphi_1(\tau) - \hphi_2(\tau)\|_{V_1}^2 \: \d \tau+ \sum_{k = 2}^{3} \mathcal{M}_k(t),
	\end{multline}
	holds for any $t \in [0,T]$, $\P$-almost surely, where the martingale processes $\{\mathcal M_k\}_{k \in \{2,3\}}$ are defined in \eqref{M1}--\eqref{M3} below.
\end{prop}
\begin{remark}
    The crucial novelty of Proposition~\ref{prop:uniqueness}
    is that the the implicit process $\Xi$ appearing in the continuous dependence estimate
    depends on higher-order norms of only one of the two solutions, namely $(\widehat\varphi_2,,\widehat\bu_2)$ only:
    the dependence of $\Xi$ on the other solution 
    $(\widehat\varphi_1,,\widehat\bu_1)$ is only through terms of order zero instead. 
    This behaviour is usually false for general evolution equations
    as it strongly depends on the geometry of the system, and 
    it is know to hold only in some specific examples, such as 
    the Navier-Stokes equation.
    In our case, such subtle refined stability is 
    absolutely non-trivial due to the coupling with the Allen-Cahn equation and will be fundamental in the subsequent proof of uniqueness of the invariant measure and asymptotic stability.
\end{remark}

\begin{proof}[Proof of Proposition~\ref{prop:uniqueness}]
	For the sake of convenience, let us define the differences
	\begin{align*}
		\hbu := \hbu_1-\hbu_2, \qquad
		\hphi := \hphi_1-\hphi_2, \qquad
		\hw  := \hw_1-\hw_2, \qquad
		\hbu_0  := \hbu_{0,1}-\hbu_{0,2}, \qquad
		\hphi_0  := \hphi_{0,1}-\hphi_{0,2}.
	\end{align*}
	In order to prove the statements, we shall need three applications of the It\^{o} lemma. Throughout the argument, the symbol $\overline v$ denotes the integral average of any $v \in L^1(\OO)$. Applying the It\^{o} lemma for the squared integral average of $\hphi$ yields
	\begin{multline} \label{eq:cd0}
		\dfrac{1}{2}|\overline{\hphi(t)}|^2 + \int_0^t \left[\overline{F'(\hphi_1(\tau))} - \overline{F'(\hphi_2(\tau))}\right]\overline{\hphi(\tau)} \: \d \tau  = 	\dfrac{1}{2}|\overline{\hphi_0}|^2 + 
        \int_0^t \overline{\hphi(\tau)} \left[\overline{G_2( \hphi_1(\tau))} - \overline{G_2( \hphi_2(\tau))}\right]\,\d  W_2(\tau) \\ + \dfrac 12\int_0^t \|\overline{G_2( \hphi_1(\tau))} - \overline{G_2( \hphi_2(\tau))}\|_{\mathcal L_{\text{HS}}(U_2, \erre)}^2 \: \d \tau,
	\end{multline}
	where
	\[
	\overline{G_2}: \BB \to \mathcal L_{\text{HS}}(U_2,\erre) \qquad \overline{G_2( \psi )}[u^2_k] := \overline{g_k(\psi)}
	\]
	for all $k\in \mathbb N$ and $\psi \in \BB$. Secondly, we apply the It\^{o} lemma (see \cite[Theorem 4.2.5]{LiuRo}) to the squared $\bHs$-norm of $\hbu$, yielding
	\begin{multline} \label{eq:cd1}
		\dfrac{1}{2}\|\hbu(t)\|^2_{\bHs} + \int_0^t \left[ \nuu\|\nabla\hbu(\tau)\|^2_{\bHs} - \left( \hw_{1}(\tau)\nabla\hphi_{1}(\tau) - \hw_{2}(\tau)\nabla\hphi_{2}(\tau), \b{u}(\tau)\right)_{\b H} \right]\mathrm{d}\tau \\ 
		+ \int_0^t \langle B(\hbu_1(\tau),\hbu_1(\tau)) - B(\hbu_2(\tau), \hbu_2(\tau)), \hbu \rangle_{\bVsd, \bVs} \: \d \tau
		= \dfrac{1}{2}\|\hbu_{0}\|^2_{\bHs} \\ +\int_0^t\left(\hbu(\tau), \left[G_{1}( \hbu_{1}(\tau)) - G_{1}( \hbu_{2}(\tau)) \right]\,\d  W_1(\tau)\right)_{\bHs}  + \dfrac{1}{2}\int_0^t \|G_{1}( \hbu_{1}(\tau)) - G_{1}( \hbu_{2}(\tau))\|^2_{\mathcal L_{\text{HS}}(U_1, \bHs)} \: \mathrm{d}\tau.
	\end{multline}
	Finally,  we apply the It\^{o} lemma (in its version for twice differentiable functionals, see \cite[Theorem 4.32]{dapratozab}) also to the squared $H$-norm of $\nabla \hphi$. This yields
	\begin{multline} \label{eq:cd2}
		\dfrac{\beta}{2}\|\nabla \hphi(t)\|^2_{H} + \int_0^t \left[ \nup^2\|\Delta \hphi(\tau)\|^2_H - \beta\left( \Delta \hphi(\tau), \nabla\hphi_{1}(\tau)\cdot \b u_1(\tau) - \nabla\hphi_{2}(\tau)\cdot \b{u}_2(\tau)\right)_{H} \right]\mathrm{d}\tau \\ 
		- \beta\int_{0}^{t} (\Delta \hphi(\tau), F'(\hphi_1(\tau)) - F'(\hphi_2(\tau)))_H \: \d \tau \\
		= \dfrac{\beta}{2}\|\nabla \hphi_{0}\|^2_{\b H} -\beta\int_0^t\left(\Delta \hphi(\tau), \left[G_{2}( \hphi_{1}(\tau)) - G_{2}( \hphi_{2}(\tau)) \right]\,\d  W_2(\tau)\right)_{\bHs}  \\+ \dfrac{\beta}{2}\int_0^t \|\nabla G_{2}( \hphi_{1}(\tau)) - \nabla G_{2}( \hphi_{2}(\tau))\|^2_{\mathcal L_{\text{HS}}(U_2, \b H)} \: \mathrm{d}\tau.
	\end{multline}
	By summing equalities \eqref{eq:cd0}-\eqref{eq:cd2}, we end up with
	\begin{multline} \label{eq:cd3}
		\dfrac{1}{2}|\overline{\hphi(t)}|^2 + 	\dfrac{\nup}{2}\|\nabla \hphi(t)\|^2_{H} + \dfrac{1}{2}\|\hbu(t)\|^2_{\bHs} + \int_0^t \left[ \nup^2\|\Delta \hphi(\tau)\|^2_H + \nuu\|\nabla\hbu(\tau)\|^2_{\bH} \right] \: \d \tau \\
		\leq  \dfrac{1}{2}|\overline{\hphi_0}|^2 + \dfrac{\nup}{2}\|\nabla \hphi_{0}\|^2_{\b H} + \dfrac{1}{2}\|\hbu_{0}\|^2_{\bHs} + \sum_{k = 1}^{7} \mathcal{D}_k(t) + \sum_{k = 1}^{3} \mathcal{M}_k(t),
	\end{multline}
	where we defined the deterministic integrals
	\begin{align*}
		\mathcal{D}_1(t) & := \left| \int_0^t \left[\overline{F'(\hphi_1(\tau))} - \overline{F'(\hphi_2(\tau))}\right]\overline{\hphi(\tau)} \: \d \tau \right|, \\
		\mathcal{D}_2(t) & := \dfrac 12\int_0^t \|\overline{G_2( \hphi_1(\tau))} - \overline{G_2( \hphi_2(\tau))}\|_{\mathcal L_{\text{HS}}(U_2, \erre)}^2 \: \d \tau, \\
		\mathcal D_3(t) & := \left| \int_0^t \left[ \nup\left( \Delta \hphi(\tau), \nabla\hphi_{1}(\tau)\cdot \b u_1(\tau) - \nabla\hphi_{2}(\tau)\cdot \b{u}_2(\tau)\right)_{H} - \left( \hw_{1}(\tau)\nabla\hphi_{1}(\tau) - \hw_{2}(\tau)\nabla\hphi_{2}(\tau), \b{u}(\tau)\right)_{\b H} \right]\mathrm{d}\tau \right|, \\
		\mathcal D_4(t) & := \left| \int_0^t \langle B(\hbu_1(\tau),\hbu_1(\tau)) - B(\hbu_2(\tau), \hbu_2(\tau)), \hbu \rangle_{\bVsd, \bVs} \: \d \tau \right|, \\
		\mathcal D_5(t) & := \dfrac{1}{2}\int_0^t \|G_{1}( \hbu_{1}(\tau)) - G_{1}( \hbu_{2}(\tau))\|^2_{\mathcal L_{\text{HS}}(U_1, \bHs)} \: \mathrm{d}\tau, \\
		\mathcal D_6(t) & := \nup\left| \int_{0}^{t} (\Delta \hphi(\tau), F'(\hphi_1(\tau)) - F'(\hphi_2(\tau)))_H \: \d \tau \right|,\\
		\mathcal D_7(t) & := \dfrac{\nup}{2}\int_0^t \|\nabla G_{2}( \hphi_{1}(\tau)) - \nabla G_{2}( \hphi_{2}(\tau))\|^2_{\mathcal L_{\text{HS}}(U_2, \b H)} \: \mathrm{d}\tau,
	\end{align*}
	and the stochastic integrals
	\begin{align}
    \label{M1}
		\mathcal M_1(t) & := \int_0^t 
        \overline{\hphi(\tau)} \left[\overline{G_2( \hphi_1(\tau))} - \overline{G_2( \hphi_2(\tau))}\right]\,\d  W_2(\tau), \\
    \label{M2}
		\mathcal M_2(t) & := \int_0^t\left(\hbu(\tau), \left[G_{1}( \hbu_{1}(\tau)) - G_{1}( \hbu_{2}(\tau)) \right]\,\d  W_1(\tau)\right)_{\bHs}, \\
    \label{M3}
		\mathcal M_3(t) & := \nup\int_0^t\left(\Delta \hphi(\tau), \left[G_{2}( \hphi_{1}(\tau)) - G_{2}( \hphi_{2}(\tau)) \right]\,\d  W_2(\tau)\right)_{H}.
	\end{align}
	Let us tackle the deterministic terms first. Owing to the mean value theorem and the quasi-convexity of $F''$, we have that
	\begin{equation} \label{eq:d1}
		\begin{split}
			\mathcal D_1(t) & = \left| |\OO|^{-2}\int_0^t \int_\OO{\hphi(\tau)} \int_0^1
            \int_\OO F''(\hphi_1(\tau)-\theta\hphi(\tau))\hphi(\tau)\:\d\theta\: \d \tau \right|\\
			 &\leq |\OO|^{-\frac{3}{2}}\int_0^t
             \sum_{i=1}^2\|F''(\hphi_i(\tau))\|_H
             \|\hphi(\tau)\|^2_{V_1} \: \d \tau.
		\end{split}
	\end{equation}
    By exploiting the definition of $\overline{G_2}$, we get
	\begin{equation} \label{eq:d2}
		\begin{split}
			\mathcal D_2(t) &\leq \int_0^t \sum_{k\in\mathbb N} \left|\overline{g_k(\hphi_1(\tau))} - \overline{g_k(\hphi_2(\tau))}    \right|^2 \: \d \tau 
            \leq |\OO|^{-2} \int_0^t \sum_{k\in\mathbb N} \left| \int_\OO g_k(\hphi_1(\tau)) - g_k(\hphi_2(\tau))    \right|^2 \: \d \tau \\
			& \leq |\OO|^{-1} \int_0^t \sum_{k\in\mathbb N} \|g_k(\hphi_1(\tau))- g_k(\hphi_2(\tau))    \|^2_H \: \d \tau 
			\leq |\OO|^{-1}L_{G_2}\int_0^t\|\hphi(\tau)\|^2_{V_1}\: \d \tau.
		\end{split}
	\end{equation}
	As for the third integral, we observe that exploiting the equality
	\[
	\hw_i\nabla\hphi_i = -\nup\Delta\hphi_i\nabla\hphi_i + \nabla \left[ F(\hphi_i) \right]
	\]
	for $i \in \{1,2\}$ and the incompressibility constraint, we get, by the Ladyzhenskaya inequality and integration by parts
	\begin{equation} \label{eq:d3}
		\begin{split}
			\mathcal D_3(t) & =
            \left| \int_0^t \nup\left[ \left( \Delta \hphi(\tau), \nabla\hphi_{1}(\tau)\cdot \b u_1(\tau) - \nabla\hphi_{2}(\tau)\cdot \b{u}_2(\tau)\right)_{H} - \left( \Delta \hphi_{1}(\tau)\nabla\hphi_{1}(\tau) - \Delta \hphi_{2}(\tau)\nabla\hphi_{2}(\tau), \b{u}(\tau)\right)_{\b H} \right]\mathrm{d}\tau \right|\\
            &=\left| \int_0^t - \nup\left(\Delta \hphi(\tau), \nabla \hphi(\tau) \cdot \hbu_2(\tau) \right)_H + \nup\left(\nabla \hphi_2(\tau), D^2 \hphi(\tau) \hbu(\tau) + \nabla \hbu(\tau)\nabla \hphi(\tau) \right)_{\bH} \mathrm{d}\tau \right| \\
			& \leq \nup 
            \int_0^t \|\Delta \hphi(\tau)\|_H\|\nabla \hphi(\tau)\|_{\b L^4(\OO)}\|\hbu_2(\tau)\|_{\b L^4(\OO)} \\
			& \hspace{0.2cm}+ \|\nabla \hphi_2(\tau)\|_{\b L^4(\OO)}\left[ K_\Delta\left( \|\Delta \hphi(\tau)\|_{H} + \|\hphi(\tau)\|_{H}\right) \|\hbu(\tau)\|_{\b L^4(\OO)} + \|\nabla \hphi(\tau)\|_{\b L^4(\OO)}\|\nabla \hbu(\tau)\|_{\bHs} \right]\: \d \tau.
		\end{split}
	\end{equation}
    On account of the fact that
    \begin{equation*}
        \begin{split}
            & \nup \|\Delta \hphi(\tau)\|_H\|\nabla \hphi(\tau)\|_{\b L^4(\OO)}\|\hbu_2(\tau)\|_{\b L^4(\OO)} \\
            & \hspace{2cm} \leq \dfrac{\nup^2}{24}\|\Delta \hphi(\tau)\|^2_H + 6K_\Delta\KL^2\|\hphi(\tau)\|_{V_1}\left( \|\Delta \hphi(\tau)\|_{H} + \|\hphi(\tau)\|_{H}\right)\|\hbu_2(\tau)\|_{\b L^4(\OO)}^2 \\
            & \hspace{2cm} \leq \dfrac{\nup^2}{12}\|\Delta \hphi(\tau)\|^2_H + 
            6 K_\Delta\KL^2\left(1+
            \frac{36}{\beta^2}K_\Delta\KL^2\right)
            \left(1+
            \|\hbu_2(\tau)\|_{\b L^4(\OO)}^4 \right)
            \|\hphi(\tau)\|_{V_1}^2\\
            & \hspace{2cm} \leq \dfrac{\nup^2}{12}\|\Delta \hphi(\tau)\|^2_H + 
            6 K_\Delta\KL^2\left(1+
            \frac{36}{\beta^2}K_\Delta\KL^2\right)
            \left(1+\KL^4
            \|\hbu_2(\tau)\|_{\bHs}^2\|\nabla\hbu_2(\tau)\|_{\bHs}^2 \right)
            \|\hphi(\tau)\|_{V_1}^2,
        \end{split}
    \end{equation*}
    while
    \begin{equation*}
        \begin{split}
            & \nup K_\Delta\|\nabla \hphi_2(\tau)\|_{\b L^4(\OO)}\left( \|\Delta \hphi(\tau)\|_{H} + \|\hphi(\tau)\|_{H}\right) \|\hbu(\tau)\|_{\b L^4(\OO)} \\
            & \hspace{2cm} \leq \dfrac{\nup^2}{12}\| \Delta \hphi(\tau)\|^2_H + \dfrac{\nup^2}{12}\|  \hphi(\tau)\|^2_H +6K_\Delta^2\KGN^2\KL^2\|\hphi_2(\tau)\|_{V_2}\|\hbu(\tau)\|_{\bHs}\|\nabla \hbu(\tau)\|_{\bHs} \\
            & \hspace{2cm} \leq \dfrac{\nup^2}{12}\| \Delta \hphi(\tau)\|^2_H + \dfrac{\nup^2}{12}\|  \hphi(\tau)\|^2_{V_1} + \dfrac{\nuu}{8}\|\nabla \hbu(\tau)\|^2_{\bHs}+\dfrac{72}{\nuu}K_\Delta^4\KGN^4\KL^4\|\hphi_2(\tau)\|_{V_2}^2\|\hbu(\tau)\|^2_{\bHs} \\
        \end{split}
    \end{equation*}
    and finally
    \begin{equation*}
        \begin{split}
            & \nup \|\nabla \hphi_2(\tau)\|_{\b L^4(\OO)}\|\nabla \hphi(\tau)\|_{\b L^4(\OO)}\|\nabla \hbu(\tau)\|_{\bHs} \\
            & \leq \dfrac{\nuu}{8}\|\nabla \hbu(\tau)\|^2_{\bHs} + \dfrac{2}{\nuu}\nup^2K_\Delta\KGN^2\KL^2\|\hphi(\tau)\|_{V_1}\left( \|\Delta \hphi(\tau)\|_{H} + \|\hphi(\tau)\|_{H}\right)\|\hphi_2(\tau)\|_{V_2}\\
            &  \leq \dfrac{\nuu}{8}\|\nabla \hbu(\tau)\|^2_{\bHs} + \dfrac{\nup^2}{12}\| \Delta \hphi(\tau)\|^2_H + 
            \dfrac{2\beta^2}{\nuu}
            K_\Delta\KGN^2\KL^2
            \left[ 1 + \dfrac{6\beta^2}{\nuu} K_\Delta\KGN^2\KL^2\right]
            \left(1+\|\hphi_2(\tau)\|_{V_2}^2\right)
            \|\hphi(\tau)\|_{V_1}^2,
        \end{split}
    \end{equation*}
    we can infer from \eqref{eq:d3} that
    \begin{multline} \label{eq:d32}
        \mathcal D_3(t) \leq \dfrac{\nuu}{4}\int_0^t\|\nabla \hbu(\tau)\|^2_{\bHs}\:\d\tau + 
        \dfrac{\nup^2}{4}\int_0^t\| \Delta \hphi(\tau)\|^2_H\:\d\tau + \dfrac{\nup^2}{12}
        \int_0^t\|  \hphi(\tau)\|^2_{V_1} \:\d\tau\\
        +\dfrac{72}{\nu}K_\Delta^4\KGN^4\KL^4
        \int_0^t
        \|\hphi_2(\tau)\|_{V_2}^2\|\hbu(\tau)\|^2_{\bHs}
        \:\d\tau\\
        + \dfrac{2\beta^2}{\nuu}
            K_\Delta\KGN^2\KL^2
            \left[ 1 + \dfrac{6\beta^2}{\nuu} K_\Delta\KGN^2\KL^2\right]
            \int_0^t
            \left(1+\|\hphi_2(\tau)\|_{V_2}^2\right)
            \|\hphi(\tau)\|_{V_1}^2\:\d\tau \\ 
        + 6 K_\Delta\KL^2\left(1+
            \frac{36}{\beta^2}K_\Delta\KL^2\right)
            \int_0^t
            \left(1+\KL^4
            \|\hbu_2(\tau)\|_{\bHs}^2\|\nabla\hbu_2(\tau)\|_{\bHs}^2 \right)
            \|\hphi(\tau)\|_{V_1}^2\:\d\tau .
    \end{multline}
	Next, we have
	\begin{equation} \label{eq:d4}
		\begin{split}
			\mathcal D_4(t) & 
			 = \left| \int_0^t \langle B(\hbu(\tau), \hbu_2(\tau)), \hbu \rangle_{\bVsd, \bVs} \: \d \tau \right| 
			 \leq \int_0^t  \|\hbu(\tau)\|_{\bVs}\|\hbu_2(\tau)\|_{\b L^4(\OO)}\|\hbu(\tau)\|_{\b L^4(\OO)} \: \d \tau \\
			& \leq\KL \int_0^t \|\hbu(\tau)\|_{\bVs}^\frac 32\|\hbu_2(\tau)\|_{\b L^4(\OO)}\|\hbu(\tau)\|_{\bHs}^\frac 12 \: \d \tau \\
            & \leq \frac \nu4\int_0^t \|\hbu(\tau)\|_{\bVs}^2 \: \d \tau + \dfrac{27\KL^8}{4\nu^3}\int_0^t\|\hbu_2(\tau)\|_{\bHs}^2\|\hbu_2(\tau)\|_{\bVs}^2\|\hbu(\tau)\|_{\bHs}^2 \: \d \tau,
		\end{split}
	\end{equation}
	while owing to Assumption \ref{hyp:diffusionNS} we conclude
	\begin{equation} \label{eq:d5}
		\mathcal D_5(t) = \int_0^t \|G_{1}( \hbu_{1}(\tau)) - G_{1}( \hbu_{2}(\tau))\|^2_{\mathcal L_{\text{HS}}(U_1, \bHs)}\:\d \tau \leq L_{G_1}\int_0^t \|\hbu(\tau)\|^2_{\bHs} \: \d \tau.
	\end{equation}
	Finally, owing to the mean value theorem, the quasi-convexity of $F''$ and the two-dimensional Sobolev embedding $V \hookrightarrow L^r(\OO)$ for any finite $r \geq 1$, as well as Lemma \ref{lem:F''},
	\begin{equation} \label{eq:d6}
		\begin{split}
			\mathcal D_6(t) &  \leq \nup\int_0^t \|\Delta \hphi(\tau)\|_H\|\hphi(\tau)\|_{L^{\frac{2\gamma}{\gamma-2}}(\OO)}\sum_{i=1}^2\|F''(\hphi_i(\tau))\|_{L^\gamma(\OO)} \: \d \tau \\
			& \leq \dfrac{\nup^2}{4}\int_0^t\|\Delta \hphi(\tau)\|^2_H \: \d \tau + 2K_{\frac{2\gamma}{\gamma-2}}^2\int_0^t \sum_{i=1}^2\|F''(\hphi_i(\tau))\|^2_{L^\gamma(\OO)}\|\hphi(\tau)\|^2_{V_1} \: \d \tau,
		\end{split}
	\end{equation}
    while on account of Assumption \ref{hyp:diffusionAC}
	\begin{equation} \label{eq:d7}
		\begin{split}
			\mathcal{D}_7(t) & = \frac \nup 2\int_0^t \sum_{k\in\enne} \|g'_k(\hphi_1(\tau))\nabla \hphi_1(\tau) -g'_k(\hphi_2(\tau))\nabla \hphi_2(\tau)\|^2_{\b H} \: \d \tau \\
			& \leq \nup\int_0^t \sum_{k\in\enne} 
            \left(
            \norm{\left[ g'_k(\hphi_1(\tau)) - g'_k(\hphi_2(\tau))\right]\nabla \hphi_2(\tau)}_{\bH}^2 + \norm{g'_k(\hphi_1(\tau))\nabla \hphi(\tau)}^2_{\b H}\right)
			\\
			& \leq \nup L_{G_2}\int_0^t\left[\|\hphi(\tau)\|^2_{L^4(\OO)}\|\nabla\hphi_2(\tau)\|_{\b L^4(\OO)}^2 + \|\nabla \hphi(\tau)\|^2_{\bH} \right]  \: \d \tau\\
			& \leq
			\nup L_{G_2}\int_0^t\left[1+K_4^4 \|\hphi_2(\tau)\|^2_{V_2}\right]\| \hphi(\tau)\|^2_{V_1} \: \d \tau.
		\end{split}
	\end{equation}
	Collecting \eqref{eq:d1}--\eqref{eq:d7} in \eqref{eq:cd3} and recalling that $v \mapsto (|\overline{v}|^2 + \|\nabla v\|^2_{\b H})^\frac 12$ is an equivalent norm in $V_1$ yield
	\begin{multline} \label{eq:cd4}
		\dfrac{\min\{1,\nup\}}{2}\|\hphi(t)\|^2_{V_1} + \dfrac{1}{2}\|\hbu(t)\|^2_{\bHs} + \int_0^t \left[ \frac {\nup^2 }{2}\|\Delta \hphi(\tau)\|^2_H + \frac \nu2\|\nabla\hbu(\tau)\|^2_{\bHs} \right] \: \d \tau 
		\leq C\left(\|\hphi_0\|^2_{V_1} + \|\hbu_{0}\|^2_{\bHs}\right)\\
        +C\left[ \int_0^t \left(1 + \sum_{i=1}^2 \|F''(\hphi_i(\tau))\|_{L^\gamma(\OO)}^2 + \|\hbu_2(\tau)\|^2_{\bHs} \|\hbu_2(\tau)\|^2_{\bVs} + \| \hphi_2(\tau)\|^2_{V_2}\right)\|\hphi(\tau)\|^2_{V_1}\: \d \tau \right. \\ \l\left. + \int_0^t \left(1 + \|\hbu_2(\tau)\|^2_{\bHs} \|\hbu_2(\tau)\|^2_{\bVs} + \| \hphi_2(\tau)\|^2_{V_2}\right)\|\hbu(\tau)\|^2_{\bHs} \: \d \tau 
		+ \sum_{k = 1}^{3} \mathcal{M}_k(t)\right],
	\end{multline}
    for a suitable constant $C > 0$ only depending on the structural parameters of the problem.
	This is the first estimate, up to applying interpolation inequalities to the right hand side. In order to show the exponential estimate, for any $\sigma > 0$, consider the auxiliary process
	\[
	\Xi(t) := \sigma\int_0^t \left(1+ \sum_{i=1}^2 \|F''(\hphi_i(\tau))\|_{L^\gamma(\OO)}^2 + \|\hbu_2(\tau)\|^2_{\bHs} \|\hbu_2(\tau)\|^2_{\bVs}+ \| \hphi_2(\tau)\|^2_{V_2} \right) \: \d \tau.
	\]
	The regularity of solutions (see Definition \ref{def:strong_sol}) as well as Lemma \ref{lem:F''} shows that $\Xi$ is well defined. Applying the It\^{o} lemma to the process
	\[
	t \mapsto \dfrac 12 e^{-\Xi (t)}\left( |\overline{\hphi(t)}|^2 + \|\nabla \hphi(t)\|^2_{\bH} + \|\hbu(t)\|^2_{\bHs}\right)
	\]
	and following the computations leading to \eqref{eq:cd3}, we have
	\begin{multline} \label{eq:cd5}
		\dfrac{1}{2} e^{-\Xi (t)}|\overline{\hphi(t)}|^2 + 	\dfrac{\nup}{2} e^{-\Xi (t)}\|\nabla \hphi(t)\|^2_{H} + \dfrac{1}{2} e^{-\Xi (t)}\|\hbu(t)\|^2_{\bHs} \\ +\dfrac{\sigma}{2} \int_0^t  \left(1+ \sum_{i=1}^2 \|F''(\hphi_i(\tau))\|_{L^\gamma(\OO)}^2 + \|\hbu_2(\tau)\|^2_{\bHs} \|\hbu_2(\tau)\|^2_{\bVs}+ \| \hphi_2(\tau)\|^2_{V_2} \right)  \\
        \times e^{-\Xi (\tau)}\left( |\overline{\hphi(\tau)}|^2 + \|\nabla \hphi(\tau)\|^2_{\bH} + \|\hbu(\tau)\|^2_{\bHs}\right) \: \d \tau \\ + \dfrac{1}{2}\int_0^t \left[ \nup^2 e^{-\Xi (\tau)}\|\Delta \hphi(\tau)\|^2_H +  \nu e^{-\Xi (\tau)}\|\nabla\hbu(\tau)\|^2_{\bHs} \right] \: \d \tau  \\
		\leq  \dfrac{1}{2}|\overline{\hphi_0}|^2 + \dfrac{\nup}{2}\|\nabla \hphi_{0}\|^2_{\b H} + \dfrac{1}{2}\|\hbu_{0}\|^2_{\bHs} + \sum_{k = 1}^{7} \widetilde{\mathcal{D}}_k(t) + \sum_{k = 1}^{3} \widetilde{\mathcal{M}}_k(t),
	\end{multline}
    where the integrals $\widetilde{\mathcal{D}}_k$ and $\widetilde{\mathcal{M}}_k$ are defined as before, with the additional factor $e^{-\Xi(\tau)}$ inside the time integral. Since $0 < e^{-\Xi(t)}  \leq 1$ for all $t \geq 0$, we can follow the computations throughout \eqref{eq:d1}-\eqref{eq:d7} to estimate the deterministic integrals, and arrive hence to the analogous relation of \eqref{eq:cd4}, that is
	\begin{multline} \label{eq:cd6}
			\dfrac{\min\{1,\nup\}}{2} e^{-\Xi (t)}\| \hphi(t)\|^2_{V_1} + \dfrac{1}{2} e^{-\Xi (t)}\|\hbu(t)\|^2_{\bHs}
            + \dfrac{1}{2}\int_0^t \left[ \nup^2 e^{-\Xi (\tau)}\|\Delta \hphi(\tau)\|^2_H +  \nu e^{-\Xi (\tau)}\|\nabla\hbu(\tau)\|^2_{\bHs} \right] \: \d \tau
            \\ + \rho \int_0^t  \left(1+ \sum_{i=1}^2 \|F''(\hphi_i(\tau))\|_{L^\gamma(\OO)}^2 + \|\hbu_2(\tau)\|^2_{\bHs} \|\hbu_2(\tau)\|^2_{\bVs}+ \| \hphi_2(\tau)\|^2_{V_2} \right)\\
            \times e^{-\Xi (\tau)}\left( \| \hphi(\tau)\|^2_{V_1} + \|\hbu(\tau)\|^2_{\bHs}\right) \: \d \tau   
			\leq  C\left[ \|\hphi_{0}\|^2_{V_1} + \|\hbu_{0}\|^2_{\bHs}   +  \sum_{k = 1}^{3} \widetilde{\mathcal{M}}_k(t)\right],
	\end{multline}
    where $\rho > 0$ is arbitrary, and we selected accordingly $\sigma = 2C$, with $C$ being the constant appearing in \eqref{eq:cd4}. Taking absolute values at right hand side, supremum in time and expectations, we are now in a position to deal with the stochastic integrals and close the estimate. 
        Owing to the Burkholder--Davis--Gundy inequality
	\begin{equation} \label{eq:m1}
		\begin{split}
			&\hE \sups \widetilde{\mathcal M}_1(s)  \leq C\hE \left| \int_0^t e^{-2\Xi(\tau)}\|\hphi(\tau)\|_H^2\|\overline{G_2( \hphi_1(\tau))} - \overline{G_2( \hphi_2(\tau))}\|_{\mathcal L_{\text{HS}}(U_2,\mathbb R)}^2\,\d \tau\right|^\frac 12 \\
			&\qquad \leq C\hE \left| \int_0^t e^{-2\Xi(\tau)}\|\hphi(\tau)\|_{V_1}^4\,\d \tau\right|^\frac 12 
			 \leq C\hE \left[ \supp e^{-\frac 12\Xi(\tau)}\|\hphi(\tau)\|_{V_1} \left|  \int_0^t e^{-\Xi(\tau)}\|\hphi(\tau)\|_{V_1}^2\,\d s\right|^\frac 12\right] \\
			&\qquad \leq \dfrac{\min \{1, \nup\}}{8} \hE \supp e^{- \Xi(\tau)}\|\hphi(\tau)\|_{V_1}^2 + C\hE  \int_0^t e^{-\Xi(\tau)}\|\hphi(\tau)\|_{V_1}^2\,\d s,
		\end{split}
	\end{equation}
	and, similarly,
	\begin{equation} \label{eq:m2}
		\begin{split}
			&\hE \sups \widetilde{\mathcal M}_2(s) 
            \leq C\hE \left| \int_0^t e^{-2\Xi(\tau)}
            \|\hbu(\tau)\|_{\bH}^2\|G_{1}( \hbu_{1}(\tau)) - G_{1}( \hbu_{2}(\tau))\|_{\mathcal L_{\text{HS}}(U_1,\mathbb \bHs)}^2\,\d \tau\right|^\frac 12 \\
			&\qquad \leq C\hE \left| \int_0^t e^{-2\Xi(\tau)}\|\hbu(\tau)\|_{\bHs}^4\,\d \tau\right|^\frac 12 
			 \leq C\hE \left[ \supp e^{-\frac 12\Xi(\tau)}\|\hbu(\tau)\|_{\bHs} \left|  \int_0^t e^{-\Xi(\tau)}\|\hbu(\tau)\|_{\bHs}^2\,\d \tau\right|^\frac 12\right] \\
			&\qquad \leq \frac 14 \hE \supp e^{- \Xi(\tau)}\|\hbu(\tau)\|_{\bHs}^2 + C\hE \int_0^t e^{-\Xi(\tau)}\|\hbu(\tau)\|_{\bHs}^2\,\d \tau,
		\end{split}
	\end{equation}
	and finally
	\begin{equation} \label{eq:m3}
		\begin{split}
			&\hE \sups \widetilde{\mathcal M}_3(s) 
            \leq C\hE \left| \int_0^t e^{-2\Xi(\tau)}\|\nabla\hphi(\tau)\|_{H}^2\|\nabla G_{2}( \hphi_{1}(\tau)) - \nabla G_{2}( \hphi_{2}(\tau))\|_{\mathcal L_{\text{HS}}(U_2,H)}^2\,\d \tau\right|^\frac 12 \\
			&\qquad \leq C\hE \left| \int_0^t e^{-2\Xi(\tau)}\left[1+\|\nabla\hphi_2(\tau)\|^4_{ \b L^4(\OO)}\right]\| \hphi(\tau)\|^4_{V_1}\,\d \tau\right|^\frac 12 \\
			&\qquad \leq C\hE \left[ \supp e^{-\frac 12\Xi(\tau)}\|  \hphi(\tau)\|_{V_1}\left|  \int_0^t e^{-\Xi(\tau)}\left[1+\|\nabla\hphi_2(\tau)\|^4_{ \b L^4(\OO)}\right]\| \hphi(\tau)\|^2_{V_1}\,\d \tau\right|^\frac 12\right] \\
			&\qquad \leq \dfrac{\min \{1, \nup\}}{8}  \hE \supp e^{- \Xi(\tau)}\|\hphi(\tau)\|_{V_1}^2 
            + C\hE  \int_0^t e^{-\Xi(\tau)}\left[1+\| \hphi_2(\tau)\|^2_{H^2(\OO)}\right]\| \hphi(\tau)\|^2_{V_1}\,\d \tau.
		\end{split}
	\end{equation}
	Choosing $\rho$ sufficiently large and collecting \eqref{eq:m1}-\eqref{eq:m3} in \eqref{eq:cd6} yield the desired continuous dependence estimate. Uniqueness follows then standardly by choosing $\hphi_{0,1} = \hphi_{0,2}$ and $\hbu_{0,1} = \hbu_{0,2}$ and recalling that the exponential is positive.\\
    It is only left to show \eqref{eq:cd4_D}, under the assumption of Dirichlet conditions. We shall follow the computations performed above and highlight modifications if need be. First of all, since $v \mapsto \|\nabla v\|^2_H$ is the norm in $V_1$, we can avoid using the It\^{o} lemma for the integral average of $\hphi$. Therefore, we only apply it to the squared $\b H$-norm of $\hbu$ and $\nabla \hphi$, obtaining exactly the equalities \eqref{eq:cd1} and \eqref{eq:cd2}. Summing the two equalities, multiplying the latter by $\beta$, we obtain the analogous version of \eqref{eq:cd3}, that is
    \begin{multline} \label{eq:cd3_D}
		\dfrac{\beta}{2}\| \hphi(t)\|^2_{V_1} + \dfrac{1}{2}\|\hbu(t)\|^2_{\bHs} + \int_0^t \left[ \beta^2\|\Delta \hphi(\tau)\|^2_H +\nu\|\nabla\hbu(\tau)\|^2_{\bHs} \right] \: \d \tau \\
		\leq   \dfrac{\beta}{2}\| \hphi_{0}\|^2_{V_1} + \dfrac{1}{2}\|\hbu_{0}\|^2_{\bHs} + \sum_{k = 3}^{7} \mathcal{D}_k(t) + \sum_{k = 2}^{3} \mathcal{M}_k(t).
	\end{multline}
    We can estimate the deterministic integrals following the same computations of \eqref{eq:d4}-\eqref{eq:d7} with very minor modifications, leading to
    \begin{align*}
			\mathcal D_4(t) & \leq \frac \nu4\int_0^t \|\hbu(\tau)\|_{\bVs}^2 \: \d \tau + \dfrac{27\KL^8}{4\nu^3}\int_0^t\|\hbu_2(\tau)\|_{\bHs}^2\|\hbu_2(\tau)\|_{\bVs}^2\|\hbu(\tau)\|_{\bHs}^2 \: \d \tau, \\
		  \mathcal D_5(t) &  \leq L_{G_1}\int_0^t \          \|\hbu(\tau)\|^2_{\bHs} \: \d \tau, \\
			\mathcal D_6(t) & \leq \dfrac{\beta^2}{4}\int_0^t\|\Delta \hphi(\tau)\|^2_H \: \d \tau + 
            2K_{\frac{2\gamma}{\gamma-2}}^2\int_0^t \sum_{i=1}^2\|F''(\hphi_i(\tau))\|^2_{L^\gamma(\OO)}\|\hphi(\tau)\|^2_{V_1} \: \d \tau, \\
			\mathcal{D}_7(t) & \leq
			\beta L_{G_2}\int_0^t\left[
            1+K_4^4\|\hphi_2(\tau)\|_{ V_2}^2\right]\| \hphi(\tau)\|^2_{V_1} \: \d \tau. 
    \end{align*}
    As for $\mathcal D_3$, computations are highly simplified as here one has directly $\|D^2\varphi\|_{\b H} \leq K_\Delta \|\Delta \varphi\|_H$ thanks to the homogeneous Dirichlet boundary conditions. Therefore, we get
    \begin{equation} \label{eq:d3_D}
		\begin{split}
			&\mathcal D_3(t)  = \left| \int_0^t - \beta\left(\Delta \hphi(\tau), \nabla \hphi(\tau) \cdot \hbu_2(\tau) \right)_H + \beta\left(\nabla \hphi_2(\tau), D^2 \hphi(\tau) \hbu(\tau) + \nabla \hbu(\tau)\nabla \hphi(\tau) \right)_{\bH} \mathrm{d}\tau \right| \\
			& \leq \beta \int_0^t \|\Delta \hphi(\tau)\|_H\|\nabla \hphi(\tau)\|_{\b L^4(\OO)}\|\hbu_2(\tau)\|_{\b L^4(\OO)} \\
			& \hspace{1cm}+ \|\nabla \hphi_2(\tau)\|_{\b L^4(\OO)}\left[ K_\Delta\|\Delta \hphi(\tau)\|_{H}\|\hbu(\tau)\|_{\b L^4(\OO)} + \|\nabla \hphi(\tau)\|_{\b L^4(\OO)}\|\nabla \hbu(\tau)\|_{\bHs} \right]\: \d \tau \\
            & \leq \dfrac{\beta^2}{16} \int_0^t \|\Delta \hphi(\tau)\|_H^2 \: \d \tau + \dfrac{\nu}{8} \int_0^t \|\nabla \hbu(\tau)\|_{\bHs}^2 \: \d \tau
            + 8\KL^2\KGN^2K_\Delta^2\int_0^t \| \hphi_2(\tau)\|_{V_2}\|\hbu(\tau)\|_{\bHs}\|\hbu(\tau)\|_{\bVs} \: \d \tau\\
            & \hspace{1cm}+ 8\KL^4K_\Delta\int_0^t \| \hphi(\tau)\|_{V_1}\| \Delta \hphi(\tau)\|_{H}\|\hbu_2(\tau)\|_{\bHs}\|\hbu_2(\tau)\|_{\bVs} \: \d \tau \\
            & \hspace{1cm}+ \dfrac{2\beta^2}{\nu}\KL^2\KGN^2K_\Delta\int_0^t \| \hphi_2(\tau)\|_{V_2}\|\hphi(\tau)\|_{V_1}\|\Delta \hphi(\tau)\|_{H} \: \d \tau  \\
            & \leq \dfrac{\beta^2}{8} \int_0^t \|\Delta \hphi(\tau)\|_H^2 \: \d \tau 
            + \dfrac{\nu}{4} \int_0^t \|\nabla \hbu(\tau)\|_{\bHs}^2 \: \d \tau 
            + \dfrac{128}{\nu}\KL^4\KGN^4K_\Delta^4\int_0^t \| \hphi_2(\tau)\|_{V_2}^2\|\hbu(\tau)\|_{\bHs}^2 \: \d \tau\\
            & \hspace{1cm}+ \int_0^t \left(\frac{512}{\beta^2}\KL^8K_\Delta^2\|\hbu_2(\tau)\|_{\bHs}^2\|\hbu_2(\tau)\|_{\bVs}^2+
            \dfrac{32\beta^2}{\nu^2}\KL^4\KGN^4K_\Delta^2\| \hphi_2(\tau)\|_{V_2}^2
            \right)\| \hphi(\tau)\|_{V_1}^2 \: \d \tau.
		\end{split}
	\end{equation}
 Combining the previous inequalities, jointly with \eqref{eq:d3_D}, into \eqref{eq:cd3_D}, we arrive at \eqref{eq:cd4_D}, up to multiplying and dividing by a suitable constant the time integrals to the right hand side. This conlcudes the proof.
\end{proof}


\section{Existence of invariant measures}
\label{sec:existence}
In this section, we investigate the \review{longtime behavior} of the system \eqref{AC_NS} in terms of \review{the} existence of invariant measures. More precisely,
we work under assumptions that ensure uniqueness of solutions 
for the system \eqref{AC_NS}, namely either Assumptions \ref{hyp:structural}--\ref{hyp:diffusionAC} with $Z=\bVsd$ or Assumptions \ref{hyp:structural}--\ref{hyp:additionalAC} with $Z=\bHs$. This allows to introduce the transition semigroup 
associated to \eqref{AC_NS}. This \review{can not be achieved by means of standard arguments} due to the nonlinear condition on the initial data needed to solve the SPDE: indeed, while the initial value of the 
Navier--Stokes solution component is assumed to \review{belong to} a linear (Hilbert) space \review{(cf. \eqref{eq:u0})}, the initial value of the Allen--Cahn \review{order parameter} must be taken in a \review{nonlinear subset of $L^p(\Omega, \cF_0;V_1)$, due to the integrability constraint on $F(\varphi_0)$ (cf. \eqref{eq:phi0})}. This forces to define the transition semigroup on functions on a general metric space, and not \review{only on} a Hilbert space as in more classical settings. \review{The proof of existence of invariant measures shall be divided in steps for better clarity. Before moving on, we point out that} from now on, for simplicity of notation, the initial \review{velocity field} will be denoted by $\b x$ instead of $\bu_0$, and 
initial \review{order parameter} will be denoted by $y$ instead of $\varphi_0$. In this direction, we set 
\begin{equation*}
    \b X:=\bHs, \qquad  Y:= \mathcal{A} \cap V_1,
\end{equation*}
\review{with $\mathcal A$ being the unitary ball of $L^\infty(\OO)$}, for the spaces of initial data. \review{Moreover,} for every $(\b x,y)\in \b X\times Y$ and for all $s\geq0$, we use the notation 
$(\bu^{s,\b x},\varphi^{s,y})$ to denote the unique \review{strong} solution
to system \eqref{AC_NS} \review{such that $(\b u (s),\,\varphi(s)) = (\b x,\,y)$.} As usual, we set $(\bu^{\b x},\,\varphi^y):=(\bu^{0,\b x},\,\varphi^{0,y})$.
\subsection{The transition semigroup}
\review{Recalling Assumption \ref{hyp:potential}, we readily have the characterization}
\begin{equation}
\label{mat_A}
\mathcal{A}=\{y \in \review{L^\infty(\OO)} : \|y\|_{L^\infty(\mathcal{O})}\leq 1\} 
= \{y \in H : F(y)\in L^1(\mathcal{O})\}=
\bigcup_{n\in\enne}\left\{y\in H: \ \int_\OO F(y) \leq n\right\},
\end{equation}
which implies that $\mathcal A$ is a Borel subset of $H$ since $F$ is convex and lower semicontinuous.
More specifically, we note that $\mathcal A$ is a closed convex subset of $H$, hence also weakly closed in $H$, and, \review{analogously,} that $\mathcal{A} \cap V_1$ is a closed convex subset of $V_1$, hence also weakly closed in $V_1$. In particular, $\mathcal A\cap V_1$ is a Borel subset of $V_1$. \review{In our framework, due to the singular nature of the problem,} the transition semigroup can only make sense as a family of operators acting on some
suitable space of functions on $\b X \times Y$.
This requires first to address the issue of measurability on $\b X\times Y$. \review{To this end, }let us \review{first} focus on the possible topologies that one can consider on $\b X\times Y$.
As for the space $\b X$, the situation is classical: \review{it can be endowed with} the strong topology $\tau_{\text s}^\b X$, the weak topology $\tau_{\text w}^\b X$, and the bounded-weak 
topology $\tau_{\text{bw}}^\b X$. \review{On the other hand, }for the space $Y$,
we consider the respective topologies of $V_1$ induced on $Y$:
more precisely, $\tau_{\text s}^Y$ \review{(resp. $\tau_{\text w}^Y$, $\tau_{\text{bw}}^Y$)} is the topology induced on $Y$ from 
the strong \review{(resp. weak, bounded-weak)} topology of $V_1$. With this notation, we introduce the topological spaces
\begin{equation*}
    \Xb:=(\b X, \tau^\b X_\bullet), \qquad  
    \Yb:=(Y, \tau^Y_\bullet),
\end{equation*}
as well as the respective topological product space 
$\b X_\bullet \times Y_\bullet$,
where $\bullet \in \{\text{s},\, \text{w},\,\text{bw}\}$. \review{Recalling Subsection \ref{ssec:weaktop}, we have that the Borel $\sigma$-algebra on the product space $\Xb \times \Yb$ does not depend on the choice of the topologies on $\b X$ and $Y$ and shall therefore be denoted simply by $\mathscr B(\b X \times Y)$.} 

We \review{are now in a position} to rigorously introduce 
the family of transition operators
associated to the system \eqref{AC_NS}.
To this end, \review{observe} that for every $(\b x,y)\in \b X\times Y$ 
one has that $(\bu^{\b x}(t),\varphi^y(t))\in \b X\times Y$
for every $t\geq0$, $\P$-almost surely,
and that $(\bu^{\b x}(t),\varphi^y(t))$ is $\mathcal F_t/\mathscr B(\b X\times Y)$-measurable.
Hence, \review{given $t > 0$}, the law of $(\bu^{\b x}(t),\varphi^y(t))$ on $\b X \times Y$ \review{is defined}
as the pushforward measure $\cL_t((\b x,y), \cdot):=\Law_\P(\bu^{\b x}(t),\varphi^y(t))(\cdot)$. \review{Indeed, letting $\cP(\b X \times Y)$ the set of all probability measures on $(\b X \times Y,\,\mathscr{B}(\b X \times Y))$, we have
\[
\cL_t: (\b X\times Y) \times \mathscr B(\b X \times Y) \to [0,1],
\]
and
\[
\cL_t((\b x,y), A):=\P((\bu^{\b x}(t),\varphi^y(t))\in A),
\quad (\b x,y)\in \b X\times Y,
\quad A\in\mathscr B(\b X\times Y).
\]}
By virtue of Theorem~\ref{th:probstrongsol}, 
pathwise uniquenss holds for the problem \eqref{AC_NS}, hence \cite[Corollary 23]{On2005} applies, ensuring that 
$(x,y)\mapsto\cL_t((x,y), A)$ is $\mathscr B(X\times Y)$-measurable
for all $A\in \mathscr B(X\times Y)$. This implies that, for every 
$t\geq0$, $\cL_t$ is a transition 
kernel from $(\b X\times Y,\, \mathscr B(\b X\times Y))$ into itself. Hence, 
we can define the family of operators
$P=\{P_t\}_{t\geq0}$ as
\review{
\begin{equation}
\label{P_t}
P_t : \mathcal{B}_b(\b X\times Y) \to \mathcal{B}_b(\b X\times Y), \quad 
(P_t\phi)(\bx,y):= 
\int_{\b X\times Y}\phi(w)\: \cL_t((\bx,y), \d w)=
\mathbb{E}[ (\phi(\bu^{\b x}(t), \varphi^y(t))]
\end{equation}
for all $(\b x,y)\in \b X\times Y$ and $\phi\in \mathcal{B}_b(\b X\times Y)$. Observe that $P_t\phi\in\mathcal B_b(\b X\times Y)$
for every $\phi \in \mathcal{B}_b(\b X \times Y)$ since $\cL_t$ is a transition kernel. Moreover,} we can also define the family of operators $P^*=\{P^*_t\}_{t\ge 0}$ as
\begin{equation}
\label{P^*}
P_t^*: \cP(\b X \times Y) \to \cP(\b X \times Y), \quad 
P_t^*\mu(A):=\int_{\b X \times Y}\cL_t((\bx, y),A)\, \mu({\rm d}\bx, {\rm d}y), 
\end{equation}
for all $A \in \mathscr{B}(\b X \times Y)$ and $\mu \in \cP(X \times Y)$. In particular, we recall that (see e.g. \cite[Proposition 11.1]{dapratozab}), for any $t\ge 0$ and $(\bx, y) \in X \times Y$
\begin{equation}
\label{ossP^*}    
P_t^*\delta_{(\bx, y)}(\cdot)=\cL_t((\bx, y), \cdot).
\end{equation}
\subsection{Feller properties}
The aim of this \review{sub}section is to prove that the family of operators $P$ \review{possesses some} Feller properties, \review{in the sense defined below}.
\begin{defin}
\label{swF}
    The family of operators $P$ has the Feller property 
    (or is Feller) if 
        \begin{equation*}
P_t: \C_b(\Xs \times \Ys) \rightarrow \C_b(\Xs \times \Ys) \qquad \forall\,t \ge 0.
    \end{equation*}
    The family of operators $P$ has the weak Feller property 
    (or is weak Feller) if
\begin{equation*}
P_t: \C_b(\Xw \times \Yw) \rightarrow \C_b(\Xw \times \Yw) \qquad \forall\,t \ge 0.
    \end{equation*}
The family of operators $P$ has the sequential weak Feller property 
    (or is sequentially weak Feller) if
    \begin{equation*}
P_t: \S\C_b(\Xw \times \Yw) \rightarrow \S\C_b(\Xw \times \Yw) \qquad \forall\,t \ge 0.
    \end{equation*}
\end{defin}
Before moving forward, we state a technical lemma about topological and sequential continuity.
\begin{lem}
\label{eq_spaces}
It holds that 
$\S\C(\Xw \times \Yw)=
\C(\Xbw \times \Ybw)$ and 
$\S\C_b(\Xw \times \Yw)=
\C_b(\Xbw \times \Ybw)$.
\end{lem}
\begin{proof}
For all $n \in \mathbb{N}$, \review{define} $B_n:= 
\{(\bx,y) \in \b X\times V_1: \ \|(\bx,y)\|_{\bHs \times  V_1} \le n\}$\review{, i.e.,} the closed ball of radius $n$ in  $\b X\times V_1$. \review{Moreover, set} $K_n:=B_n\cap(\b X\times\mathcal A)\subset \b X\times Y$.
Since $\b X \times  V_1$ is a separable Hilbert space, 
one has that \review{the product topology} $\tw{\b X} \otimes \tw{V_1}$ is metrizable on $B_n$ \review{for all $n \in \enne$} (see \cite[Theorem 3.29]{brez2011}).
By defintion of $\tw Y$, this implies in turn that 
$\tw{\b X} \otimes \tw{Y}$ is metrizable on $K_n$, for all $n\in\mathbb N$. 
This implies that for every $f\in \S\C(\Xw \times \Yw)$ it holds that
$f_{|K_n}\in \C(K_n, \tw{\b X}\otimes\tw {Y})$ for all $n\in\mathbb N$. Hence, for every open subset $A\subset \mathbb R$
and for every $n\in\mathbb N$, one has that 
$f_{|K_n}^{-1}(A)=f^{-1}(A)\cap K_n \in \tw{\b X}\otimes\tw{Y}$.
Recalling that $K_n=B_n\cap(\b X\times\mathcal A)$, 
by definition of $\tw Y$ this implies that 
$f^{-1}(A)\cap B_n\in \tw{\b X}\otimes\tw{V_1}$ 
for all $n\in\mathbb N$, which is equivalent to 
$f^{-1}(A)\in \tbw{\b X} \otimes \tbw{V_1}$.
Since $f^{-1}(A)=f^{-1}(A)\cap(\b X\times Y)$ by choice of $f$,
by definition of $\tbw Y$ this is equivalent to 
$f^{-1}(A)\in \tbw{\b X}\otimes\tbw{Y}$,
so that $f\in \C(\mathcal \Xbw \times  \Ybw)$.
Conversely, for every 
$f\in \C( \b X_{bw}\times  Y_{bw})$ the same argument
shows that 
$f_{|K_n}\in \C(K_n, \tw{\b X}\otimes\tbw Y)$ for all $n\in\mathbb N$.
Hence, for every \review{sequence} $\{(\b x_k, y_k)\}_{k \in \enne} \subset \b X \times Y$ \review{such that} $(\b x_k,y_k)\to(\b x,y) \in \b X \times Y$ with respect to \review{the weak topology} $\tw{\b X}\otimes\tw Y$,
there exists $N\in\mathbb N$ large enough so that 
$(\bx_k,y_k)$ \review{and} $(\bx,y)$ both lie in $K_N$ for all $k\in\mathbb N$.
Since $f_{|K_N}\in \C(K_N, \tw{\b X}\otimes\tw Y)$,
one has that $f(\bx_k,y_k)\to f(\bx,y)$, so that $f\in \S\C(\Xw \times\Yw)$.
\end{proof}\noindent

\review{
\begin{remark}
When working with the bounded-weak topology, 
by virtue of Lemma \ref{eq_spaces}, it follows immediately that the
sequential weak Feller property is equivalent to the 
Feller property with respect to the bounded-weak topology, namely
$P$ is sequential weak Feller if and only if the operator
$P_t: \C_b(\Xbw \times \Ybw) \rightarrow \C_b(\Xbw \times \Ybw)$ for all $t \ge 0$.
\end{remark}}\noindent
\review{The main result of this subsection is the following.}
\begin{thm}\label{th:Feller}
Let either Assumptions \ref{hyp:structural}--\ref{hyp:diffusionAC} hold with $Z=\bVsd$ or Assumptions \ref{hyp:structural}--\ref{hyp:additionalAC} hold with $Z=\bHs$. Then, the family of operators $P$ satisfies the Feller and the sequential weak Feller properties.
\end{thm}
\begin{proof}
\review{Let us start from the latter property. Recalling Definition \ref{swF}, it amounts to prove the operator $P_t$ maps $\S\C_b(\Xw \times \Yw)$ into itself for all $t \geq 0$. If $t=0$, then the result is trivial.} Let then $t > 0$ and 
$\phi \in \S\C_b(\Xw \times \Yw)$.
Let \review{further} $\{(\b x_k,y_k)\}_{k \in \enne} \subset \b X\times Y$
and $(\bx,y)\in \b X\times Y$ be
such that $(\bx_k,y_k) \to (\bx,y)$ in 
$\Xw \times \Yw$. For convenience of notation, let us \review{also} set $(\bu_k, \varphi_k):= (\bu^{x_k},\varphi^{y_k})$
and $(\bu,\varphi):=(\bu^x,\varphi^y)$.
We have to prove that $(P_t\phi)(x_k,y_k)\to (P_t\phi)(x,y)$, namely that 
\begin{equation}\label{mean-convergence}
\lim_{k\to +\infty}\E\left[ \phi(\bu_k(t), \varphi_k(t))\right]=\E \left[ \phi(\bu(t), \varphi(t))\right].
\end{equation}
\review{By convergence of the sequence $\{(\b x_k,y_k)\}_{k \in \enne}$}, there exists a constant $C>0$ independent of $k$ such that 
\[
\|\b x_k\|_{\bHs} +\|y_k\|_{V_1} \leq C \qquad\forall\,k\in\mathbb N.
\]
Let $T>t$ be arbitrarily fixed. From Proposition~\ref{a_priori_est}, we infer the following uniform estimates 
\begin{equation*}
\sup_{k \in \mathbb{N}}
\|\varphi_k\|_{L^p(\Omega;L^\infty(0, T);V_1)\cap L^2(0,T; V_2)
\cap W^{\gamma,p}(0,T; V_1^*))}  
+\sup_{k\in\mathbb N}\|F'(\varphi_k)\|_{L^p(\Omega; L^2(0,T; H))}
< +\infty
\end{equation*}
and 
\begin{equation*}
\sup_{k \in \mathbb{N}}
\|\bu_k\|_{L^p(\Omega;C^0([0,T];\bHs)\cap L^2(0,T; \bVs)
\cap W^{\gamma,p}(0,T; \b V_\sigma^*))} 
 < +\infty.
\end{equation*}
Therefore, by virtue of Proposition \ref{prop_tight}, 
we deduce that the sequence
$\{\text{Law}_{\mathbb{P}}(\bu_k,\varphi_k)\}_{k\in\mathbb N}$
 is tight in the \review{trajectory} space $\trajU \times \trajF$.  Therefore, recalling that 
 the spaces $\trajU$ and $\trajF$ are locally convex spaces, we can apply the Prokhorov theorem  (see \cite[Theorem II.6.7]{Parthasarathy_1967}) and the generalization to nonmetric spaces of the Skorokhod theorem given by Jakubowski (see \cite{Jak} and \cite{Brz+Ondr_2011}) to infer the existence of a probability space $(\widetilde{\Omega},\widetilde{\mathcal{F}},\widetilde{\mathbb{P}})$ and a family of therein defined $\trajU\times\trajF$-valued random variables $\{\widetilde \bu_{k}, \widetilde\varphi_{k}\}_{k\in\mathbb N}$ 
and $(\widetilde \bu, \widetilde\varphi)$ such that
\begin{equation}\label{skor1}
\text{Law}_{\mathbb{P}}(\bu_{k},\varphi_{k})=
\text{Law}_{\widetilde{\mathbb{P}}}(\widetilde\bu_{k},\widetilde\varphi_{k})
\qquad\forall\,k\in\mathbb N,
\end{equation}
and, along a \review{suitable} subsequence, 
$(\widetilde\bu_{k_j},\widetilde\varphi_{k_j})\to (\widetilde \bu, \widetilde\varphi)$
in $\trajU \times \trajF$, $\widetilde\P$-almost surely.
In particular, since $(\widetilde \bu_{k_j}, \widetilde\varphi_{k_j})\to(\widetilde\bu, \widetilde\varphi)$ in $\C_{\text w}([0,T]; \bHs)\times
\C_{\text w}([0,T]; V_1)$, it holds that 
\[
\left(\widetilde{ \bu}_{k_j}(t), \widetilde {\varphi}_{k_j}(t)\right) \to  \left( \widetilde\bu(t),\widetilde \varphi (t)  \right)\qquad \text{in }
\Xw \times \Yw,\, \widetilde{\mathbb{P}}\mbox{-a.s.}
\]
By repeating the proof of \cite[Section 3.4]{DPGS}, it follows that 
the stochastic system 
$( \widetilde{\Omega} ,  \widetilde{\mathcal{F}}, \widetilde\P ,  \widetilde \varphi,\widetilde  \bu)$ is a martingale solution to system \eqref{AC_NS} with \review{respect to} the initial data $(\bx,y)$,
in the sense of \cite[Definition 2.4]{DPGS}.
Since pathwise uniqueness holds for the system \eqref{AC_NS} 
(cf. Theorem~\ref{th:probstrongsol}), also uniqueness in law holds (see for instance \cite{Rockner08}): this implies both that 
\begin{equation}  \label{skor3}
\widetilde \E [\phi((\widetilde \bu(t),\widetilde{\varphi}(t))] = \E [\phi((\bu(t),\varphi(t))],
\end{equation}
and that the convergences hold along the entire sequence, namely 
\begin{equation}
\label{skor2}
  \left(\widetilde{\bu}_{k}(t), \widetilde {\varphi}_{k}(t)\right) \to  \left( \widetilde\bu(t),\widetilde \varphi (t)  \right)\qquad \text{in }
\Xw \times \Yw,\,\widetilde{\mathbb{P}}\mbox{-a.s.}
\end{equation}
Consequently, recalling that $\phi\in \S\C_b({\Xw \times \Yw})$,
by exploiting the equivalences of laws in \eqref{skor1} 
and \eqref{skor3}, the convergence \eqref{skor2}, 
and the dominated convergence theorem, we deduce that 
\review{\begin{align*}
    \lim_{k\to+\infty}(P_t\phi)(\bx_k,y_k)&=
    \lim_{k\to+\infty}\E [\phi(\bu_k(t),\varphi_k(t))]\\&=\lim_{k\to+\infty}\widetilde\E[\phi(\widetilde\bu_k(t),\widetilde\varphi_k(t))]\\
    &=\widetilde\E[\phi(\widetilde\bu(t),\widetilde\varphi(t))]
    \\& =\E[\phi(\bu(t),\varphi(t))]
    \\& =(P_t\phi)(\bx,y).
\end{align*}}This proves the validity of the sequential weak Feller property. \\
As for the Feller property, it is enough to prove that $P_t$ maps $\C_b(\Xs \times \Ys)$ into itself for every $t \geq 0$. Retaining the same notation from the previous argument, we now assume that the sequence $(\b x_k, y_k) \to (x, y)$ in $\Xs \times \Ys$. Once again, if $t = 0$, the claim follows immediately. Fix then $t > 0$ and $\phi \in \C_b(\Xs \times \Ys)$. Performing the change of probability via the tightness argument illustrated above, it is enough to show that
\[
(\widetilde \bu_k(t), \, \widetilde \varphi_k(t)) \to (\widetilde \bu(t), \, \widetilde \varphi(t)) \quad \text
{ in } \Xs \times \Ys, \, \tP\text{-a.s.} 
\]
to conclude the claim analogously as before. Since the above convergence already holds in $\Xw \times \Yw$, we show convergence of norms. Arguing as in \cite[Subsection 3.3]{DPGS}, an application of the It\^{o} lemma yields the equality
\begin{multline*}
    \dfrac 12 \tE \|\widetilde{\b u}_k(t)\|_{\bHs}^2 + \dfrac{\beta}{2} \tE \|\nabla \widetilde \varphi_k(t)\|^2_{\b H} = \dfrac 12 \|\b x_k\|^2_{\bHs} + \dfrac{\beta}{2}\|\nabla y_k\|^2_{\bH} + \|F(y_k)\|_{L^1(\OO)} - \tE \|F(\widetilde \varphi_k(t))\|_{L^1(\OO)} \\
    - \tE \int_0^t \left( \nu\|\nabla \widetilde {\b u}_k(\tau)\|^2_{\b H} + \|\widetilde \mu_k(\tau)\|^2_H \right) \: \d \tau + \dfrac \beta 2 \tE \int_0^t \sum_{j \in \enne} \int_\OO F''(\widetilde \varphi_k(\tau))|g_j(\widetilde \varphi_k(\tau))|^2 \: \d \tau \\ + \dfrac 12 \tE \int_0^t \|G_1(\widetilde{\b u}_k(\tau))\|^2_{\mathcal L_{\text{HS}}(U_1,\bHs)}\: \d \tau + \dfrac \beta 2 \tE \int_0^t \|\nabla G_2(\widetilde{\varphi}_k(\tau))\|^2_{\mathcal L_{\text{HS}}(U_2,H)}\: \d \tau.
\end{multline*}
Taking superior limits as $k\to+\infty$ and exploiting the weak sequential lower semicontinuity of the norms as well as the known convergence properties and the dominated convergence theorem, we have
\begin{multline*}
    \limsup_{k\to+\infty}  \left[ \dfrac 12 \tE \|\widetilde{\b u}_k(t)\|_{\bHs}^2 + \dfrac{\beta}{2} \tE \|\nabla \widetilde \varphi_k(t)\|^2_{\b H} \right] \leq \dfrac{1}{2} \| \b x \|^2_{\bHs} + \dfrac{\beta}{2}\|\nabla y\|^2_{\bH} + \|F(y)\|_{L^1(\OO)} - \tE \|F(\widetilde \varphi(t))\|_{L^1(\OO)} \\
    - \tE \int_0^t \left( \nu\|\nabla \widetilde {\b u}(\tau)\|^2_{\b H} + \|\widetilde \mu(\tau)\|^2_H \right) \: \d \tau + \dfrac \beta 2 \tE \int_0^t \sum_{j \in\enne } \int_\OO F''(\widetilde \varphi(\tau))|g_j(\widetilde \varphi(\tau))|^2 \: \d \tau \\ + \dfrac 12 \tE \int_0^t \|G_1(\widetilde{\b u}(\tau))\|^2_{\mathcal L_{\text{HS}}(U_1,\bHs)}\: \d \tau + \dfrac \beta 2 \tE \int_0^t \|\nabla G_2(\widetilde{\varphi}(\tau))\|^2_{\mathcal L_{\text{HS}}(U_2,H)}\: \d \tau.
\end{multline*}
However, it is straightforward to observe, by a further a posteriori application of the It\^{o} lemma, that the right hand side equals the limiting quantity of interest, i.e.,
\begin{equation*}
    \limsup_{k\to+\infty}  \left[ \dfrac 12 \tE \|\widetilde{\b u}_k(t)\|_{\bHs}^2 + \dfrac{\beta}{2} \tE \|\nabla \widetilde \varphi_k(t)\|^2_{\b H} \right] =  \left[ \dfrac 12 \tE \|\widetilde{\b u}(t)\|_{\bHs}^2 + \dfrac{\beta}{2} \tE \|\nabla \widetilde \varphi(t)\|^2_{\b H} \right],
\end{equation*}
which in turn implies, multiplying by two, by weak lower semicontinuity
\[
\begin{split}
\limsup_{k\to+\infty}  \tE \left[ \|\widetilde{\b u}_k(t)\|_{\bHs}^2 +  \beta \|\nabla \widetilde \varphi_k(t)\|^2_{\b H} \right] & = \limsup_{k \to +\infty} \tE \, \|(\widetilde{\b u}_k,\, \sqrt \beta \nabla \widetilde \varphi_k)\|^2_{\b X \times \bH} \\ & = \tE \|(\widetilde{\b u},\, \sqrt \beta \nabla \widetilde\varphi)\|^2_{\b X \times \bH} \\
& \leq \liminf_{k\to+\infty} \tE \|(\widetilde{\b u}_k,\, \sqrt \beta \nabla \widetilde \varphi_k)\|^2_{\b X \times \bH}.
\end{split}
\]
This then implies convergence of norms, and in turn, strong convergence (i.e., convergence in the topological space $\Xs \times \Ys$). The Feller property then follows analogously as the weak sequential case. The proof is complete.
\end{proof}

\begin{remark}
    It is important to underline that, for the sake of what follows, it is not necessary for the Feller and the sequential weak Feller properties to hold simultaneously: instead, either one is sufficient.
    However, the Feller scenario is widely represented in the literature, whereas the sequential weak Feller scenario is less studied.
    In order to prove the existence of invariant measures one can exploit the Feller property by appealing to the Krylov--Bougoliubov criterion, or the sequentially weak Feller property by appealing to the Maslowski--Seidler criterion. In what follows we will show the existence of an invariant measure in both ways.
    Furthermore, we stress that in our setting the Feller property 
    does not follow directly from a continuous dependence on
    the initial data, as it is common for classical evolution equations, but has to be proven independently.
\end{remark}

\subsection{The Markov property}
\review{Let us now turn our attention to proving the} Markov property for the solution process $(\bu,\, \varphi)$ to system \eqref{AC_NS}. When the solution is a continuous process taking values in a separable Banach space, endowed with the strong topology, \review{the result follows from the general theory (see, e.g., \cite[Theorem 9.14]{dapratozab}). In weaker settings,} one can argue similarly as in \cite[Proposition 7.3]{BFZ}. \review{Nonetheless, in the case of interest, the arguments have to be carefully revised in order to account for the absence of a linear structure on $Y$.}

\begin{thm}\label{th:Markov}
Let either Assumptions \ref{hyp:structural}--\ref{hyp:diffusionAC} hold with $Z=\bVsd$ or Assumptions \ref{hyp:structural}--\ref{hyp:additionalAC} hold with $Z=\bHs$. Let also $(\b x, y) \in \b X \times Y$ and $\phi\in \S\C_b(\Xw \times \Yw)$. Then
\begin{equation*}
  {\mathbb E}\left[  
  \phi(\bu^\bx(t+s),\, \varphi^y(t+s))
  \,|\,{\mathcal F}_s\right]
  =
  (P_{t}\phi)  
  (\bu^\bx(s),\,\varphi^y(s))  \quad {\mathbb P}\text{-a.s.}
\end{equation*}
for any $t,\,s > 0$. \review{In particular, the family $P$ defines a Markov semigroup on $\mathcal{B}_b(\b X \times Y)$.}
\end{thm}
\begin{proof}
\review{Let $t$ and $s$ be arbitrary positive times.} By the pathwise uniqueness of the solution, we have that
\begin{equation*}
   ( \bu^\bx(t+s),\,\varphi^y(t+s)) = 
   ( \bu^{s,\eta}(t+s),\,\varphi^{s,\zeta}(t+s))  \quad \mathbb{P} \text{-a.s.}
\end{equation*}
where \review{we set} $(\eta, \zeta):=(\bu^\bx(s),\varphi^y(s))$.
Therefore, \review{the claim is equivalent to}
\begin{equation}
\label{Markov_proof}
    {\mathbb E}\left[  \phi( \bu^{s,\eta}(t+s),\,\varphi^{s,\zeta}(t+s))
    \,|\,{\mathcal F}_s\right]
  =
  (P_{t}\phi) (  \eta,\,\zeta)  \qquad {\mathbb P}\text{-a.s.}
\end{equation}
\review{Now, }suppose at first that $\left(  \eta, \zeta\right)$ is a \review{simple} random variable, \review{i.e.,}
\begin{equation*}
\left(  \eta, \zeta\right) = \sum_{j=1}^N (\bx_j, y_j)\pmb{1}_{A_j},
\end{equation*}
\review{for some $N \in \enne$, for some finite family $\{(\bx_j, y_j)\}_{j=1}^N \subset \b X \times Y$ and with $\{A_j\}_{j=1}^N \subset \mathcal F_s$ being a partition of $\Omega$.} Arguing \review{as in} \cite[Theorem 9.14]{dapratozab},
again by pathwise uniqueness one has that
\begin{equation*}
    \phi( \bu^{s,\eta}(t+s),\varphi^{s,\zeta}(t+s)) = 
    \sum_{j=1}^N  \phi( \bu^{s,\bx_j}(t+s),\varphi^{s,y_j}(t+s))\pmb{1}_{A_j} \quad \mathbb{P}\text{-a.s.}
\end{equation*}
Using now the fact that $\pmb{1}_{A_j}$ are $\mathcal{F}_s$-measurable,
that the random variables $(\bu^{s,\bx_j}(t+s), \varphi^{s,y_j}(t+s))$ 
are \review{all} independent of $\mathcal{F}_s$, 
and that $(\bu^{s,\bx_j}(t+s), \varphi^{s,y_j}(t+s))=(\bu^{\bx_j}(t),\varphi^{y_j}(t))$ \review{since the system is autonomous},
we obtain
\begin{align*}
\E\left[
\phi( \bu^{s,\eta}(t+s),\varphi^{s,\zeta}(t+s))\,|\,{\mathcal F}_s
\right]&=
\sum_{j=1}^N \mathbb{E}\left[
\phi( \bu^{s,\bx_j}(t+s),\,\varphi^{s,y_j}(t+s))\pmb{1}_{A_j}\,|\,\mathcal{F}_s\right]\\
&=\sum_{j=1}^N \pmb{1}_{A_j} \mathbb{E}
\left[\phi( \bu^{s,\bx_j}(t+s),\,\varphi^{s,y_j}(t+s))\right]
\\
&=\sum_{j=1}^N \pmb{1}_{A_j} 
\mathbb{E}\left[\phi( \bu^{\bx_j}(t),\,\varphi^{y_j}(t))\right] 
= (P_t\phi)(\eta, \zeta),
\end{align*}
and this proves \eqref{Markov_proof} for simple random variables. Let now $(\eta, \zeta)$ \review{be} a general $\b X \times Y$-valued $\mathcal{F}_s$-measurable random variable such that $\eta\in L^p(\Omega,\mathcal F_s; \bH)$ and $\zeta\in L^p(\Omega,\mathcal F_s; V_1)$ for some fixed $p>2$. \review{By density,} there exists a sequence 
$\{(\eta_k, \zeta_k)\}_{k \in \mathbb{N}}$ of simple $\b X\times Y$-valued
random variables as above
such that $(\eta_k,\zeta_k)\to(\eta,\zeta)$ in $L^p(\Omega; \bHs\times V_1)$
and in $\Xs \times \Ys$, $\P$-almost surely. For every $k\in\mathbb N$, we know that 
\begin{equation}
\label{Markov_app}
   {\mathbb E}\left[  \phi( \bu^{s,\eta_k}(t+s),\varphi^{s,\zeta_k}(t+s))
    |{\mathcal F}_s\right]
  =(P_{t}\phi) (  \eta_k, \zeta_k)  \quad {\mathbb P}\text{-a.s.}
\end{equation}
\review{By the sequential weak Feller property (see Theorem \ref{th:Feller}), } $P_t\phi \in \S\C_b(\Xw \times \Yw)$: hence, the right-hand side of \eqref{Markov_app} satisfies
\[
  (P_{t}\phi) (  \eta_k, \zeta_k)\to(P_{t}\phi) (  \eta, \zeta)
  \quad {\mathbb P}\text{-a.s.}
\]
Since
\[
  \sup_{k\in\mathbb N}\|\eta_k\|_{L^p(\Omega; \bHs)}
  +\sup_{k\in\mathbb N}\|\zeta\|_{L^p(\Omega; V_1)}<+\infty,
\]
one can argue exactly as in the proof of Theorem~\ref{th:Feller}
(set on the time interval $[s,t+s]$) to deduce in particular that
\[
  ( \bu^{s,\eta_k}(t+s),\varphi^{s,\zeta_k}(t+s))\to
  ( \bu^{s,\eta}(t+s),\varphi^{s,\zeta}(t+s)) \quad\text{in } \mathcal \Xw \times \Yw, \quad\P\text{-a.s.}
\]
Hence, thanks to the dominated convergence theorem for the conditional expectation
and the fact that $\phi \in \S\C_b(\Xw \times \Yw)$, 
we obtain also that 
\[
{\mathbb E}\left[  \phi( \bu^{s,\eta_k}(t+s),\varphi^{s,\zeta_k}(t+s))
    \,|\,{\mathcal F}_s\right] \to 
    {\mathbb E}\left[  \phi( \bu^{s,\eta}(t+s),\varphi^{s,\zeta}(t+s))
    \,|\,{\mathcal F}_s\right] \quad\P\text{-a.s.}
\]
Hence, the claim follows by letting $k\to+\infty$ in \eqref{Markov_app}. \review{By exploiting} the properties of the conditional expectation, 
we directly deduce that the family of operators 
$P=\{P_t\}_{t\ge0}$ is a Markov semigroup on $\mathcal B_b(\b X\times Y)$,
namely that $P_{t+s}=P_tP_s$ for all $s,\,t\ge0$.
\end{proof}

\subsection{Existence of invariant measures}
\label{MS_sec}
\review{This subsection is devoted to showing that} the family of operators $P$ given in \eqref{P_t} admits at least an invariant measure. \review{
\begin{defin} Let either Assumptions \ref{hyp:structural}--\ref{hyp:diffusionAC} hold with $Z=\bVsd$ or Assumptions \ref{hyp:structural}--\ref{hyp:additionalAC} hold with $Z=\bHs$.
An invariant measure for the Markov semigroup $P$ is a probability measure $\mu\in\cP(\b X \times Y)$ such that 
  \[
  \int_{\b X\times Y} \phi(\bx,y)\,\mu(\d \bx, \d y) = 
  \int_{\b X\times Y} P_t\phi(\bx,y)\,\mu(\d \bx, \d y) \quad\forall\,t\geq0,\quad\forall\,\phi \in \B_b(\b X \times Y).
  \]
\end{defin}}\noindent
The main result of this \review{sub}section is the following.
\begin{thm}\label{th:ex_inv}
Let either Assumptions \ref{hyp:structural}--\ref{hyp:diffusionAC} hold with $Z=\bVsd$ or Assumptions \ref{hyp:structural}--\ref{hyp:additionalAC} hold with $Z=\bHs$. Then, there exists at least an invariant measure for the semigroup $P$.
\end{thm}
\begin{proof}
For the sake of completeness, we provide two possible proofs, by exploiting either the sequential weka Feller property or the Feller property.
\paragraph{\review{\textit{The sequential weak Feller case.}}} The main idea is to use a version of the Krylov-Bogoliubov method \review{involving} the weak topology, as proposed by Maslowski and Seidler \review{in} \cite{MS}, in combination with a tightness result, which is a consequence of Proposition~\ref{a_priori_est}. Recalling that  $\cL_t((0,0),\cdot)$ is
the law of the random variable $(\bu^{\b 0}(t), \varphi^0(t))$ on $\b X \times Y$, since the mapping $(\omega, t) \mapsto (\bu^{\b 0}, \varphi^0)(\omega,t)$ 
is jointly measurable, we can integrate with respect to both variables and define the probability measure
 \begin{equation}
 \label{mu_n}
     \mu_n(A) := \frac 1n \int_0^n \cL_t((\b 0,0), A)\, {\rm d}t , \qquad A \in \cB(\b X \times Y).
 \end{equation}
 on $(\b X \times Y,\, \cB(\b X \times Y))$, for all positive $n\in\mathbb N$.
\review{Our first aim is showing that the family of measures} $\{\mu_n\}_{n \in \enne_+}$ is tight on $\Xbw \times \Ybw$.
To this end, 
given $R>0$, let $B_R:=\{(x, y) \in \bHs \times V_1: \  \|(x, y)\|_{\bHs \times V_1}\le R\}$ and define the set 
\begin{equation}
\label{B_R}
\K_R:=B_R \cap (\bHs\times\mathcal{A})\subset \b X\times Y.
\end{equation}
\review{Let} $\K_R^C:=(\b X\times Y)\setminus \K_R$ denote its complement set in $\b X\times Y$. 
Since $B_R$ is a closed bounded set in $\bHs\times V_1$, 
by reflexivity of and by definition of the topologies on $\b X \times Y$
it follows that $\K_R$ is compact in $\Xw \times \Yw$ and in $\Xbw \times \Ybw$.
Bearing in mind the definition \eqref{mu_n}, by means of the Chebychev inequality and \review{Proposition} \ref{a_priori_est}, for all $n>0$ and $R>0$ we infer that
\begin{align*}
    \mu_n(\K_R^C) 
    &=\frac 1n \int_0^n \mathbb{P}\left(\|(\bu^{\b 0}(t), \varphi^0(t))\|_{\bHs\times  V_1}>R\right)\, {\rm d}t
    \le \frac{1}{nR^2}\int_{0}^n \mathbb{E} \left[ 
    \|(\bu^{\b 0}(t)\|^2_{\bHs}+ \|\varphi^0(t))\|_{V_1}^2
    \right]\, {\rm d}t \leq \frac C{R^2},
\end{align*}
for some constant $C>0$ independent of $n$. This direclty
implies the tightness of 
$\{\mu_n\}_{n \in \enne_+}$ on
$\Xbw \times \Ybw$. \review{Next, we} exploit tightness to extract a subsequence \review{of measures weakly}
converging to some probability measure $\mu$,
which will be our candidate invariant measure for $P$. \review{In the non-metric setting, the Prokhorov--Jakubowski Theorem (see \cite[Theorem~3]{Jak})} requires that the space $\Xbw \times \Ybw$ is countably separated, that is, \review{that} there exists a countable family $\{f_i:\b X \times Y\rightarrow [-1,1] \}_{i \in \mathbb{N}}$ with $f_i \in \C(\Xbw \times \Ybw)$ for all $i \in \mathbb{N}$ \review{such that} for every \review{pair of distinct points} $(\b x_1,y_2),\, (\b x_2,y_2)\in \b X\times Y$ there is $i\in\mathbb N$
such that $f_i(\b x_1,y_1)\neq f_i(\b x_2,y_2)$. \review{In our setting, since} the dual space $(\bHs\times V_1)^*$ is separable, there exists a countable dense family $\{(\b h_i, l_i)\}_{i \in \enne} \subset 
(\bHs\times V_1)^*$. Thus, for any two \review{distinct} elements $(\b x_1,y_2), (\b x_2,y_2)\in \b X\times Y$, one clearly has that 
$(\b x_1,y_1)\neq (\b x_2,y_2)$ as elements in $\bHs\times V_1$.
Hence, as a consequence of the Hahn--Banach theorem, 
there exists \review{some index $i \in \enne$} such that 
\[\langle (\b h_i,l_i),(\b x_1,y_1)\rangle_{
(\bHs \times V_1)^*,\bHs \times V_1} \ne 
\langle ((\b h_i,l_i), (\b x_2,y_2)\rangle_{
(\bHs \times V_1)^*,\bHs \times V_1}. \] 
\review{For any $i \in \enne$, define}
\[
f_i:\b X \times Y\rightarrow [-1,1] \qquad 
(\b x,y) \mapsto \tanh(\langle (\b h_i,l_i),\,(\b x,y)\rangle)_{
(\bHs \times V_1)^*,\bHs \times V_1}.
\]
\review{It is straightforward to check that $f_i \in \S\C(\Xw \times \Yw)=\S\C(\Xbw \times \Ybw)$ for all $i \in \enne$ and that, in turn, $\b X\times Y$ is countably separated}.
Therefore, by the Prokhorov--Jakubowski Theorem 
we infer the existence of a subsequence $\{\mu_{n_k}\}_{k\in \enne}$ and a probability measure $\mu$ on $(\b X \times Y, \cB(\b X \times Y)$ such that $\mu_{n_k}$ converges narrowly to $\mu$ as $k\rightarrow +\infty$, that is 
\[
    \int_{\b X \times Y} \phi \,{\rm d}\mu_{n_k} \rightarrow 
    \int_{\b X \times Y} \phi \, {\rm d}\mu \qquad \forall\,\phi \in \C_b(\Xbw \times \Ybw).
\]
Now, given $t\geq0$ and $\phi\in \C_b(\Xbw \times \Ybw)$, we write
\begin{align*}
    \int_{\b X \times Y} P_t \phi\, \d \mu_{n_k} 
    &=\frac1{n_k}\int_0^{n_k}\int_{\b X \times Y}P_t\phi\,\d\cL_{s}((\b 0,0), \cdot)\,\d s\\
    &=\frac1{n_k}\int_0^{n_k}(P_{t+s}\phi)(\b 0,0)\,\d s
    =\frac1{n_k}\int_t^{t+n_k}(P_{r}\phi)(\b 0,0)\,\d r\\
    &=\int_{\b X \times Y} \phi\, \d\mu_{n_k}
    +\frac1{n_k}\int_{n_k}^{t+n_k}(P_{r}\phi)(\b 0,0)\,\d r
    -\frac1{n_k}\int_0^{t}(P_{r}\phi)(\b 0,0)\,\d r,
    \end{align*}
where
\[
  \left|\frac1{n_k}\int_{n_k}^{t+n_k}(P_{r}\phi)(\b 0,0)\,\d r
    -\frac1{n_k}\int_0^{t}(P_{r}\phi)(\b 0,0)\,\d r\right|
  \leq \frac{2t}{n_k}\|\phi\|_{\infty}\rightarrow0.
\]
From the sequential weak Feller property, proved in Theorem~\ref{th:Feller}, and Lemma \ref{eq_spaces}, 
we know that $P_t \phi \in \C_b(\Xbw \times \Ybw)$ if $\phi \in \C_b(\Xbw \times \Ybw)$. Hence, by letting $k\to+\infty$ we get 
\[  
    \int_{\b X \times Y} P_t \phi\, \d \mu
    =\int_{\b X \times Y} \phi\, \d \mu
\]
for all $\phi \in \C_b(\Xbw \times \Ybw)$. \review{Since, by density arguments, the set $\C_b(\Xbw \times \Ybw)=\S\C_b(\Xw \times \Yw)$ is a determining set for the measure $\mu$ (see e.g. \cite{BFZ}), we infer that $\mu$ in an invariant measure for $P$.}
\review{\paragraph{\textit{The Feller case.}}
In this case, the proof is based on the Krylov--Bougoliubov theorem (see \cite[Theorem 3.3]{SZ}). Consider the family of probability measures} $\{\mu_t\}_{t > 0} \subset \cP(\b X \times Y)$ given by 
\[
\mu_t(A):=\frac 1t \int_0^t \cL_s(( \b 0,0),A)\, {\rm d}s, \qquad A \in \cB(\b X \times Y),\,\, t>0.
\]
Let now $B_n := \{ (\bx,y) \in \bVs \times V_2 \ : \ \|(\bx,y)\|_{\bVs \times V_2} \le n\}$ i.e., the closed ball of radius $n$ in $\bVs \times V_2$, and let $\K_n:= B_n \cap (\bVs \times \mathcal{A}) \subset \bVs \times (\mathcal{A} \cap V_2)$. Then $\K_n$ is a compact \review{sub}set of $\b X_s \times Y_s$ since the embedding $\bVs \times V_2 \subset \bHs \times V_1$ is compact. Hence, Proposition \ref{a_priori_est} and the Chebychev inequality yield, for any $t > 0$,
\begin{multline*}
\mu_t(\K_n^C)
=\frac1t\int_0^t \cL((\b 0,0),\K_n^C)\, {\rm d}s 
=\frac 1t \int_0^t \mathbb{P}\left(\|(\bu^{\b 0}(s), \varphi^0(s))\|_{\bVs \times V_2} > n\right)\,{\rm d}s  
\\
\le \frac{1}{tn^2}\int_0^t \mathbb{E} \left[\|\bu^{\b 0}(s)\|^2_{\bVs}+ \|\varphi^0(s)\|^2_{V_2}\right]\,{\rm d}s 
\le 
\frac{C}{n^2}\left(1 + \frac 1t \right),
\end{multline*}
for some positive constant $C$ only depending on the structural parameters and independent of $t$. \review{As before, it follows that the family of probability measures $\{\mu_t\}_{t>0}$ is tight, and hence the thesis follows from the Krylov--Bougoliubov theorem and the fact that the set $\C_b(\Xs \times \Ys)$ is a determining set for the measure $\mu$.}
\end{proof}
\section{Support of the invariant measures and ergodicity}
\label{sec:erg}
The aim of this section is to show some qualitative properties of the
invariant measures, specifically concerning their moments and their support.
This will allow to prove the existence of an ergodic invariant measure. In this direction, we introduce the spaces
\begin{align}
\label{new_spaces'}
    \mathcal A_s&:=
    \left\{y\in \mathcal A:\; \Psi_s(y) \in L^1(\OO)\right\}, \quad
    s\in[1,+\infty),\\
    \label{new_spaces''}
    \mathcal A_\infty
    &:=\left\{
    y\in \mathcal A: \;\exists\,\delta>0:\;
    |y| \leq 1 - \delta \text{ a.e.~in } \OO
    \right\}=\left\{
    y\in L^\infty(\OO): \;\|y\|_{L^\infty(\OO)}<1
    \right\},\\
    \label{new_spaces}
    Y_s&:=V_1\cap\mathcal A_s, \quad s\in[1,+\infty].
\end{align}
By exploiting the continuity of $\Psi_s$,
it is not difficult to show that $\mathcal A_s, \mathcal A_\infty\in\mathscr B(Y)$ for all $s\geq1$.
Furthermore, it holds that $\mathcal A_{s_1}\subset\mathcal A_{s_2}$ for every $s_1,s_2\in[1,+\infty]$ with 
$s_1\geq s_2$.
With this notation, we define 
\begin{equation}
\label{new_spaces_reg}
  \b X_{reg}:=\bVs\,, \qquad Y_{reg}:=V_2\cap \mathcal A_\infty,
\end{equation}
and note that
$\b X_{reg}\in\mathscr B(\b X)$ and 
$Y_{reg}\in\mathscr B(Y)$. Moreover, 
as a consequence of the the H\"older continuity of 
elements in $V_2$, 
by following the 
arguments of \cite[Section~4]{OS23} we point out that 
$V_2 \cap \mathcal A_{s} \subset \mathcal A_\infty$
for all $s>2$,
which implies the characterisation of $Y_{reg}$ as
\[
  Y_{reg}= V_2 \cap \mathcal A_s \quad\forall\,s>2.
\]

We also recall here the definition of ergodicity
for the the transition semigroup $P$.
First, note that if $\mu$ is an invariant measure for $P$,
then the semigroup $P$ can be extended to 
$L^p(\b X\times Y,\, \mathscr B(\b X\times Y),\, \mu)$:
this can be shown, for instance, by 
density and by definition of invariance.
\begin{defin}
  \label{def:erg}
  An invariant measure $\mu$ for the semigroup $P$ is ergodic if
  \[
  \lim_{t \rightarrow +\infty}\frac 1t \int_0^tP_s \phi\, {\rm d}s 
  = \int_{\b X\times Y}\phi \, \d\mu \qquad \text{in } 
  L^2(\b X\times Y,\, \mu)\quad \forall\, \phi \in L^2(\b X\times Y,\, \mu).
  \]
\end{defin}

Our first result concerns the support of invariant measures. Let $\Pi\subset \mathscr P(\b X\times Y)$ be the 
(nonempty) set of invariant measures for the Markov semigroup $P$.

\begin{thm}
\label{th:support} 
Let either Assumptions \ref{hyp:structural}--\ref{hyp:diffusionAC} hold with $Z=\bVsd$ or Assumptions \ref{hyp:structural}--\ref{hyp:additionalAC} hold with $Z=\bHs$. Then,
there exists a constant $C>0$,
only depending on $\nu$, $\beta$, $L_{G_1}$, $L_{G_2}$, $|\OO|$,
$L_F$, $\|F\|_{\mathcal C([-1,1])}$, such that
\begin{equation}
\label{eq:support}
\sup_{\mu\in\Pi}\int_{\b X\times Y}
\left(\|\bx\|_{\bVs}^2 + \|y\|_{V_2}^2 + \|F'(y)\|_{H}^2\right)\, 
\mu({\rm d}\b x, {\rm d} y)\le C\,.
\end{equation}
Moreover, if Assumptions \ref{hyp:structural}--\ref{hyp:additionalAC} hold with $Z=\bHs$, then for every $s\geq1$
there exists a constant $C>0$,
only depending on $s, L_{G_1}$, $L_{G_2}$, $|\OO|$,
$L_F$, $\|F\|_{\mathcal C([-1,1])}$, such that
\begin{equation}
\label{eq:support2}
\sup_{\mu\in\Pi}\int_{\b X\times Y}
\norm{\Psi_s(y)}_{L^1(\OO)}\, 
\mu({\rm d}\b x, {\rm d} y)\le C,
\end{equation}
and
every invariant measure $\mu$ is supported in $\b X_{reg}\times Y_{reg}$, 
i.e.~$\mu(\b X_{reg}\times Y_{reg})=1$.
\end{thm}
\begin{proof}
The last assertion of the theorem is an immediate consequence of the estimates
\eqref{eq:support}--\eqref{eq:support2}
and the fact that 
$Y_{reg}=V_2\cap \mathcal A_s$ for any $s>2$.
Let us prove now \eqref{eq:support}, so 
let $\mu$ be an arbitrary invariant measure for the transition semigroup $P$. \review{Let us stress that the symbol $C$ may denote any generic constant depending only on the structural assumptions, but not on $\mu$, whose value may be updated from line to line. For the sake of clarity, we split the argument into steps.} 
\paragraph{\textit{Step 1.}} We first show that there exists $C>0$, independent of $\mu$, such that 
\begin{equation}
    \label{supp1}
    \int_{\b X\times Y} \|y\|_{V_1}^2\,\mu(\d \b x, \d y)\leq C.
\end{equation}
Since $Y\subset \mathcal A$, one trivially has
\[
\int_{\b X\times Y} \|y\|_H^2\,\mu(\d \b x, \d y) \leq |\OO|\,.
\]
Now, by arguing on the lines of \eqref{stima1} as in the proof of Proposition~\ref{a_priori_est},
one obtains, after taking expectations, the following estimate 
\begin{align*}
\E\int_0^1
\beta\|\nabla \varphi^y(s)\|_{\b H}^2 \,\d s
\leq \frac12\|y\|_H^2+
2L_F|\OO| + \frac12L_{G_2}|\OO|\,.
\end{align*}
Let us define now, for every $n\in\mathbb N$,
the functions $\Upsilon_n:\b X\times Y\to[0,+\infty)$ as
\begin{equation*}
\Upsilon_n(\b x,y):=\|\nabla y\|^2_{\b H}\pmb{1}_{\{\|\nabla y\|_{\b H}\leq n\}}
+n^2\pmb{1}_{\{\|\nabla y\|_{\b H}> n\}}, \quad (\b x,y)\in \b X\times Y.
\end{equation*}
It is not difficult to check that 
$\Upsilon_n\in\mathcal{B}_b(\b X\times Y)$ for every $n\in\enne$.
Hence, thanks to 
the invariance of $\mu$, the boundedness of $\Upsilon_n$, the definition of $P$,
and the Fubini-Tonelli Theorem we have that
\begin{align}
\notag
&\int_{\b X\times Y} \Upsilon_n(\bx,y)\, \mu(\d \bx, \d y)= 
\int_0^1 \int_{\b X\times Y} \Upsilon_n(\bx,y)\, \mu(\d \bx, \d y)\, {\rm d}s=
\int_0^1 \int_{\b X\times Y} (P_s\Upsilon_n)(\bx,y)\, \mu(\d \bx, \d y)\, {\rm d}s\\
\notag
&= \int_0^1 \int_{\b X\times Y}
\E\left[\Upsilon_n(\bu^{\bx}(s), \varphi^y(s))\right]\, \mu(\d \bx, \d y)\, {\rm d}s
\leq\int_{\b X\times Y}\E\int_0^1
\|\nabla \varphi^y(s)\|_{\b H}^2 \,\d s\, \mu(\d\bx, \d y)\\
&\leq\frac{1}{2\beta}\int_{\b X\times Y}\|y\|^2_H\, \mu(\d \bx, \d y)
+{\frac{2L_F}{\beta}}|\OO| + {\frac{1}{2\beta}L_{G_2}|\OO|}\leq
{\frac{|\OO|}{2\beta}\left(1+4L_F+L_{G_2}\right)}.
\label{aux_supp}
\end{align}
Hence, by letting $n\to+\infty$, by the monotone convergence theorem 
we infer \eqref{supp1}.
\paragraph{\textit{Step 2.}}
We show now that there exists $C>0$, independent of $\mu$, such that 
\begin{equation}
    \label{supp2}
    \int_{\b X\times Y} \|\bx\|_{\bHs}^2\,\mu(\d \b x, \d y)\leq C.
\end{equation}
To this end, for $(\b x,y)\in \b X\times Y$, the idea is 
apply the It\^o formula to 
\[
\mathcal E_\varepsilon(\bu^{\b x}, \varphi^y):=
\chi_\varepsilon\left(\mathcal E(\bu^{\b x}, \varphi^y)\right):=
\chi_\varepsilon\left(\frac12\|\bu^{\b x}\|^2_{\bHs}+
\frac\beta2\|\nabla\varphi^y\|^2_{\b H}+
\|F(\varphi^y)\|_{L^1(\OO)}\right),
\]
where for all $\varepsilon>0$ the function
  $\chi_\varepsilon \in C^2_b(\erre_+)$ is defined as $\chi_\varepsilon(r)=\frac{r}{1+\varepsilon r}$, $r\geq0$, so that $\chi'_\varepsilon(r)=\frac{1}{(1+\varepsilon r)^2}$ and $\chi''_\varepsilon(r)=-\frac{2\varepsilon}{(1+\varepsilon r)^3}$
for all $r\geq0$. 
By arguing on the lines of \eqref{stima6} (possibly employing 
suitable Yosida-type approximations in order to rigorously justify the arguments)
one obtains, after taking expectations, 
\begin{align*}
&\E\mathcal E_\varepsilon(\bu^{\b x}(t), \varphi^y(t))+
\E\int_0^t
\chi'_\varepsilon\left(
\mathcal E(\bu^{\b x}(s), \varphi^y(s))\right)
\left[\|\nabla\bu^{\b x}(s)\|^2_{\bH}
+\beta^2\|\Delta \varphi^y(s)\|^2_H
+\|F'(\varphi^y(s))\|_H^2
\right] \,\d s
\\
&\qquad+2\beta \E\int_0^t
\chi'_\varepsilon\left(\mathcal E(\bu^{\b x}(s), \varphi^y(s))\right)\int_{\OO} 
F^{\prime \prime}(\varphi^y(s))|\nabla \varphi^y(s)|^2 \, {\rm d}s\\
&\leq \mathcal E_\varepsilon({\b x}, y)
+\frac{1}{2}\E\int_0^t 
\chi'_\varepsilon\left(
\mathcal E(\bu^{\b x}(s), \varphi^y(s))\right)
\|G_{1}(\b{u}^{\b x}(s))\|^2_{\LL_{HS}(U_1, \bHs)} \, \mathrm{d}s\\
&\qquad+ 
\frac 12 \E\int_0^t 
\chi'_\varepsilon\left(\mathcal E(\bu^{\b x}(s), \varphi^y(s))\right)
\sum_{k\in\enne}\int_{\OO}\left[
|g_k'(\varphi^y(s))\nabla \varphi^y(s)|^2+
F''(\varphi^y(s))|g_k(\varphi^y(s))|^2\right]\,\d s\\
&\qquad+\frac12\E\int_0^t
\chi''_\varepsilon\left(\mathcal E(\bu^{\b x}(s), \varphi^y(s))\right)
\|(\bu^{\b x}(s), G_1(\bu^{\b x}(s)))_{\bHs}\|^2_{\LL_{HS}(U_1,\mathbb R)}\,\d s\\
&\qquad+\frac12\E\int_0^t
\chi''_\varepsilon\left(\mathcal E(\bu^{\b x}(s), \varphi^y(s))\right)
\|(-\beta\Delta\varphi^y(s)+F'(\varphi^y(s)), G_2(\varphi^y(s)))_{H}\|^2_{\LL_{HS}(U_2,\mathbb R)}\,\d s.
\end{align*}
By exploiting the fact that $0<\chi_\varepsilon'<1$ and $\chi_\varepsilon''\leq0$, arguing as in the proof of 
Proposition~\ref{a_priori_est}, we deduce that 
\begin{align*}
&\E\mathcal E_\varepsilon(\bu^{\b x}(t), \varphi^y(t))+
\E\int_0^t
\chi'_\varepsilon\left(\mathcal E(\bu^{\b x}(s), \varphi^y(s))\right)
\|\b{u}^\bx(s)\|^2_{\bVs}\,\d s\\
&\leq \mathcal E_\varepsilon(\bx, y)
+{\beta\left(2L_F+\frac{L_{G_2}}{2}\right)}\E\int_0^t
\chi'_\varepsilon\left(\mathcal E(\bu^{\b x}(s), \varphi^y(s))\right)
\|\nabla\varphi^y(s)\|_{\b H}^2\,\d s
+\frac12(L_{G_1}+|\OO|L_{G_2})t.
\end{align*}
Note now that $\mathcal E_\varepsilon\in\mathcal B_b(\b X\times Y)$
and $\E[\mathcal E_\varepsilon(\bu^{\b x}(t), \varphi^y(t))]=(P_t \mathcal E_\varepsilon)(\b x,y)$.
Hence, denoting by $k_0$ the norm of the continuous inclusion 
$\bVs\embed \bHs$, by integrating with respect to $\mu$ one gets for $t=1$ that 
\begin{align*}
&\frac1{k_0}
\int_{\b X\times Y}\E\int_0^1
\chi'_\varepsilon\left(\mathcal E(\bu^{\b x}(s), \varphi^y(s))\right)
\|\b{u}^{\b x}(s)\|^2_{\bHs}\,\d s\,\mu(\d \b x, \d y)\\
&\leq
{\beta\left(2L_F+\frac{L_{G_2}}{2}\right)}
\int_{\b X\times Y}
\E\int_0^1
\chi'_\varepsilon\left(\mathcal E(\bu^{\b x}(s), \varphi^y(s))\right)\|\nabla\varphi^y(s)\|_{\b H}^2\,\d s\,\mu(\d \b x, \d y)
+\frac12(L_{G_1}+|\OO|L_{G_2}).
\end{align*}
Now, setting $\mathcal U_\varepsilon(\bx,y):=
\chi'_\varepsilon(\mathcal E(\b x,y))\|\bx\|_{\bHs}^2$
and $\mathcal Y_\varepsilon(\b x,y):=
\chi'_\varepsilon(\mathcal E(\b x,y))
\|\nabla y\|_{\b H}^2$,
for $(\b x,y)\in \b X\times Y$,
one has that $\mathcal U_\varepsilon, \mathcal Y_\varepsilon \in\mathcal B_b(\b X\times Y)$, 
so that the Fubini-Tonelli theorem and the invariance of $\mu$ yield
\begin{align*}
&\frac1{k_0}\int_{\b X\times Y}
\mathcal U_\varepsilon(\b x,y)\,\mu(\d \b x, \d y)
=\frac1{k_0}\int_{\b X\times Y}\int_0^1
(P_s\mathcal U_\varepsilon)(\b x,y)\,\d s\,\mu(\d\b x, \d y)\\
&=\frac1{k_0}
\int_{\b X\times Y}\E\int_0^1
\mathcal U_\varepsilon(\bu^{\b x}(s), \varphi^y(s))\,\d s\,\mu(\d\b x, \d y)\\
&\leq
{\beta\left(2L_F+\frac{L_{G_2}}{2}\right)}
\int_{\b X\times Y}
\E\int_0^1
\mathcal Y_\varepsilon(\bu^{\b x}(s), \varphi^y(s))\,\d s\,\mu(\d\b x, \d y)
+\frac12(L_{G_1}+|\OO|L_{G_2})\\
&{\beta\left(2L_F+\frac{L_{G_2}}{2}\right)}
\int_{\b X\times Y}
\int_0^1
(P_s\mathcal Y_\varepsilon)(\b x,y)\,\d s\,\mu(\d\b x, \d y)
+\frac12(L_{G_1}+|\OO|L_{G_2})\\
&{\beta\left(2L_F+\frac{L_{G_2}}{2}\right)}
\int_{\b X\times Y}
\mathcal Y_\varepsilon(\b x,y)\,\mu(\d\b x, \d y)
+\frac12(L_{G_1}+|\OO|L_{G_2}).
\end{align*}
By letting now $\varepsilon\searrow0$, thanks to the dominated convergence theorem,
we infer that 
\begin{align*}
\frac1{k_0}\int_{\b X\times Y}
\|\bx\|_{\bHs}^2\,\mu(\d\b x, \d y)
\leq
{\beta\left(2L_F+\frac{L_{G_2}}{2}\right)}
\int_{\b X\times Y}
\|\nabla y\|_{\b H}^2\,\mu(\d\b x, \d y)
+\frac12(L_{G_1}+|\OO|L_{G_2}),
\end{align*}
so that \eqref{supp2} follows from the estimate \eqref{supp1}.
\paragraph{\textit{Step 3.}}
We show now that there exists $C>0$, independent of $\mu$, such that 
\begin{equation}
    \label{supp3}
    \int_{\b X\times Y} \left(
    \|\bx\|_{\bVs}^2+
    \|\Delta y\|_{H}^2+
    \|F'(y)\|_H^2\right)\,\mu(\d\b x, \d y)\leq C.
\end{equation}
By arguing on the lines of \eqref{stima6} as in the proof of Proposition~\ref{a_priori_est},
one obtains, after taking expectations, 
\begin{align*}
&\E\int_0^1
\left({\nu}\|\bu^\bx(s)\|_{\bVs}^2+
\beta^2\|\Delta \varphi^y(s)\|_{H}^2
+\|F'(\varphi^y(s))\|_H^2\right)\,\d s\\
&\leq \frac12\|\bx\|_{\bHs}^2
+{\frac{\beta}{2}}\|\nabla y\|_{\b H}^2+
{\|F(y)\|_{L^1(\OO)}}+\frac12(L_{G_1}+|\OO|L_{G_2})
+{\beta\left(2L_F+\frac{L_{G_2}}{2}\right)}\E\int_0^1\|\nabla\varphi^y(s)\|_{\b H}^2\,\d s.
\end{align*}
Define now, for every $n\in\mathbb N$,
the functions $\widetilde\Upsilon_n:\b X\times Y\to[0,+\infty)$ as
\begin{equation*}
\widetilde\Upsilon_n(\b x,y):=
\begin{cases}
\|\bx\|_{\bVs}^2+\beta^2\|\Delta y\|^2_{H}+\|F'(y)\|_H^2 \quad&\text{if }
\|\bx\|_{\bVs}^2+\beta^2\|\Delta y\|^2_{H}+\|F'(y)\|_H^2\leq n,\\
n \quad&\text{otherwise}.
\end{cases}
\end{equation*}
It is not difficult to check that 
$\widetilde\Upsilon_n\in\mathcal{B}_b(\b X\times Y)$ for every $n\in\enne$.
Hence, thanks to 
the invariance of $\mu$, the boundedness of $\widetilde\Upsilon_n$, the definition of $P$,
and the Fubini-Tonelli Theorem we have that
\begin{align*}
&\int_{\b X\times Y} \widetilde\Upsilon_n(\bx,y)\, \mu(\d\bx, \d y)
= 
\int_0^1 \int_{\b X\times Y} \widetilde\Upsilon_n(\bx,y)\, \mu(\d\bx, \d y)\, {\rm d}s
=
\int_0^1 \int_{\b X\times Y} (P_s\widetilde\Upsilon_n)(\bx,y)\, \mu(\d\bx, \d y)\, {\rm d}s\\
&= \int_0^1 \int_{\b X\times Y}
\E\left[\widetilde\Upsilon_n(\bu^\bx(s), \varphi^y(s))\right]\, \mu(\d\bx, \d y)\, {\rm d}s\\
&{\le}\int_{\b X\times Y}\E\int_0^1
\left(\|\bu^\bx(s)\|_{\bVs}^2+
\beta^2\|\Delta \varphi^y(s)\|_{H}^2
+\|F'(\varphi^y(s))\|_H^2\right)\,\d s\,\mu(\d\bx,\d y)\\
&\leq\frac12\int_{\b X\times Y}
\left(\|\bx\|^2_\bHs+{\beta}\|\nabla y\|_{\b H}^2 {+\|F(y)\|_{L^1(\OO)}}\right)\, \mu(\d \bx, \d y)
+\frac12(L_{G_1}+|\OO|L_{G_2})\\
&\qquad+{\beta\left(2L_F+\frac{L_{G_2}}{2}\right)}
\int_{\b X\times Y}\E\int_0^1\|\nabla\varphi^y(s)\|_{\b H}^2\,\d s\,\mu(\d\bx,\d y).
\end{align*}
This implies, together with \eqref{supp1}, \eqref{aux_supp}, and \eqref{supp2}
that 
\[
\int_{\b X\times Y} \widetilde\Upsilon_n(\bx,y)\, \mu(\d\b x, \d y) \leq C, \quad\forall\,n\in\mathbb N,
\]
and \eqref{supp3} follows by the monotone convergence theorem. This proves the estimate \eqref{eq:support}.
\paragraph{\textit{Step 4.}}
We show eventually that for every $s\geq1$ there exists $C>0$, independent of $\mu$, such that 
\begin{equation}
    \label{supp4}
    \int_{\b X\times Y} \norm{\Psi_{s+1}(y)}_{L^1(\OO)}\,\mu(\d\b x, \d y)\leq C\,.
\end{equation}
To this end, for every integer $n\geq2$, we define
$\Psi_{s,n}:\mathbb R\to\mathbb R_+$ as
\[
  \Psi_{s,n}(r):=\Psi_s(r)\pmb{1}_{\{|r|<\sqrt{1-n^{-1}}\}}
  +n^s\pmb{1}_{\{|r|\geq \sqrt{1-n^{-1}}\}}\,, \quad r\in\mathbb R,
\]
so that $\Psi_{s,n}\in\mathcal{B}_b(\mathbb R)$, 
$\Psi_{s,n}=\Psi_s$ in $[-\sqrt{1-n^{-1}},\sqrt{1-n^{-1}}]$,
$\Psi_{s,n}\leq\Psi_s$ in $\mathbb R$, and 
$\Psi_{s,n}\nearrow\Psi_s$ in $(-1,1)$.
By arguing as in the proof of Lemma~\ref{lem:F''} 
to deduce \eqref{eq:gs0} but
on the function $\Psi_{s,n}$ instead, we infer that there exists a constant $C>0$, independent of $n$, such that 
\[
\E\int_\OO\Psi_{s,n}(\varphi^y(t))
+\frac1 C\E\int_0^t\int_\OO\Psi_{s+1}(\varphi^y(\tau))
\pmb{1}_{\{|\varphi^y|<\sqrt{1-n^{-1}}\}}
\,\d \tau
\leq \int_\OO\Psi_{s,n}(y) + C\left(1 + t\right)
\quad\forall\,t>0.
\]
Now, we set $\widetilde\Psi_{s,n}, \overline\Psi_{s+1,n}:
\b X\times Y\to\erre$ as
\[
\widetilde\Psi_{s,n}(\b x, y):= \|\Psi_{s,n}(y)\|_{L^1(\OO)}\,,
\qquad
\overline\Psi_{s+1,n}(\b x,y):=\|\Psi_{s+1}(y)
\pmb{1}_{\{|y|<\sqrt{1-n^{-1}}\}}\|_{L^1(\OO)}\,,
\qquad (\b x, y) \in \b X\times Y,
\]
so that $\widetilde\Psi_{s,n}, \overline\Psi_{s+1,n}\in\mathcal{B}_b(\b X\times Y)$.
By taking $t=1$ and by integrating with respect to $\mu$, we get 
\begin{align*}
&\int_{\b X\times Y} 
\widetilde\Psi_{s,n}(\bx,y)\, \mu(\d\b x, \d y)
+\frac1 C\int_{\b X\times Y} 
\overline\Psi_{s+1,n}(\bx,y)\, \mu(\d\b x, \d y)\\
&= 
\int_{\b X\times Y} 
\widetilde\Psi_{s,n}(\b x,y)\, \mu(\d\b x, \d y)
+\frac1 C\int_0^1 \int_{\b X\times Y} 
\overline\Psi_{s+1,n}(\b x,y)\, \mu(\d\b x, \d y)\, {\rm d}\tau\\
&=\int_{\b X\times Y} 
(P_1\widetilde\Psi_{s,n})(\b x,y)\, \mu(\d\b x, \d y)
+\frac1 C\int_0^1 \int_{\b X\times Y} 
(P_\tau\overline\Psi_{s+1,n})(\b x,y)\, \mu(\d\b x, \d y)\, {\rm d}\tau
\\
&= \int_{\b X\times Y}
\E\left[\Psi_{s,n}(\varphi^y(1))\right]\, \mu(\d\b x, \d y)
+\frac1 C\int_0^1 \int_{\b X\times Y}
\E\int_\OO\Psi_{s+1}(\varphi^y(\tau))
\pmb{1}_{\{|\varphi^y(\tau)|<\sqrt{1-n^{-1}}\}}
\, \mu(\d\b x, \d y)\,\d \tau
\\
&\leq \int_{\b X\times Y}\widetilde\Psi_{s,n}(\b x,y)\,\mu(\d\b x, \d y)
+2C,
\end{align*}
from which we deduce, thanks to the definition of $\overline\Psi_{s+1,n}$, that 
\[
\frac1 C\int_{\b X\times Y} 
\norm{\Psi_{s+1,n}(y)\pmb{1}_{\{|y|<\sqrt{1-n^{-1}}\}}}_{L^1(\OO)}\, \mu(\d\b x, \d y)\leq 2C, \quad\forall\,n\geq2.
\]
Hence, \eqref{supp4} follows by letting $n\to+\infty$ 
by using the monotone convergence theorem.
\end{proof}\noindent
The main consequence of Theorem~\ref{th:support}
concerns the existence of ergodic invariant measures.
\begin{thm} 
\label{th:ergodic} 
Let either Assumptions \ref{hyp:structural}--\ref{hyp:diffusionAC} hold with $Z=\bVsd$ or Assumptions \ref{hyp:structural}--\ref{hyp:additionalAC} hold with $Z=\bHs$. Then,
there exists an ergodic invariant measure for the transition semigroup $P$.
\end{thm}
\begin{proof}
The set $\Pi\subset \cP(\b X\times Y)$ of all
invariant measures for the Markov semigroup $P$ 
is not empty by Theorem~\ref{th:ex_inv}, 
and convex by definition of invariance.
By the Krein-Milman Theorem (see e.g.\cite[ Theorem 7.68]{Ali}), 
any convex compact set of a locally convex Hausdorff spaces
\review{possesses extreme points and equals the} closed convex hull of its extreme points.
Moreover, by \cite[Theorem 19.25]{Ali} one has that
the ergodic measures for $P$
are the extreme points of $\Pi$.
Hence, we only have to show that the closure of $\Pi$ in the
narrow topology of $\mathscr P(\b X\times Y)$ is compact, 
or equivalently that $\Pi$ is tight on $\Xbw\times\Ybw$. 
Thanks to \eqref{eq:support},
there exists a constant $C$ such that
\begin{equation*}
\int_{ \b X\times Y}\left(\|\bx\|^2_{\bVs} + \|y\|_{V_2}^2\right)
\, \mu( \d x, \d y) \le C \qquad \forall\, \mu \in \Pi.
\end{equation*}
\review{For any} $n\in\mathbb N$, \review{we define} $D_n:=\{(\b x,y)\in \bVs\times V_2:
\|\bx\|_{\bVs}+\|y\|_{V_2}\leq n\}$ and $D'_n:=D_n\cap(\b X\times Y)$, 
\review{so} that $D_n$ is compact in $\bHs\times V_1$ \review{and} $D_n'$ is compact in $\Xs\times \Ys$ (hence also in $\Xbw\times \Ybw$).
By the Markov inequality we infer
that
\begin{align*}
\sup_{\mu \in \Pi} \mu \left((\b X\times Y)\setminus D_n'\right)
&=\sup_{\mu \in \Pi} \mu (\{ (\b x,y) \in \b X\times Y: \ \|\bx\|_{\bVs}+\|y\|_{V_2} >n\})\\
&\le \frac{2}{n^2}\sup_{\mu \in \Pi} 
\int_{\b X\times Y}\left(\|\bx\|^2_\bVs+\|y\|^2_{V_2}\right)\,
\mu({\rm d}\b x, \d y)\le \frac{2C}{n^2} \rightarrow 0,
\end{align*}
as $ n \rightarrow +\infty$. \review{Therefore,} $\Pi$ is tight and admits extreme points, which are \review{precisely the} ergodic invariant measures for $P$.
\end{proof}

\section{Uniqueness and asymptotic stability of the invariant measure} \label{sec:uniqueness}

Throughout this Section we assume $\alpha_d = 1$ and $\alpha_n = 0$ in Assumption \ref{hyp:structural}, that is, Dirichlet boundary condition in the Allen-Cahn equation. 
Here we prove the uniqueness of the invariant measure for \eqref{AC_NS} and the weak convergence to it, also named asymptotic stability of the invariant measure. These results are based on the abstract results in \cite{GHMR17} and \cite{KS}, relying on generalized coupling techniques; we recall them in Section \ref{KS_sec}. In Section \ref{FP_sec} we prove some Foias-Prodi-type estimates in expected value. These estimates are the crucial ingredient to apply the abstract results of \cite{GHMR17} and \cite{KS} in our framework and infer the uniqueness of the invariant measure and its asymptotic stability. This is the content of Sections \ref{uniq_sec} and \ref{asym_sec}.

Let us start by recalling the generalized Poincar\'e inequality, that will be a crucial ingredient in what follows.
We denote by $\{\lambda_n\}_{n \in \mathbb{N}}$ the eigenvalues of $\b A$ and by $\{\be_n\}_{n \in \mathbb{N}}$ the corresponding eigenvectors of $\b A$ that form a complete orthonormal system in $\bHs$. Moreover, $0<\lambda_1\le \lambda_2 \le \cdots,$ and 
\[
\lim_{n \rightarrow \infty}\lambda_n=+ \infty.
\]
Denoting by $P_N$ and $Q_N$ the orthogonal projection in $\bHs$ onto the space $\operatorname{Span}\{\be_n\}_{1\le n\le N}$ and onto its orthogonal space, respectively, we have the generalized Poincar\'e inequalities 
\begin{equation}
\label{gen_Poin}
\|P_N\b v\|_{\bVs}^2 \le \lambda_N\|P_N \b v\|^2_{\bHs}, \qquad \|Q_N\b v\|^2_{\bHs} \le \frac{1}{\lambda_N}\|Q_N \b v\|^2_{\bVs},
\end{equation}
that hold for all $\b v\in \bVs$ and any $N\ge 1$.

In addition to \ref{hyp:structural}-\ref{hyp:additionalAC}, from now on, we assume the following.
\begin{enumerate}[start = 7, label=\textbf{(A\arabic*)}]
\itemsep0.3em
\item\label{hyp:nondeg}
There exists a measurable map $G_1^{-1}:\bHs \rightarrow \mathcal L(\bHs,U_1)$ and a positive integer $M\in\enne$ such that
\begin{equation}
\label{bound_g}
\sup_{\b v \in \bHs}\|G_1^{-1}(\b v)\|_{L(\bHs,U_1)} < \infty
\end{equation}
and
\begin{equation}
\label{GgM}
G_1(\b v)G_1^{-1}(\b v)=P_M \qquad \forall\, \b v \in \bHs.
\end{equation}
\end{enumerate}

\begin{remark}
The existence of a map $G_1^{-1}: \bHs \rightarrow \mathcal L(\bHs,U_1)$  fulfilling \eqref{GgM} is equivalent to the following property 
\begin{equation*}
P_M\bHs \subseteq \text{{Im}}\ G_1(\b v) \qquad \forall\, \b v \in \bHs.
\end{equation*}
Assumption \ref{hyp:nondeg} can be seen as a non degeneracy condition on the low modes. We refer to \cite[Remarks 3.1 and 3.2]{Oda2008} for more details. For a concrete example of an operator satisfying Assumptions \ref{hyp:diffusionNS} and \ref{hyp:nondeg} see \cite[Example 2.4]{FZ23}. 
\end{remark}

From now on we will reserve the symbol $M$ to denote the integer that appears in Assumption \ref{hyp:nondeg}. To infer the uniqueness and asymptotic stability of the invariant measure, at some point, we require $M$ to be sufficiently large. This means that the range of the covariance operator $G_1$ contains the so-called unstable directions. This is usually referred as the ``effectively elliptic setting''.

\subsection{The abstract results in \cite{GHMR17} and \cite{KS}}
\label{KS_sec}

The proof of the uniqueness of the invariant measure relies on the abstract results in \cite{GHMR17}, whereas the proof of its asymptotic stability relies on the abstract results in \cite{KS}. In these two papers the authors provide sufficient conditions in terms of generalized couplings for the uniqueness and the weak convergence to the invariant measure.
The statement of these results requires to fix some notation.

For a measurable space $(S, \mathcal{S})$ we denote the set of all probability measures on $(S, \mathcal{S})$ by $\cP(S)$. For two given measures $\vartheta_1, \vartheta_2 \in \cP(S)$, we define
\begin{equation*}
    \mathcal{C}(\vartheta_1, \vartheta_2):= \{ \xi \in \cP\left( S \times S\right) \ : \pi_1(\xi)=\vartheta_1, \pi_2(\xi)=\vartheta_2\},
\end{equation*}
where $\pi_i(\xi)$ denotes the $i$-th marginal distribution of $\xi$, $i=1,2$. Any element $\xi \in \mathcal{C}(\vartheta_1, \vartheta_2)$ is called a \textit{coupling} for $\vartheta_1, \vartheta_2$. We now introduce an extension of the notion of coupling. Recall that $\vartheta_1 \ll \vartheta_2$ means that $\vartheta_1$ is absolutely continuous with respect to $\vartheta_2$. 
 We define 
 \[
     \widehat{\mathcal{C}}(\vartheta_1, \vartheta_2):= \{ \xi \in \cP\left( S \times S \right) \ :\pi_1(\xi)\ll\vartheta_1, \pi_2(\xi)\ll\vartheta_2\},
 \]
 and call any probability measure from the class 
 $\widehat{\mathcal{C}}(\vartheta_1, \vartheta_2)$ a \textit{generalized coupling} for $\vartheta_1, \vartheta_2$.

Let now $(S, \rho)$ be a Polish space
. We recall that the 1-Wasserstein distance 
is defined as
\[
\mathcal{W}_\rho(\vartheta_1, \vartheta_2):= \inf_{\xi \in \mathcal{C}(\vartheta_1, \vartheta_2)} \int_{S \times S} \rho(x, y) \, \xi({\rm d}x, {\rm d} y), \qquad \vartheta_1, \vartheta_2 \in \cP(S).
\]
Moreover, we denote by $\text{L}_b(S)$ the space of bounded Lipschitz functions on $S$, that is the space of bounded $\rho$-continuous functions $\phi: S\rightarrow \mathbb{R}$ for which 
\[
\text{Lip}(\phi):= \sup_{(x, y) \in S,\atop x\ne y }\frac{\left|\phi(x)-\phi(y)\right|}{\rho(x,y)} < \infty,
\]
and note that $\text{L}_b(S)$ determines measures in $\mathscr P(S)$.
We recall that from the Kantorovich–Rubinstein theorem (see e.g. \cite[Theorem 1.2.32]{KS_BOOK}) it holds that
\begin{equation}
\label{KW}
\|\vartheta_1-\vartheta_2\|_{K}=\mathcal{W}_\rho(\vartheta_1, \vartheta_2), \qquad \vartheta_1, \vartheta_2 \in \mathscr{P}(S),
\end{equation}
where the Kantorovich distance is defined as 
\begin{equation*}
\|\vartheta_1-\vartheta_2\|_{K}:=\sup_{\phi \in \text{L}_b(S),\atop\text{Lip}(\phi) \le 1}\left|\int_{S}\phi\, {\rm d}\vartheta_1-\int_{S}\phi\, {\rm d}\vartheta_2  \right|, \qquad \vartheta_1, \vartheta_2 \in \mathscr{P}(S).
\end{equation*}

Together with the initial metric $\rho$ on $S$, let $d$ be another metric on $S$ such that $d \le 1$ and $d$ is continuous with respect to the metric $\rho$. The metric $d$ is
not assumed to be complete. When $d \ne \rho$, we denote by $\overline S^d$ the completion of $S$ with respect to $d$, and regard $S$ as a subset in $\overline S^d$. Notice that $(\overline S^d, d)$ is a Polish space (we denote the extended metric again by $d$). We also assume that  for any $y \in S$ there exists a sequence of $d$-continuous functions $\rho_n^y: \overline S^d \rightarrow [0, \infty)$ such that for $z \in \overline S^d$
\[
\rho_n^y(z) \rightarrow 
\begin{cases}
\rho(z,y), & z \in S,
\\
+ \infty, & \text{otherwise},
\end{cases}
\]
as $n \rightarrow \infty$. 
We use the notation $\overset{d}{\rightharpoonup}$ for the weak convergence in $\cP(S)$ with respect to the metric $d$: namely, given a sequence $\{\vartheta_k\}_k\subset \cP(S)$, we say that  $\vartheta_k \overset{d}{\rightharpoonup} \vartheta$ if 
\begin{equation*}
    \int_S \phi \, {\rm d} \vartheta_k \rightarrow \int_S \phi \, {\rm d}\vartheta, \quad \text{as} \ k \rightarrow \infty,
\end{equation*}
for any $\phi: S \rightarrow \mathbb{R}$ which is bounded and $d$-continuous.

Let $Z=\{Z_n\}_{n \in \mathbb{N}}$ be a Markov chain with state space $(S, \mathcal{S})$ and let $P=(P_n)_{n\in\enne}$ be its 
transition discrete semigroup. 
For $z\in S$ we denote the $n$-step transition probability 
for $Z$ as 
$P_n^*\delta_{z}: (S, \mathcal{S})\rightarrow [0,1]$,
and recall that 
$
P_n^*\delta_{z}(A)
:=\mathbb{P}(Z_n\in A|Z_0=z)$ for $A \in \mathcal{S}$.
Analogously, for
$\mu\in\mathscr P(S)$ we denote by 
$P_n^*\mu\in\mathscr P(S)$ the push-forward measure of $\mu$
through $P_n$, namely 
$
P_n^*\mu(A)
:=\int_S\mathbb{P}(Z_n\in A|Z_0=z)
\:\mu(\d x)$ for $A \in \mathcal{S}$,
and recall that $\mu\in\mathscr P(S)$ is invariant for $P$ if $P_n^*\mu = \mu$ for all $n\in\enne$.
We consider the space of one-sided infinite sequences $S^{\mathbb{N}}$ with its Borel $\sigma$-field $\mathcal{S}^{\otimes \mathbb{N}}$. The law of the sequence $\{Z_n\}_{n \in \mathbb{N}}$ in $(S^{\mathbb{N}},\mathcal{S}^{\otimes \mathbb{N}})$ with initial datum $z \in S$ is denoted by $\mathbb{P}_{z} \in \cP\left(S^{\mathbb{N}}\right)$. We call the Markov chain $d$-Feller if for each bounded and $d$-continuous function $\phi:S \rightarrow \mathbb{R}$, the map $z \mapsto P_n \phi(z)$ is $d$-continuous
for every $n\in\enne$.

We now state the three main abstract results from \cite{GHMR17} and \cite{KS}. The first one provides the uniqueness of the invariant measure, the second one provides some asymptotic stability results for an invariant  measure, whereas the last one provides at the same time uniqueness of the invariant measure and weak convergence to it. We state these results in the best form that fits our framework, thus not considering the most possible general case as in \cite{GHMR17} and \cite{KS}.

Let us start by recalling \cite[Corollary 2.1]{GHMR17} that provides the \textit{uniqueness} of the invariant measure.
\begin{thm}
    \label{GHMRthm}
    Let $Z$ be a Markov chain with state space $(S,\mathcal{S})$.
If $Q \in \mathcal S$ is a measurable set such that for any $z, \widetilde z \in Q$ there exists some $\xi:=\xi_{z, \widetilde z} \in \widehat{\mathcal{C}}(\mathbb{P}_{z}, \mathbb{P}_{\widetilde z})$ with 
\begin{equation}
\label{cond_GHMR}
\xi \left((\mathsf z, \widetilde{\mathsf z}) \in S^\mathbb{N} \times S^{\mathbb{N}} \ : \ 
 \lim_{n \rightarrow \infty} 
 \rho\left(\mathsf z_n, \widetilde{\mathsf z}_n\right)=0
\right)>0,
\end{equation}
then there exists at most one ergodic invariant probability measure $\vartheta$ for $Z$ with $\vartheta(Q)>0$.
\end{thm}

We now state \cite[Corollary 2]{KS}. Assuming the existence of an invariant measure $\vartheta$, this result provides a sufficient condition for weak convergence of transition probabilities for $\vartheta$-almost all initial conditions. Convergence is in probability, i.e.~the measure $\vartheta$ of the set of initial
conditions for which the $\mathcal{W}_\rho$-distance of the transition probability to the invariant measure $\vartheta$ after $n$ steps is larger than $\varepsilon$ converges to $0$ for every $\varepsilon>0$.

\begin{thm}
\label{KS_thm_weak}
Let $Z$ be a $\rho$-Feller Markov chain with state space $(S,\mathcal{S})$. Let $\vartheta$ be an invariant probability measure for $Z$ and let $Q\in \mathcal S$ be a measurable set with $\vartheta(Q)=1$. Assume that for any $z,\widetilde z \in Q$ there exist some $\xi:=\xi_{z, \widetilde z} \in \widehat{\mathcal{C}}(\mathbb{P}_{z},\mathbb{P}_{\widetilde z})$ and $\alpha_{z, \widetilde z}>0$ such that 
\begin{equation}
\label{coupling_cond_weak}
\liminf_{n \rightarrow \infty} \xi\left((\mathsf z, \widetilde{\mathsf z}) \in S^\mathbb{N} \times S^{\mathbb{N}} \ : \ \rho\left(\mathsf z_n, \widetilde{\mathsf z}_n\right) \le \varepsilon\right)\ge \alpha_{z, \widetilde z},
\end{equation}
for any $\varepsilon>0$.
Then,
\begin{equation}
\label{conv_in_prob}
\lim_{n\to\infty}\vartheta\left\{z\ : \ \mathcal{W}_{\rho}\left(P^*_n\delta_z, \vartheta\right)>\varepsilon\right\}= 0 \quad\forall\, \varepsilon >0,
\end{equation}
i.e., $(P^*_n\delta_{z})_{n\in\enne}$, considered as a sequence of $\cP(S)$-valued random elements on
$(S, \mathcal{S}, \vartheta)$, $\mathcal{W}_{\rho}$-converges in probability to the random element identically equal to $\vartheta$.
\end{thm}
\begin{remark}
Notice that, in view of Remark \ref{KW}, \eqref{conv_in_prob} can be restated as follows
\begin{equation*}
\lim_{n\to\infty}\vartheta\left\{ z\ : \ \|P^*_n\delta_z- \vartheta\|_K>\varepsilon\right\}= 0
\quad \forall\,\varepsilon >0.
\end{equation*}
\end{remark}

To conclude, we state \cite[Corollary 4]{KS}, which provides the uniqueness of the invariant measure and the weak convergence of transition probabilities starting from a given point $z \in S$ to it.
\begin{thm}
\label{KS_thm_strong}
Let $Z$ be a $d$-Feller Markov chain on the state space $(S, \mathcal{S})$. Assume that for any $z, \widetilde z \in S$ there exists some $\xi:=\xi_{z, \widetilde z} \in \widehat{\mathcal{C}}(\mathbb{P}_z, \mathbb{P}_{\widetilde{z}})$ such that $\pi_1(\xi) \sim \mathbb{P}_z$ and 
\begin{equation}
\label{coupling_cond_strong}
\lim_{n \rightarrow \infty} \xi \left( (\mathsf z, \widetilde{\mathsf z}) \in S^{\mathbb{N}} \times S^{\mathbb{N}} \ : \ d(\mathsf z_n, \widetilde{\mathsf z}_n) \le \varepsilon \right)=1 \quad\forall\,\varepsilon>0.
\end{equation}
Then, there exists at most one invariant probability measure for $Z$ and, if such a measure $\vartheta$ exists, then 
\[
P^*_n\delta_z \overset{d}{\rightharpoonup} \vartheta
\quad\forall\, z\in S.
\]
\end{thm}

\subsection{Foias-Prodi estimates in expected value}
\label{FP_sec}
As we have anticipated, in this section we assume that \ref{hyp:structural}--\ref{hyp:additionalAC} hold with $Z=\bHs$
and that $\alpha_d=1$ and $\alpha_n=0$, i.e.~we only consider Dirichlet boundary conditions.
Moreover, we will always assume that the initial data belong to the space $\b X \times Y_{s_0}$; we recall that this space have been introduced in \eqref{new_spaces}.

This Section is devoted to establishing a Foias-Prodi type estimate for the Allen-Cahn-Navier-Stokes system \eqref{AC_NS} that holds in expected value. This result will serve as a crucial technical tool to apply the abstract results recalled in Section \ref{KS_sec}. 
In our framework the Foias-Prodi estimates will describe the behaviour for which any two solutions of system \eqref{AC_NS}, with different initial data, converge to each other as $t\to+\infty$ if a control acts on a sufficient  finite number $N$ of modes in the Navier-Stokes system and if the viscosity coefficient $\beta$ in the Allen-Cahn equation is sufficiently large.

Given $(\bx,y) \in \b X \times Y$, let $(\bu, \varphi)=(\bu^{\bx},\varphi^y)$ denote the  solution of  the Allen-Cahn-Navier-Stokes system \eqref{AC_NS}.
Given $\eta>0, N>0$ and $(\widetilde \bx,\widetilde y) \in X \times Y$, let $(\widetilde \bu, \widetilde \varphi)=(\widetilde\bu^{\widetilde \bx},\widetilde\varphi^{\widetilde y})$ denote the solution of the system 
\begin{equation}
\label{AC_NS_nud}
\begin{cases}
    {\rm d} \widetilde{\bu}+ \left[ -\nu\Delta \widetilde{\bu}+(\widetilde{\bu} \cdot \nabla )\widetilde{\bu}+ \nabla \pi-\widetilde{w} \nabla \widetilde{\varphi}\right]\, {\rm d}t=G_1(\widetilde{\bu})\, {\rm d}W_1 -\eta P_N(\widetilde{\bu}-\bu)\,{\rm d}t & \text{in} \ (0,T) \times \mathcal{O},
    \\
    \nabla \cdot \widetilde{\bu}=0 & \text{in} \ (0,T) \times \mathcal{O},
    \\
    {\rm d}\widetilde{\varphi} +\left[\widetilde \bu \cdot \nabla \widetilde{\varphi}+ \widetilde w\right]\, {\rm d}t =G_2(\widetilde \varphi)\, {\rm d}W_2 & \text{in} \ (0,T) \times \mathcal{O},
    \\
    \widetilde w=-\beta\Delta \widetilde{\varphi} +F'(\widetilde{\varphi}) & \text{in} \ (0,T) \times \mathcal{O},
    \\
    \widetilde{\bu}=0, \  \widetilde {\varphi} =0, & \text{in} \ (0,T) \times \partial\mathcal{O},
    \\
    \widetilde{\bu}(0)=\widetilde x, \ \widetilde{\varphi}(0)=\widetilde y & \text{in}\ \mathcal{O},
    \end{cases}
\end{equation}
where $P_N$ is the orthogonal projection from  $\bHs$ onto the space $\operatorname{Span}\{\be_n\}_{1\le n\le N}$. Here $\eta>0$ is a parameter to be  suitably chosen later on, possibly depending on $N$.
From now on, we reserve the symbol $N$ to indicate the dimension of the projected space $P_N\bHs$ where the control $P_N(\widetilde{\bu}-\bu)$  acts. At some point we will require $N$ to be sufficiently large. 

We will refer to \eqref{AC_NS_nud} as the \textit{nudged system} corresponding to the Allen-Cahn-Navier-Stokes system \eqref{AC_NS}. The well posedness of \eqref{AC_NS_nud} can be easily proved  since the additional term $\eta P_N(\widetilde{\bu}-\bu)$ 
does not crucially impact the well-posedness estimates (see \cite[Remark 8]{KS} for further details).

The effect of the \textit{nudging term} $\eta P_N(\widetilde{\bu}-\bu)$ in the Navier-Stokes equations is to drive $\bu$ towards $\widetilde \bu$ on the space $P_N\bHs$, i.e.~on the low modes. 
The Foias-Prodi estimates (in expectation) will quantify how many modes $N$ need to be activated and how large
the viscosity coefficient $\beta$ needs to be in order to drive $(\bu, \varphi)$ towards $(\widetilde \bu, \widetilde \varphi)$ on the full (infinite-dimensional) space  $\bHs \times V_1$. 
More formally, we will show that $\mathbb{E}\left[\beta\|\widetilde  \varphi(t)-  \varphi(t)\|^2_{V_1} + \|\widetilde \bu(t)-\bu(t)\|^2_{\bHs} \right]$ vanishes as $t \rightarrow + \infty$, provided the number of modes $N$ and the viscosity coefficient $\beta$ are taken sufficiently large. By sufficiently large we mean that $N$ and $\beta$ has to be larger to some integer parameter $\overline N \ge 1$ and some real parameter $\overline \beta>0$ satisfying the following joint condition.

\begin{condition}
\label{conditionFP}
Let $\overline N \in\enne_+$ and $\overline \beta > 0$ satisfy
\begin{multline}
    \label{conditionNbeta}
 \left(\frac 14 \min\left\{ \nu \lambda_{\overline N}, \frac{\overline{\beta}^2}{K_\Delta^2}\right\} - \bC_G \right)
    \min\left\{ \frac{\nu\min\{\nu,\overline\beta\}^3}
    {\bK_1(1+\nu^{-2})},
    \frac{\overline\beta^2\nu^2}
 {L_{G_2}K_4^4\nu^2 +
 \bK_2\max\{\nu,\overline\beta\}}, \frac{\overline \beta}{\bK_3}\right\}
    \\
    \ge 3\max((\mathfrak{C}_1+2),\mathfrak{C}_2+2,2(\mathfrak{C}_3+1)),
\end{multline}
where 
\begin{equation}
    \label{C_noise}
\bC_G:= 2\max\left\{ L_{G_1}, L_{G_2}\right\},
\end{equation}
and
\begin{equation}
\label{K_123}
    \bK_1:= \max\left\{27K_L^8, 1024 K_L^8K_\Delta^2\right\},
    \quad
    \bK_2:= 
    256 K_L^4K_{GN}^4K_\Delta^4, \quad
    \bK_3:= 4 K_{\frac{2\gamma}{\gamma-2}}^2,
\end{equation}
and $\mathfrak{C}_1, \mathfrak{C}_2, \mathfrak{C}_3$ are the positive constants that appear in Lemma \ref{lemma1}, \ref{lemma2}, \ref{lemma3} and depend on the structural parameters $s_0, s_F, \lambda_1, L_{G_1},L_{G_2},L_F,|\OO|,\|F\|_{\C([-1,1])}$ only.
\end{condition}

We do not expect Condition \ref{conditionFP} to be sharp. However the meaning of the condition is clear: the minimum numbers of modes $\overline N$ on which the finite-dimensional control acts and the minimum value of the viscosity coefficient $\overline \beta$ depend on the viscosity coefficient $\nu>0$, the intensity of the noises $L_{G_1}, L_{G_1}$ and the domain. We obtain a \textit{joint} condition on $\overline N$ and $\overline\beta$, which should be intended as follows:
provided that $\overline\beta$ is sufficiently large, for every 
$\nu>0$ arbitrarily small there exists $\overline N$ such that 
Condition \ref{conditionFP} holds.
A joint condition on the parameters is not surprising since we are dealing with a coupled system of equations.

The Foias-Prodi-type estimates we derive read as follows.

\begin{thm}[Foias-Prodi estimates]
\label{FP_thm_complete}
Let Assumptions \ref{hyp:structural}-\ref{hyp:additionalAC} hold with $Z=\bHs$, $\alpha_d = 1$, $\alpha_n = 0$,
and let $\overline\beta, \overline N$ satisfy Condition~\ref{conditionFP}.
Let $(\bx, y), (\widetilde \bx, \widetilde y) \in \b X \times Y_{s_0}$, $N\geq\overline N$, $\beta\geq\overline\beta$, and $\eta=\frac{\lambda_N\nu}2$.
Let also $(\bu,\varphi)$ be the unique strong solution
to the system \eqref{AC_NS} starting from the initial state
$(\b x,y)$, and
let 
$(\widetilde \bu,\widetilde \varphi)$ be the unique strong solution
to the system \eqref{AC_NS_nud} starting from the initial state
$(\widetilde{\b x},\widetilde y)$.
Then, for every $p>0$, there exists a constant $C>0$,
depending on the structural coefficients and on the initial data,
such that
\begin{equation*}
    \mathbb{E}\left[\beta\|\widetilde  \varphi(t)-  \varphi(t)\|^2_{V_1} + \|\widetilde \bu(t)-\bu(t)\|^2_{\bHs} \right] \le \frac{C}{t^p} \quad \forall\,t \ge 0.
\end{equation*}
\end{thm}

The proof of the above result is inspired by \cite{GHMR21} and \cite{FZ23}, and will follow after some intermediate steps.
First, in Proposition~\ref{FP_thm} we prove some
preliminary stability estimates by exploiting the techniques of 
Proposition~\ref{prop:uniqueness}.
Secondly, in Corollary~\ref{corFP} we show that for a suitable stopping time $\tau_{R, \varepsilon}^{N, \beta}$, that controls the growth of the solutions to \eqref{AC_NS} and \eqref{AC_NS_nud}, the expectation 
\[
\mathbb{E}\left[\pmb{1}_{\left(\tau_{R,\varepsilon}^{N, \beta}=+ \infty\right)} \left(\beta\|\widetilde  \varphi(t)-  \varphi(t)\|^2_{V_1} + \|\widetilde \bu(t)-\bu(t)\|^2_{\bHs}\right)\right]
\]
decays with an exponential rate in time. Eventually,  in Proposition~\ref{stimeP} we exploit energy estimates for the solutions to \eqref{AC_NS} and \eqref{AC_NS_nud} to prove that the probability such stopping time $\tau_{R, \varepsilon}^{N, \beta}$ remains finite decays (with a polynomial rate) in the cut-off parameter $R$, when $N$ and $\beta$ are large enough (i.e. are larger to some $\overline N$, $\overline \beta$ satisfying Condition \ref{conditionFP}) and for a suitable choice of the parameter $\varepsilon$. 
The Foias-Prodi estimates then easily follow, by combining the results of Corollary~\ref{corFP} and Proposition~\ref{stimeP}.

Let us start with the first preliminary result.
\begin{prop}
\label{FP_thm}
Let Assumptions \ref{hyp:structural}-\ref{hyp:additionalAC} hold with $Z=\bHs$, $\alpha_d = 1$, $\alpha_n = 0$,
and let $\overline\beta, \overline N$ satisfy Condition~\ref{conditionFP}.
Let $(\bx, y), (\widetilde \bx, \widetilde y) \in \b X \times Y_{s_0}$, $N\geq\overline N$, $\beta\geq\overline\beta$, and $\eta=\frac{\lambda_N\nu}2$.
Let also $(\bu,\varphi)$ be the unique strong solution
to the system \eqref{AC_NS} starting from the initial state
$(\b x,y)$, and
let 
$(\widetilde \bu,\widetilde \varphi)$ be the unique strong solution
to the system \eqref{AC_NS_nud} starting from the initial state
$(\widetilde{\b x},\widetilde y)$.
Then, defining the process $\Lambda:\Omega\times[0, + \infty)\rightarrow \mathbb{R}$ as 
\begin{align}
\nonumber
\Lambda(t)&:= 
\left(\frac 12 \min\left\{\nu \lambda_N, \frac{\beta^2}{K^2_\Delta}\right\} - \bC_G\right)t -\frac{\bK_1}{\min\{\nu,\beta\}^3}\int_0^t\|\bu(s)\|_{\bHs}^2\|\bu(s)\|_{\bVs}^2  \,{\rm d}s
\\ \nonumber
&- \left(L_{G_2}C_4^4+\frac{\bK_2}\nu\max\left\{1,\frac\beta{\nu}\right\}\right)
\int_0^t\|\varphi(s)\|^2_{V_2}\, {\rm d}s\\
\label{Gamma}
&- \frac{\bK_3}{\beta} \int_0^t \left(\|F''(\varphi(s))\|^2_{L^\gamma(\OO)}
+\|F''(\widetilde{\varphi}(s))\|^2_{L^\gamma(\OO)}\right)\, {\rm d}s,
\end{align}
the estimate 
\begin{multline}
\label{FP_est}
\mathbb{E} \left[e^{\Lambda(t \wedge \sigma)}\left(\beta\|\widetilde  \varphi(t\wedge \sigma)-  \varphi(t\wedge \sigma)\|^2_{V_1} + \|\widetilde \bu(t\wedge \sigma)-\bu(t\wedge \sigma)\|^2_{\bHs}\right)\right]
\\
+ \frac{1}{2}\min\left\{ \nu \lambda_N, \frac{\beta^2}{K_\Delta^2}\right\} \mathbb{E}\left[ \int_0^{t \wedge \sigma} e^{\Lambda(\tau)} \left(\beta\| \widetilde \varphi(\tau)-\varphi(\tau)\|^2_{V_1} + \|\widetilde \bu(\tau)-\bu(\tau)\|^2_{\bHs}\right) \, {\rm d}\tau\right]\\
\le \beta\|\widetilde y-y\|^2_{V_1} + \|\widetilde \bx-\bx\|^2_{\bHs},
\end{multline}
holds for every stopping time $\sigma\ge 0$ and every $t \ge 0$.
\end{prop}
\begin{proof}
For the sake of convenience, let us define the differences
	\begin{align*}
		\hbu := \widetilde \bu-\bu, \quad
		\hphi := \widetilde \varphi-\varphi, \quad
		\hw  := \widetilde w-w, \quad
		\hbu_0  := \widetilde{\bx}-\bx, \quad
		\hphi_0  := \widetilde y -y.
	\end{align*}
Arguing as in the proof of Proposition \ref{prop:uniqueness} we obtain the estimate
 \begin{multline} \label{eq:cd4_D_FP}
		\dfrac{\beta}{2}\| \hphi(t)\|^2_{V_1} + \dfrac{1}{2}\|\hbu(t)\|^2_{\bHs} + \dfrac{1}{2}\int_0^t \left[ \beta^2\|\Delta \hphi(\tau)\|^2_H +\nu\|\nabla\hbu(\tau)\|^2_{\bHs} \right] \: \d \tau 
  + \eta \int_0^t \|P_N\hbu(\tau)\|^2_{\bHs}\, {\rm d}\tau
  \\
		\leq \dfrac{\beta}{2}\|\hphi_{0}\|^2_{V_1} 
        + \dfrac{1}{2}\|\hbu_{0}\|^2_{\bHs} \\
        + 2\int_0^t \left[ L_{G_1} + \dfrac{27\KL^8}{4\nu^3}\|\bu(\tau)\|_{\bHs}^2\|\bu(\tau)\|_{\bVs}^2
        +\dfrac{128}{\nu}\KL^4\KGN^4K_\Delta^4\| \varphi(\tau)\|_{V_2}^2 \right]\dfrac{1}{2}\|\hbu(\tau)\|_{\bHs}^2 \: \d \tau \\
        + 2\int_0^t \left[L_{G_2} + \frac{512}{\beta^3}\KL^8K_\Delta^2\|\bu(\tau)\|_{\bHs}^2\|\bu(\tau)\|_{\bVs}^2 + \dfrac{2}{\beta}K_{\frac{2\gamma}{\gamma-2}}^2 \left(\|F''(\varphi(\tau))\|^2_{L^\gamma(\OO)}
        +\|F''(\widetilde{\varphi}(\tau))\|^2_{L^\gamma(\OO)}\right) \right. \\
        \left.+ \left( L_{G_2}K_4^4 + \dfrac{32\beta}{\nu^2}\KL^4\KGN^4K_\Delta^4 \right)\|\varphi(\tau)\|^2_{ V_2} \right]\dfrac{\beta}{2}\|\hphi(\tau)\|_{V_1}^2 \: \d \tau+ \sum_{k = 2}^{3} \mathcal{M}_k(t),
	\end{multline}
that holds for any $t \ge 0$, $\P$-almost surely, where the martingale processes $\{\mathcal M_k\}_{k \in \{2,3\}}$ are defined as in Proposition \ref{prop:uniqueness}. The only difference with respect to the estimate \eqref{eq:cd4_D} is the presence of the additional term $\eta \int_0^t \|P_N\hbu(\tau)\|^2_{\bHs}\, {\rm d}\tau$ in the left hand side.
Thanks to the generalized Poincar\'e inequality \eqref{gen_Poin}, we obtain 
\begin{multline}
\label{P1}
\frac{\nu}{2}\|\nabla\hbu(\tau)\|^2_{\bHs} + \eta \|P_N \hbu(\tau)\|^2_{\bHs} 
\ge 
\frac{\nu}{2}\|\nabla Q_N\hbu(\tau)\|^2_{\bHs} + \eta \|P_N \hbu(\tau)\|^2_{\bHs} 
\\
\ge 
\frac{\nu \lambda_N}{2}\|Q_N\hbu(\tau)\|^2_{\bHs} + \eta \|P_N \hbu(\tau)\|^2_{\bHs}
=\frac{\nu \lambda_N}{2}\|\hbu(\tau)\|^2_{\bHs}.
\end{multline}
Moreover, the embedding $V_2 \subset V_1$ yields,
together with the Dirichlet boundary conditions, the estimate 
\begin{equation}
\label{P2}
\frac{\beta^2}{2} \|\Delta \hphi(\tau)\|^2_{H} \ge \frac{\beta^2}{2K_\Delta^2}\|\hphi(\tau)\|^2_{V_1}.
\end{equation}
By inserting estimates \eqref{P1} and \eqref{P2} in estimate \eqref{eq:cd4_D_FP} we obtain
\begin{multline*} 
		\beta\| \hphi(t)\|^2_{V_1} + \|\hbu(t)\|^2_{\bHs}
   \\
        + 2\int_0^t \left[\frac{\nu \lambda_N}2 - L_{G_1} - \dfrac{27\KL^8}{2\nu^3}\|\bu(\tau)\|_{\bHs}^2\|\bu(\tau)\|_{\bVs}^2-\dfrac{128}{\nu}\KL^4\KGN^4K_\Delta^4\| \varphi(\tau)\|_{V_2}^2 \right]\|\hbu(\tau)\|_{\bHs}^2 \: \d \tau \\
        + 2\int_0^t \left[\frac{\beta^2}{2K_\Delta^2}
        -L_{G_2} - \frac{512}{\beta^3}\KL^8K_\Delta^2
        \|\bu(\tau)\|_{\bHs}^2\|\bu(\tau)\|_{\bVs}^2 
        - \dfrac{2}{\beta} K_{\frac{2\gamma}{\gamma-2}}^2 \left(\|F''(\varphi(\tau))\|^2_{L^\gamma(\OO)}
        +\|F''(\widetilde{\varphi}(\tau))\|^2_{L^\gamma(\OO)}\right) \right. \\
        \left.- \left( L_{G_2}K_4^4 + \dfrac{32\beta}{\nu^2}\KL^4\KGN^4K_\Delta^4 \right)\|\varphi(\tau)\|^2_{ V_2} \right]\beta\|\hphi(\tau)\|_{V_1}^2\: \d \tau
		\leq \beta\|\hphi_{0}\|^2_{V_1} + \|\hbu_{0}\|^2_{\bHs} + 2\sum_{k = 2}^{3} \mathcal{M}_k(t),
	\end{multline*} 
from which we obtain the differential estimate 
\begin{multline*}
{\rm d}\left(\beta\| \hphi(t)\|^2_{V_1} 
+ \|\hbu(t)\|^2_{\bHs}\right)
\\
+\left[ \min\left\{\nu \lambda_N, 
\frac{\beta^2}{K_\Delta^2}\right\} - \bC_G -\frac{\bK_1}{\min\{\nu,\beta\}^3}\|\bu(t)\|_{\bHs}^2\|\bu(t)\|_{\bVs}^2  
-\left(L_{G_2}K_4^4+\frac{\bK_2}\nu\max\left\{1,\frac\beta{\nu}\right\}\right)
\|\varphi(t)\|^2_{V_2}
\right.
\\
\left.
- \frac{\bK_3}{\beta} \left(\|F''(\varphi(t))\|^2_{L^\gamma(\OO)}
+\|F''(\widetilde{\varphi}(t))\|^2_{L^\gamma(\OO)}\right) 
\right] \left( \beta\| \hphi(t)\|^2_{V_1} + \|\hbu(t)\|^2_{\bHs} \right)\, {\rm d}t \le 2\sum_{k = 2}^{3} {\rm d}\mathcal{M}_k,
	\end{multline*}
 where the constants $\bC_G, \bK_1, \bK_2$ and $\bK_3$ are the ones appearing in Condition \ref{conditionFP}.
Recalling the definition of the functional $\Lambda$ in \eqref{Gamma}, we rewrite the above estimate as 
\begin{equation*}
{\rm d}\left(\beta\| \hphi(t)\|^2_{V_1} + \|\hbu(t)\|^2_{\bHs}\right) + \left(\frac12 \min\left\{\nu \lambda_N, \frac{\beta^2}{K_\Delta^2}\right\} + \Lambda'(t)\right)\left(\beta\| \hphi(t)\|^2_{V_1} + \|\hbu(t)\|^2_{\bHs}\right)\,{\rm d}t \le 2\sum_{k = 2}^{3} {\rm d}\mathcal{M}_k(t),
\end{equation*}
which yields
\begin{multline*}
{\rm d}\left(e^{\Lambda(t)}\left(\beta\| \hphi(t)\|^2_{V_1} + \|\hbu(t)\|^2_{\bHs}\right)\right)
+ \frac{1}{2}\min\left\{\nu \lambda_N, \frac{\beta^2}{K_\Delta^2}\right\} 
e^{\Lambda(t)} \left(\beta\| \hphi(t)\|^2_{V_1} + \|\hbu(t)\|^2_{\bHs}\right)\,{\rm d}t\\
\le 2 e^{\Lambda(t)}\sum_{k = 2}^{3} {\rm d}\mathcal{M}_k(t).
\end{multline*}
The proofs follows then
by integrating in time up to a stopping time $\sigma$ and taking expectations.
\end{proof}

For our purpose, we will have to make a suitable choice of the stopping time $\sigma$ appearing in \eqref{FP_est}. 
To this end, for every $R,\varepsilon,\beta>0$ and $N \ge 1$, let
\begin{equation}
\label{def_st}
    \tau_{R, \varepsilon}^{N, \beta}:= \inf\left\{r \ge 0 \ : \ \frac14\min\left\{\nu \lambda_N, \frac{\beta^2}{K_\Delta^2}\right\}r - \Lambda(r)-\varepsilon \ge R\right\},
\end{equation}
with the usual convention that $\tau_{R, \varepsilon}^{N, \beta}:=+\infty$ if the set on the right-hand side of \eqref{def_st} is empty, i.e.
\[
\{\tau_{R, \varepsilon}^{N, \beta}=+\infty\}
=\left\{\frac14\min\left\{\nu\lambda_N, \frac{\beta^2}{K_\Delta^2}\right\}r- \Lambda(r)-\varepsilon < R, \quad \forall \ r   \ge 0\right\}.
\]
The parameter $\varepsilon$ will be useful to track the dependence of the initial data $(\bx, y), (\widetilde \bx, \widetilde y)$ in subsequent estimates of $\tau_{R, \varepsilon}^{N, \beta}$; see Proposition \ref{stimeP} below.
With $\tau_{R, \varepsilon}^{N, \beta}$ in hand we draw the following immediate Corollary of Proposition \ref{FP_thm}.
\begin{cor}
\label{corFP}
    Under the conditions of Proposition~\ref{FP_thm}, for every $R, \varepsilon > 0$ it holds that
\begin{multline*}
\mathbb{E}\left[\pmb{1}_{\left(\tau_{R,\varepsilon}^{N, \beta}=+ \infty\right)} \left(\beta\|\widetilde  \varphi(t)-  \varphi(t)\|^2_{V_1} + \|\widetilde \bu(t)-\bu(t)\|^2_{\bHs}\right)\right]\\
\le
e^{R +\varepsilon-\frac 14
\min\left\{\nu \lambda_N, \frac{\beta^2}{K_\Delta^2}\right\}t} \left(\beta\|\widetilde y-y\|^2_{V_1} + \|\widetilde \bx-\bx\|^2_{\bHs}\right) \qquad\forall\,t\geq0.
\end{multline*}
\end{cor}
\begin{proof}
    It is an immediate consequence of the proof of Proposition~\ref{FP_thm}, provided to notice that, 
    if $\tau_{R, \varepsilon}^{N, \beta}=+\infty$, then 
\[
\frac14\min\left\{\nu \lambda_N, \frac{\beta^2}{K_\Delta^2}\right\}t - \Lambda(t)-\varepsilon \le R \quad \forall \ t   \ge 0,
\]
i.e.
\[
\Lambda(t) \ge \frac14\min\left\{\nu \lambda_N, \frac{\beta^2}{K_\Delta^2}\right\}t -\varepsilon - R \quad \forall \ t   \ge 0.
\qedhere
\]
\end{proof}
 
We now estimate the probability $\mathbb{P}(\tau_{R, \varepsilon}^{N, \beta}< \infty)$ in terms of the parameter $R$, when $N$ and $\beta$ are large enough and for a suitable choice of the parameter $\varepsilon$. We will obtain a polynomial decay in $R$.
\begin{prop}
\label{stimeP}
Let Assumptions \ref{hyp:structural}-\ref{hyp:additionalAC} hold with $Z=\bHs$, $\alpha_d = 1$, $\alpha_n = 0$,
and let $\overline\beta, \overline N$ satisfy Condition~\ref{conditionFP}.
Let $(\bx, y), (\widetilde \bx, \widetilde y) \in \b X \times Y_{s_0}$, $N\geq\overline N$, $\beta\geq\overline\beta$, and $\eta=\frac{\lambda_N\nu}2$.
Let also $(\bu,\varphi)$ be the unique strong solution
to the system \eqref{AC_NS} starting from the initial state
$(\b x,y)$,
let 
$(\widetilde \bu,\widetilde \varphi)$ be the unique strong solution
to the system \eqref{AC_NS_nud} starting from the initial state
$(\widetilde{\b x},\widetilde y)$, and set
\begin{align}
\nonumber
    \varepsilon &:= 
\max\left\{\frac{3 \bK_1}{\nu \min\{\nu,\overline\beta\}^3}(\mathfrak{C}_1+2)(1 + \|\bx\|^4_{\bHs}+\|y\|^4_{H} + \beta^2\|\nabla y\|_\bH^4),
\right.
\\ \nonumber
&\qquad\qquad\left.
\frac{3(L_{G_2}C_4^4\nu^2 +
 \bK_2\max\{\nu,\beta\})}{\beta^2\nu^2}
 (\mathfrak{C}_2+2)(1 + 
\|\bx\|^2_{\bHs}+\|y\|^2_{H} + \beta\|\nabla y\|_\bH^2, \right. \\
\label{conditionvareps}
&\qquad\qquad\left.\frac{6\bK_3}{\overline\beta}
(\mathfrak{C}_3+1)\left(1+ \int_\OO \Psi_{s_0}(y) + \int_\OO \Psi_{s_0}(\widetilde y)\right) \right\}.
\end{align}
Then, for every $p>0$ there exists a constant $C$, 
depending on the structural parameters, but independent of the initial data, such that
\begin{equation*}
 \mathbb{P}\left(\tau_{R, \varepsilon}^{N, \beta}<+ \infty\right) \le \frac{C\left( 1 +  \|\bx\|^{6(p+1)}_{\bHs}+\|y\|^{6(p+1)}_{V_1} \right)}{R^p}, \qquad\forall\,R>0.
\end{equation*}
\end{prop}
\begin{proof}
Keeping in mind the definition of the stopping time $\tau_{R, \varepsilon}^{N, \beta}$, we introduce the set 
    \[
    A_{R,\varepsilon}^{N ,\beta}:=\left\{ \sup_{r \ge 0}\left(\frac 14\min\left(\nu \lambda_N, \frac{\beta^2}{K_\Delta^2}\right)r - \Lambda(r)-\varepsilon \right) \ge R\right\}\in\cF.
    \]
Recalling the definition  \eqref{Gamma} of the functional $\Lambda$, we obtain
\begin{equation*}
    A_{R, \varepsilon}^{N, \beta} \subseteq A_{R,\varepsilon}^{N, \beta,1} \cup A_{R,\varepsilon}^{N, \beta,2} \cup A_{R,\varepsilon}^{N, \beta,3}, 
\end{equation*}
where 
\begin{multline*}
A_{R,\varepsilon}^{N, \beta,1}:=\left\{\sup_{r \ge 0} \left(\frac{\bK_1}{\min\{\nu,\beta\}^3}
\int_0^r \|\bu(\tau)\|^2_{\bHs}\|\bu(\tau)\|^2_{\bVs}\,{\rm d}\tau 
\right.\right.
\\ \left.\left.+ \frac 13 \left[\bC_G- \frac 14 \min\left\{\nu \lambda_N, \frac{\beta^2}{K_\Delta^2}\right\} \right]r-\frac{\varepsilon}{3}\right) \ge \frac{R}{3}
\right\},
\end{multline*}
\begin{multline*}
A_{R,\varepsilon}^{N, \beta,2}:=\left\{\sup_{r \ge 0} \left(\left(L_{G_2}C_4^4+\frac{\bK_2}\nu\max\left\{1,\frac\beta{\nu}\right\}\right)\int_0^r\|\varphi(\tau)\|^2_{V_2}\, {\rm d} \tau
\right.\right.
\\
\left. \left. +  \frac 13 \left[\bC_G
- \frac 14 \min\left\{\nu \lambda_N, \frac{\beta^2}{K_\Delta^2}\right\} \right]r -\frac{\varepsilon}{3}\right) \ge \frac{R}{3}
\right\},
\end{multline*}
\begin{multline*}
A_{R,\varepsilon}^{N, \beta,3}:=\left\{\sup_{r \ge 0} \left(\frac{\bK_3}{\beta}\int_0^r\left(\|F''\left(\varphi(\tau)\right)\|^2_{L^\gamma(\OO)}+\|F''\left(\widetilde{\varphi}(\tau)\right)\|^2_{L^\gamma(\OO)}\right)\, {\rm d} \tau
\right.\right.
\\
\left. \left. +  \frac 13 \left[\bC_G
- \frac 14 \min\left\{\nu \lambda_N, \frac{\beta^2}{K_\Delta^2}\right\} \right]r -\frac{\varepsilon}{3}\right) \ge \frac{R}{3}
\right\}.
\end{multline*}
Thus we infer
\begin{equation}
\label{Pr1}
\mathbb{P}(\tau_{R, \varepsilon}^{N, \beta} < + \infty)
\le \mathbb{P}(A_{R, \varepsilon}^{N, \beta}) \le \mathbb{P}(A^{N, \beta,1}_{R, \varepsilon})+\mathbb{P}(A^{N, \beta,2}_{R, \varepsilon})+\mathbb{P}(A^{N, \beta,3}_{R, \varepsilon}).
\end{equation}
We write the complementary set of $A_{R, \varepsilon}^{N, \beta,1}$ as
\begin{multline*}
(A_{R, \varepsilon}^{N, \beta,1})^C = \left\{
\nu \int_0^r \|\bu(\tau)\|^2_{\bHs}\|\bu(\tau)\|^2_{\bVs}\,{\rm d}\tau  \right.
\\\notag
 \left. < \frac{\nu \min\{\nu,\beta\}^3}{3 \bK_1}\left[\left(\frac 14 \min\left(\nu \lambda_N, \frac{\beta^2}{K_\Delta^2}\right) - \bC_G \right)r+\varepsilon +R\right] \quad \forall \ r \ge 0
\right\}.
\end{multline*}
Since $\overline N , \overline \beta$ are as in \eqref{conditionNbeta}, then for $N\ge \overline N$ and $\beta \ge \overline\beta$ it holds
\[
 \frac{\nu \min\{\nu,\beta\}^3}{3 \bK_1} \left(\frac 14 \min\left\{\nu \lambda_{N}, \frac{{\beta}^2}{K_\Delta^2}\right\} - \bC_G \right)
    \ge (\mathfrak{C}_1+2)(1+\nu^{-2}) ,
\]
while the definition \eqref{conditionvareps} of
$\varepsilon>0$ implies that
\[
\frac{\nu \min\{\nu,\beta\}^3}{3 \bK_1}\varepsilon\ge (\mathfrak{C}_1+2)(1+
\|\bx\|^4_{\bHs}+\|y\|^4_{H} + \beta^2\|\nabla y\|_\bH^4).
\]
For this choice of the parameters and setting $\overline R:=\frac{\nu \min\{\nu,\beta\}^3}{3 \bK_1}R$, 
we have the inclusion 
\begin{multline*}
(A_{r, \varepsilon}^{N, \beta,1})^C \supseteq \left\{ \nu \int_0^R  \|\bu(\tau)\|^2_{\bHs}\|\bu(\tau)\|^2_{\bVs}\,{\rm d}\tau \right.\\
\left.<  (\mathfrak{C}_1+2)\left( 1+(1+\nu^{-2}) r
+\|\bx\|^4_{\bHs}+\|y\|^4_{H} + \beta^2\|\nabla y\|_\bH^4\right) + \overline R
\quad \forall \ r\ge 0\right\}.
\end{multline*}
Therefore, from Proposition \ref{prop1} we infer that, for any $p>0$,
\begin{equation}
\label{Pr2}
\mathbb{P}(A^{N, \beta,1}_{R, \varepsilon}) \le \frac{C\left(1+\|\bx\|_{\bHs}^{6(p+1)}+ \|y\|_{H}^{6(p+1)}+\beta^{3(p+1)}\|\nabla y\|_\bH^{6(p+1)}\right)}{R^p},
\end{equation}
with $C$ a positive constant depending on  $p$ and the structural parameters of the problem, but not depending on $R$ and the initial data.
We argue in a similar way for the sets $A_{R, \varepsilon}^{N, \beta,2}$ and $A_{R, \varepsilon}^{N, \beta,3}$. 
We write the complementary set of $A_{R, \varepsilon}^{N, \beta,2}$ as
\begin{multline*}
(A_{R, \varepsilon}^{N,\beta, 2})^C = \left\{
\beta^2 \int_0^r \|\varphi(\tau)\|^2_{V_2}\,{\rm d}\tau   \right.
\\\notag
 \left.<\frac{\beta^2\nu^2}
 {3(L_{G_2}C_4^4\nu^2 +
 \bK_2\max\{\nu,\beta\})}
\left[\left(\frac 14 \min\left(\nu \lambda_N, \frac{\beta^2}{K_\Delta^2}\right) - \bC_G\right)r+\varepsilon +R\right] \quad \forall \ r \ge 0
\right\}
\end{multline*}
Since $\overline N , \overline \beta$ are as in Condition \ref{conditionFP}, then for $N\ge \overline N$ and $\beta \ge \overline\beta$ it holds
\[
\frac{\beta^2\nu^2}
 {3(L_{G_2}C_4^4\nu^2 +
 \bK_2\max\{\nu,\beta\})}\left(\frac 14\min\left\{\nu \lambda_{ N}, \frac{{\beta}^2}{K_\Delta^2}\right\} - \bC_G \right)
    \ge (\mathfrak{C}_2+2),
\]
while the definition \eqref{conditionvareps} of
$\varepsilon>0$ implies that
\[
\frac{\beta^2\nu^2}
 {3(L_{G_2}C_4^4\nu^2 +
 \bK_2\max\{\nu,\beta\})}\varepsilon\ge (\mathfrak{C}_2+2)(1+ \|\bx\|^2_{\bHs}
 +\|y\|^2_{H} + \beta\|\nabla y\|_\bH^2).
\]
For this choice of the parameters and setting $\overline R:=\frac{\beta^2\nu^2}
 {3(L_{G_2}C_4^4\nu^2 +
 \bK_2\max\{\nu,\beta\})}R$, 
we have the inclusion 
\[
(A_{R, \varepsilon}^{N, \beta, 2})^C \supseteq \left\{ \beta^2 \int_0^t\|\varphi(\tau)\|_{V_2}^2\,\d \tau < (\mathfrak{C}_2+2)\left(1 +\|\bx\|^2_{\bHs}+ \|y\|^2_{H} + \beta\|\nabla y\|_\bH^2 +r\right)  + \overline R \quad \forall \ r\ge 0\right\}
\]
so that from Proposition \ref{prop2} we infer that, for any $p>0$,
\begin{equation}
\label{Pr3}
\mathbb{P}(A^{N, \beta,2}_{R, \varepsilon}) \le \frac{C\left(1+\|\bx\|_{\bHs}^{4(p+1)}+ \|y\|_{H}^{4(p+1)} + \beta^{2(p+1)}\|\nabla y\|_\bH^{4(p+1)}\right)}{R^p},
\end{equation}
with $C$ a positive constant depending on  $p$ and the structural parameters of the problem  and not depending on $R$ and the initial data.
Eventually, we write the complementary set of $A_{R, \varepsilon}^{N, \beta,3}$ as
\begin{multline*}
(A_{R, \varepsilon}^{N,\beta, 3})^C = \left\{
 \int_0^r  \left(\|F''({\varphi}(\tau))\|_{L^{\gamma}(\OO)}^2+
 \|F''(\widetilde{\varphi}(\tau))\|_{L^{\gamma}(\OO)}^2\right) \: \d \tau  
 \right.
 \\
 \left. <\frac{\beta}{3\bK_3}\left[
 \left(\frac 14 \min\left(\nu \lambda_N, 
 \frac{\beta^2}{K_\Delta^2}\right) - \bC_G \right])r+ \varepsilon +R\right] \quad \forall \ r \ge 0
\right\}.
\end{multline*}
Since $\overline N , \overline \beta$ are as in Condition \ref{conditionFP}, then for $N\ge \overline N$ and $\beta \ge \overline\beta$ it holds
\[
\frac{ \beta }{3  \bK_3}\left(\frac 14 \min\left\{\nu \lambda_{\overline N}, \frac{{\beta}^2}{K_\Delta^2}\right\} - \bC_G \right)
    \ge 2(\mathfrak{C}_3+1),
\]
while the definition \eqref{conditionvareps} of
$\varepsilon>0$ implies that
\[
\frac{\beta \varepsilon }{3\bK_3}\ge 2(\mathfrak{C}_3+1)\left(1+  \int_\OO \Psi_{s_0}(y) + \int_\OO \Psi_{s_0}(\widetilde y)\right).
\]
For this choice of the parameters and setting $\overline R:=\frac{\beta }{3\bK_3}R$, 
we have the inclusion 
\begin{multline*}
(A_{r, \varepsilon}^{N, \beta, 3})^C \supseteq 
\left\{\int_0^r \left( \|F''(\varphi(\tau))\|_{L^{\gamma}(\OO)}^2+ \|F''(\widetilde{\varphi}(\tau))\|_{L^{\gamma}(\OO)}^2 \right) \,\d \tau \right.
\\
\left. \le 2({\mathfrak{C}_3}+1)\left(1 + r+ \int_\OO \Psi_{s_0}(y)  +  \int_\OO \Psi_{s_0}(\widetilde y)  \right) + \overline R \quad \forall \ r\ge 0\right\},
\end{multline*}
hence from Proposition \ref{prop3} we infer that, for any $p>0$,
\begin{equation}
\label{Pr4}
\mathbb{P}(A^{N, \beta, 3}_{R, \varepsilon}) \le \frac{C}{R^p},
\end{equation}
with $C$ a positive constant depending on  $p$ and the structural parameters of the problem  and not depending on $R$ and the initial data.
The thesis follows from \eqref{Pr1}-\eqref{Pr4} and the Young inequality.
\end{proof}

\begin{remark}
We highlight that, in order to prove Proposition \ref{stimeP}, we need to control the growth of both the solutions to systems \eqref{AC_NS} and \eqref{AC_NS_nud} (see Proposition \ref{prop1}-\ref{prop3}).  
We observe that when dealing with the Navier-Stokes equations, and in general with the many examples considered in \cite{GHMR17} and \cite{KS}, it is sufficient to control the growth of suitable norms of the solution of the equation without the nudging term. This allows to have estimates independent of the parameter $N$ and consequently to easily derive a condition on the minimal number of modes that must be stochastically forced (see, for example, condition (23) in \cite{FZ}). The difficulty with the system \eqref{AC_NS} lies in the fact that we have to also control the growth of the solution of the nudged system \eqref{AC_NS_nud} via the term $\int_0^t \|F''(\widetilde\varphi(s))\|_{L^\gamma(\OO)}^2\, {\rm d}s$. Luckily (see the proof of Lemmata \ref{lem:F''} and \ref{lemma3}) the estimates we obtain are independent of $N$. In fact, the convective term in the Allen–Cahn equation disappears when computing the estimates for $\|F''(\widetilde\varphi(t))\|_{L^\gamma(\OO)}^2$. Therefore, we do not see the dependence on the first equation in \eqref{AC_NS_nud} and obtain estimates independent of the parameter $N$.
\end{remark}

Thanks to Corollary \ref{corFP} and Proposition \ref{stimeP} we have now all the ingredients to prove Theorem \ref{FP_thm}.

\begin{proof}[Proof of Theorem \ref{FP_thm_complete}]
Set $\varepsilon$ as in \eqref{conditionvareps}.
Thanks to the H\"older and Young inequalities, as a consequence of Corollary \ref{corFP} and Propositions \ref{a_priori_est2}  and \ref{a_priori_est2_bis}  we obtain the estimate 
\begin{multline*}
\mathbb{E}\left[\beta\|\widetilde  \varphi(t)-  \varphi(t)\|^2_{V_1} + \|\widetilde \bu(t)-\bu(t)\|^2_{\bHs} \right] 
\le
\mathbb{E}\left[\pmb{1}_{\left(\tau_{R,\varepsilon}^{N, \beta}=+ \infty\right)} \left(\beta\|\widetilde  \varphi(t)-  \varphi(t)\|^2_{V_1} + \|\widetilde \bu(t)-\bu(t)\|^2_{\bHs}\right)\right]
\\
+
\mathbb{E}\left[\pmb{1}_{\left(\tau_{R,\varepsilon}^{N, \beta}<+ \infty\right)} \left(\beta\|\widetilde  \varphi(t)-  \varphi(t)\|^2_{V_1} + \|\widetilde \bu(t)-\bu(t)\|^2_{\bHs}\right)\right]
\\
\le
C e^{R +\varepsilon-\frac 14\min\left\{\nu \lambda_N, \frac{\beta^2}{K_\Delta^2}\right\}t} \left(\beta\|\widetilde y-y\|^2_{V_1} + \|\widetilde \bx-\bx\|^2_{\bHs}\right)
\\
+ \sqrt{\mathbb{P}\left( \tau_{R,\varepsilon}^{N, \beta}<+ \infty\right)}\sqrt{\mathbb{E}\left[\beta^2\|\widetilde  \varphi(t)-  \varphi(t)\|^4_{V_1} + \|\widetilde \bu(t)-\bu(t)\|^4_{\bHs}\right]}
\\
\le C\left(1+ \|\bx\|^2_{\bHs}+\|\widetilde{\bx}\|^2_{\bHs} + \beta\|y\|^2_{V_1}+ \beta\|\widetilde y\|^2_{V_1}\right)
\left(e^{R +\varepsilon-\frac 14\min\left\{\nu \lambda_N, \frac{\beta^2}{K_\Delta^2}\right\}t} +  \sqrt{\mathbb{P}\left( \tau_{R,\varepsilon}^{N, \beta}<+ \infty\right)}\right),
\end{multline*}
where $C$ is a positive constant depending on the structural coefficients, but independent of $R$
and the initial data.
Furthermore, since 
$N \ge \overline N$ and $\beta>\overline \beta$,
thanks to Condition~\ref{conditionFP}
and Proposition~\ref{stimeP} we get,
possibly updating the constant $C$,
\begin{multline*}
\mathbb{E}\left[\beta\|\widetilde  \varphi(t)-  \varphi(t)\|^2_{V_1} + \|\widetilde \bu(t)-\bu(t)\|^2_{\bHs} \right] 
\\
\le 
C\left(1+ \|\bx\|^2_{\bHs}+\|\widetilde{\bx}\|^2_{\bHs} + \beta\|y\|^2_{V_1}+ \beta\|\widetilde y\|^2_{V_1}\right)
\\
\times\left(e^{R +\varepsilon-\frac 14\min\left\{\nu \lambda_N, \frac{\beta^2}{K_\Delta^2}\right\}t} + \sqrt{\frac{C\left( 1 +  \|\bx\|^{6(p+1)}_{\bHs}+\|y\|^{6(p+1)}_{H}
+\beta^{3(p+1)}\|\nabla y\|_\bH^{6(p+1)}\right)}{R^p}}\right).
\end{multline*}
Selecting now $R= \frac{1}{8}\min\left\{\nu \lambda_N,\frac{\beta^2}{K_\Delta^2}\right\}t$ for any $t>0$, the thesis follows.
\end{proof}

\subsection{Uniqueness of the invariant measure}
\label{uniq_sec}

We prove the uniqueness of the invariant measure by relying on Theorem \ref{GHMRthm}. 
In order to  
exploit this result, the idea is to consider the nudged system \eqref{AC_NS_nud} and verify that:
\begin{itemize}
    \item [(i)] the law of the solution to system \eqref{AC_NS_nud} on the space of trajectories is absolutely continuous w.r.t. the law of the solution to system \eqref{AC_NS}; 
    \item [(ii)]  for any pair of distinct initial conditions, there is a positive probability that solutions to these systems converge at large times, when evaluated on a infinite sequence of evenly spaced times.
\end{itemize}
The key ingredient to verify the above two facts is given by the Foias-Prodi estimates of Theorem \ref{FP_thm_complete}. These require to consider the parameters $N$ and $\beta$, that appear in systems \eqref{AC_NS} and \eqref{AC_NS_nud}, to be larger to some $\overline N$ and $\overline \beta$ as in Condition \ref{conditionFP}.
Condition (ii) is a consequence of Theorem \ref{FP_thm_complete}.
On the other hand, to ensure condition (i) we impose 
Assumption \ref{hyp:nondeg} to hold for some $M \ge N$. 
This is enough to apply Girsanov theorem and ensure that the law of the solution to system \eqref{AC_NS_nud} on the space of trajectories is absolutely continuous w.r.t. the law of the solution to system \eqref{AC_NS}; see Proposition \ref{Girsanov_prop} below. 
This non-degenaracy condition of the noise on the low modes resembles the one  usually assumed  in presence of an \textit{additive} noise (see e.g. the many examples in \cite{GHMR17} and \cite{KS}). Roughly speaking the condition requires the noise to span the (unstable) directions controlled in the Foias-Prodi estimates. 

Thus we prove the uniqueness of the invariant measure provided that $\beta \ge \overline \beta$ and Assumption \ref{hyp:nondeg} holds for some $M \ge \overline N$. In particular, the last condition allow to consider \textit{any} viscosity coefficient $\nu$ in the Navier-Stokes equations, provided the operator $G_1$ is non-degenerate on a sufficiently large number of direction and the viscosity parameter $\beta$ in the Allen-Chan equation is sufficiently large. Dropping the non-degeneracy assumption on $G_1$ would result in requiring the viscosity $\nu$ to be sufficiently large. 
On the other hand, we observe that the noise structure in the Allen-Cahn equation prevents to impose any kind of non-degeneracy condition on $G_2$. Therefore, the viscosity coefficient $\beta$ must necessarily be sufficiently large, as required in Condition~\ref{conditionFP}.

Our main result reads as follows.
\begin{thm}
   \label{asinvmeas} 
Let Assumptions \ref{hyp:structural}-\ref{hyp:additionalAC} hold with $Z=\bHs$, $\alpha_d = 1$, $\alpha_n = 0$,
and let $\overline\beta, \overline N$ satisfy Condition~\ref{conditionFP}.
Let also $N\geq\overline N$, $\beta\geq\overline\beta$,  
and let
Assumption~\ref{hyp:nondeg} hold for some $M \ge \overline N$.
Then, the transition semigroup $P$ associated with system \ref{AC_NS} possesses at most one invariant measure $\vartheta$, which is in particular ergodic and such that  $\vartheta(\b X_{reg} \times Y_{reg})=1$.
\end{thm}

We start by proving that the law of the solution of system \eqref{AC_NS_nud} is absolutely continuous w.r.t.~the law of the solution to system \eqref{AC_NS}.

\begin{prop}
\label{Girsanov_prop}
Let Assumptions \ref{hyp:structural}-\ref{hyp:additionalAC} hold with $Z=\bHs$, $\alpha_d = 1$, $\alpha_n = 0$,
and let $\overline\beta, \overline N$ satisfy Condition~\ref{conditionFP}.
Let $(\bx, y), (\widetilde{\b x},\widetilde y) \in \b X \times Y_{s_0}$, $N\geq\overline N$, $\beta\geq\overline\beta$,  $\eta=\frac{\lambda_N\nu}2$, and let
Assumption~\ref{hyp:nondeg} hold for some $M \ge \overline N$.
Let also $(\bu,\varphi)$ and $(\widetilde \bu,\widetilde \varphi)$
be the unique strong solutions
to the systems \eqref{AC_NS} and \eqref{AC_NS_nud}, repsectively,  starting from the same initial state
$(\widetilde{\b x},\widetilde y)$.
Then,  the law of $(\widetilde \bu,\widetilde \varphi)$
 is absolutely continuous w.r.t.~the law of $(\bu,\varphi)$
 as measures on $\mathcal C(\erre_+;\b X \times Y)$. 
\end{prop}
\begin{proof}
We start by rewriting equations \eqref{AC_NS} and \eqref{AC_NS_nud} in a more compact form. To this aim we formally introduce the operators
\[
(\bv, \psi) \mapsto \mathbb C(\bv, \psi):=\left(- \Delta \bv + (\bv\cdot \nabla)\bv+ \nabla \pi-(-\beta \Delta \psi +F'(\psi))\nabla \psi, \bv \cdot \nabla \psi -\beta\Delta\psi +F'(\varphi) \right),
\]
and
\[
(\bv, \widetilde {\bv}, \psi) \mapsto \mathbb L(\bv,\widetilde {\bv}, \psi):=\left(\eta P_N(\widetilde \bv - \bv),0\right).
\]
Moreover, we denote by $\mathbb{U}$ the separable Hilbert space $U_1 \times U_2$, endowed with the norm
\[
\|(u_1, u_2)\|_{\mathbb{U}}:= \sqrt{\|u_1\|^2_{U_1}+\|u_2\|^2_{U_2}}, \qquad (u_1, u_2) \in \mathbb{U}.
\]
On the normal filtered probability space $(\Omega,\cF,(\cF_t)_{t\in[0,T]},\P)$, we consider the $\mathbb U$-cylindrical Wiener processes $\mathbb{W}:=(W_1,W_2)$.
We introduce the operator $\mathbb{G}: \bHs \times B^\infty_1 \rightarrow \LL_{HS}(\mathbb{U}, \bHs \times H)$
\[
\mathbb G(\bv, \psi)[(u_1, u_2)]:=(G_1(\bv)[u_1], G_2(\psi)[u_2]), \qquad \forall \ (u_1, u_2) \in \mathbb{U}.
\]
We thus rewrite the Allen-Cahn-Navier-Stokes system \eqref{AC_NS} as an evolution equation on $\bHs \times H$ as
\begin{equation}
    \label{AC_NS_com}
\begin{cases}
    {\rm d}(\bu, \varphi)+\mathbb C(\bu, \varphi)\,{\rm d}t=\mathbb G(\bu, \varphi)\,{\rm d}\mathbb{W}(t), & \text{in} \ (0,T) \times \mathcal{O},
    \\
    \bu=0, \ \varphi =0, & \text{in} \ (0,T) \times \partial\mathcal{O},
    \\
    (\bu, \varphi)(0)=(\widetilde{\b x},\widetilde y) & \text{in}\ \mathcal{O},
\end{cases}    
\end{equation}
and similarly the nudged Allen-Cahn-Navier-Stokes system \eqref{AC_NS_nud}  as 
\begin{equation}
    \label{AC_NS_nud_com}
\begin{cases}
    {\rm d}(\widetilde \bu, \widetilde \varphi)+\mathbb C(\widetilde \bu, \widetilde \varphi)\,{\rm d}t=-\mathbb L(\bu, \widetilde \bu, \varphi)\, {\rm d}t +\mathbb G(\widetilde \bu, \widetilde \varphi)\,{\rm d}\mathbb{W}(t), & \text{in} \ (0,T) \times \mathcal{O},
    \\
    \widetilde \bu=0, \ \widetilde \varphi=0, & \text{in} \ (0,T) \times \partial\mathcal{O},
    \\
    (\widetilde \bu, \widetilde \varphi)(0)=(\widetilde{\b x},\widetilde y) & \text{in}\ \mathcal{O}.
\end{cases}    
\end{equation}
We introduce the functions $h_1:\Omega\times\erre_+\to U_1$ and $h_2:\Omega\times\erre_+\to U_2$ as
\[
h_1(t):=-\eta G_1^{-1}(\widetilde \bu(t))P_N(\widetilde \bu(t)-\bu(t)) , \qquad h_2(t):=0_{|U_2}, \quad t \ge 0,
\]
with $\eta=\frac{\nu \lambda_N}{2}$ and define the ``Girsanov shift'' $\mathbb{H}:\Omega\times\erre_+\to\mathbb U$ by 
\[
\mathbb{H}(t):=(h_1(t), h_2(t)), \qquad t \ge 0.
\]
Notice that the definition of $h_1$ (and thus of $\mathbb{H}$) does make sense, thanks to Assumption \ref{hyp:nondeg} and the fact we required $M\ge N$, and
we have that $\mathbb{H} \in \mathcal C([0, \infty),\mathbb{U})$. We set 
\[
\widetilde{\mathbb{W}}(t):= \mathbb{W}(t)+ \int_0^t \mathbb{H}(s)\, {\rm d}s, \qquad t \ge 0,
\]
and we notice that system \eqref{AC_NS_nud_com} can be then equivalently rewritten as follows 
\begin{equation}
    \label{AC_NS_nud_com2}
\begin{cases}
    {\rm d}(\widetilde \bu, \widetilde \varphi)+\mathbb C(\widetilde \bu, \widetilde \varphi)\,{\rm d}t=\mathbb G(\widetilde \bu, \widetilde \varphi)\,{\rm d} \widetilde{\mathbb{W}}(t), & \text{in} \ (0,T) \times \mathcal{O},
    \\
    \widetilde \bu=0, \ \widetilde \varphi=0, & \text{in} \ (0,T) \times \partial\mathcal{O},
    \\
    (\widetilde \bu, \widetilde \varphi)(0)=( \bx, y) & \text{in}\ \mathcal{O}.
\end{cases}    
\end{equation}
Let us show that
\begin{equation}
\label{Girsanov}
\mathbb{P}\left( \int_0^\infty\|\mathbb{H}(s)\|^2_{\mathbb{U}}\, {\rm d}s < \infty\right)=1.
\end{equation}
In fact, noting that $\|\mathbb{H}(t)\|^2_{\mathbb{U}}=\|h_1(t)\|^2_{U_1}$ for any $t\ge 0$, given any positive constant $c$ by the Markov inequality we infer 
\begin{align}
\label{Gir}
    \mathbb{P}\left(\int_0^\infty \|\mathbb{H}(s)\|^2_{\mathbb{U}}\,{\rm d}s > c\right)
    &\le \frac{\eta^2}{c^2} \mathbb{E}\left[ \int_0^\infty \|G_1^{-1}(\widetilde \bu(s))P_N(\widetilde \bu(s)-\bu(s))\|^2_{U_1}\, {\rm d}s\right]
    \notag\\
    &\le \frac{\eta^2}{c^2} \sup_{\bv \in \bHs}
    \|G_1^{-1}(\bv)\|_{\mathcal L(U_1,\bHs)} \mathbb{E}\left[ \int_0^\infty \|\widetilde \bu(s)-\bu(s)\|^2_{\bHs}\, {\rm d}s\right].
\end{align}
The term depending on $G_1^{-1}$ is finite thanks to Assumption \ref{hyp:nondeg}, and Theorem \ref{FP_thm_complete} yields 
\[
\|\widetilde \bu(t)-\bu(t)\|^2_{\bHs} \le   \mathbb{E}\left[\beta\|\widetilde  \varphi(t)-  \varphi(t)\|^2_{V_1} + \|\widetilde \bu(t)-\bu(t)\|^2_{\bHs} \right] \le \frac{C}{t^p}, \quad t> 0,
\]
where $p>0$ is any positive number and $C$ a positive constant depending on the structural parameters of the system but independent on $t$. Thus, by exploiting the continuity of $\bu$ and $\widetilde \bu$ we have
\[
\mathbb{E}\left[ \int_0^\infty \|\widetilde \bu(s)-\bu(s)\|^2_{\bHs}\, {\rm d}s\right] \le C
\left(1+ 
\mathbb{E}\left[ \int_1^\infty \|\widetilde \bu(s)-\bu(s)\|^2_{\bHs}\, {\rm d}s\right]\right)
\leq C\left(1+\int_1^\infty \frac{1}{s^p}\, {\rm d}s\right),
\]
which is finite by choosing $p>1$.
Thus, by letting $c$ going to infinity in \eqref{Gir} we infer \eqref{Girsanov}.
Therefore, from the Girsanov theorem we infer that the law of the Wiener process $\widetilde{\mathbb{W}}$ on $\mathcal C([0, \infty); \mathbb{U})$ is
absolutely continuous w.r.t.~the law of the Wiener process $\mathbb{W}$ (see \cite[Theorem 7.4]{LS01} for the finite dimensional case and \cite{F12} for an extension to the infinite-dimensional case). Since $(\widetilde \bu, \widetilde \varphi)$ is the unique strong solution to \eqref{AC_NS_nud_com2}, this yields immediately that the law of $(\widetilde \bu, \widetilde \varphi)$ is absolutely continuous w.r.t.~the law of 
$(\bu,\varphi)$ as measures on the space 
$\mathcal C(\erre_+;\b X \times Y)$. This concludes the proof.
\end{proof}

We have now all the ingredients to prove Theorem \ref{asinvmeas}.

\begin{proof}
[Proof of Theorem \ref{asinvmeas}]
The proof relies on Theorem \ref{GHMRthm}. 
Let $\mathsf{d}$ be the metric on $\b X \times Y$ induced by the norm 
$
 \|\cdot\|^2_{\bHs} +\beta \|\cdot\|^2_{V_1}.
$
Bearing in mind the notation of Section \ref{KS_sec}, we choose $(S, \mathcal{S})=(\b X \times Y, \mathscr{B}(\b X\times Y))$ and the metric $\rho=\mathsf{d}$.
We notice that $(\b X \times Y, \mathsf{d})$ is a Polish space. Moreover, we choose $Q= \b X \times Y_{s_0}$.
We want to apply Theorem \ref{GHMRthm} with the choice
of the Markov chain $\{\bu(n), \varphi(n)\}_{n\in\enne}$
on $(S, \mathcal{S})$.
For $(\bx, y), (\widetilde \bx, \widetilde y)$ in $\b X \times Y_{s_0}$, we denote by $(\bu, \varphi)$ and $(\widetilde \bu, \widetilde \varphi)$ the unique strong solutions to the systems \eqref{AC_NS} and \eqref{AC_NS_nud}, respectively, with corresponding initial data $(\bx, y)$ and  $(\widetilde \bx, \widetilde y)$.
We consider the measure $\xi:=\xi_{(\bx, y), (\widetilde \bx, \widetilde y)}$ on $(\b X \times Y)^{\mathbb{N}} \times (\b X \times Y)^{\mathbb{N}}$ given by the law of the associated random vector $\left\{\left(\bu(n), \varphi(n) \right),\left( \widetilde {\bu}(n),\widetilde{\varphi}(n)\right)\right\}_{n\in \mathbb{N}}$. 
In view of Proposition \ref{Girsanov_prop}, 
the joint law $\xi$ of the pair $(\bu, \varphi)$, $(\widetilde \bu, \widetilde \varphi)$ is a generalized coupling from the class $\widehat{\mathcal{C}}\left( \mathbb{P}_{(\bx, y)}, \mathbb{P}_{(\widetilde{\bx}, \widetilde y)}\right)$ which satisfies the additional condition $\pi_1(\xi)=\mathbb{P}_{(\bx, y)}$. Hence, 
it only remains to prove condition \eqref{cond_GHMR} in Theorem \ref{GHMRthm}, that is 
\[
\mathbb{P} \left(\lim_{n \rightarrow 0} \|\bu(n)- \widetilde{\bu}(n)\|^2_{\bHs}+ \beta \|\varphi(n)- \widetilde{\varphi}(n)\|^2_{V_1} =0 \right)>0.
\]
To this end, for every $n \in \mathbb{N}$ we introduce the events 
\[
B_n:= \left\{  \|\bu(n)- \widetilde{\bu}(n)\|^2_{\bHs}+ \beta \|\varphi(n)- \widetilde{\varphi}(n)\|^2_{V_1} > \frac{1}{n^2}\right\}
\]
and we set 
\[
B:= \bigcap_{m=1}^{+\infty} \bigcup_{n=m}^{+\infty} B_n.
\]
We consider the stopping time $\tau_{R, \varepsilon}^{N, \beta}$ introduced in \eqref{def_st} and write 
\[
\mathbb{P}(B)=\mathbb{P}\left(B \cap \{\tau_{R, \varepsilon}^{N,\beta}=+ \infty\}\right)+\mathbb{P}\left(B \cap \{\tau_{R, \varepsilon}^{N,\beta}<+ \infty\}\right).
\]
We now observe that $\mathbb{P}\left(B \cap \{\tau_{R, \varepsilon}^{N,\beta}=+ \infty\}\right)=0$ for any $R, \varepsilon>0$. In fact, thanks to the Markov inequality and Corollary \ref{corFP}, for every $n \in \mathbb{N}$ it holds
\begin{multline*}
\mathbb{P}\left(B_n \cap \{\tau_{R, \varepsilon}^{N,\beta}=+ \infty\}\right)
\le n^2 \mathbb{E}\left[\pmb{1}_{\left(\tau_{R,\varepsilon}^{N, \beta}=+ \infty\right)} \left(\beta\|\widetilde  \varphi(t)-  \varphi(t)\|^2_{V_1} + \|\widetilde \bu(t)-\bu(t)\|^2_{\bHs}\right)\right]
\\
\le
n^2 e^{-\frac 14\min\left\{\nu \lambda_N, \frac{\beta^2}{K_\Delta^2}\right\}n} e^{R +\varepsilon} \left(\beta\|\widetilde y-y\|^2_{V_1} + \|\widetilde \bx-\bx\|^2_{\bHs}\right).
\end{multline*}
Thus, 
\[
\sum_{n\in\enne_+} \mathbb{P} \left(B_n \cap  \{\tau_{R, \varepsilon}^{N,\beta}=+\infty\} \right) < \infty.
\]
Hence, by the Borel-Cantelli Lemma it follows $\mathbb{P} \left(B \cap  \{\tau_{R, \varepsilon}^{N,\beta}=+\infty \} \right)=0$. We thus infer 
\[
\mathbb{P}(B) = \mathbb{P}\left(B \cap  \{\tau_{R, \varepsilon}^{N,\beta}< + \infty\} \right)\le \mathbb{P} \left(\tau_{R, \varepsilon}^{N,\beta} < + \infty\right).
\]
Since $\varepsilon$ satisfies condition \eqref{conditionvareps}, from Proposition \ref{stimeP} we infer that, for any $R>0$, it holds
\[
 \mathbb{P}\left(\tau_{R, \varepsilon}^{N, \beta}<+ \infty\right) \le \frac{C\left( 1 +  \|\bx\|^{6(p+1)}_{\bHs}+\|y\|^{6(p+1)}_{V_1} \right)}{R^p},
\]
for any $p>0$ and where the constant $C$ does not depend on $R$.
Hence, we can fix $\overline R>0$ such that $\P(B)<\frac12$.
With such choice of $\overline R$, 
 from the continuity from below we can find $m^*>0$ sufficiently large such that 
\[
\mathbb{P}\left( \bigcap_{n=m^*}^{+\infty}B_n^C\right)>\frac 12.
\]
Now we observe that 
\[
\left\{ \lim_{n \rightarrow 0} \|\bu(n)- \widetilde{\bu}(n)\|^2_{\bHs}+ \beta \|\varphi(n)- \widetilde{\varphi}(n)\|^2_{V_1} =0 \right\} \supseteq \bigcap_{n=m^*}^{+\infty}B_n^C,
\]
hence,
\[
\mathbb{P} \left( \lim_{n \rightarrow 0} \|\bu(n)- \widetilde{\bu}(n)\|^2_{\bHs}+ \beta \|\varphi(n)- \widetilde{\varphi}(n)\|^2_{V_1} =0 \right) \ge \mathbb{P} \left(  \bigcap_{n=m^*}^{+\infty}B_n^C\right)> \frac 12 >0.
\]
Since Assumptions of Theorem \ref{GHMRthm} are verified, we conclude that there exists at most one invariant measure $\vartheta$ for the Markov chain $\{(\bu(n), \varphi(n))\}_{n\in\enne}$ with $\vartheta(\b X \times Y_{s_0})>0$. 
Since invariant measures for the semigroup $P$
associated to \eqref{AC_NS} are 
invariant also for such Markov chain,
this implies the desired uniqueness result.
Since by virtue fo Theorem \ref{th:support} all the invariant measures are supported in $\b X_{reg}\times Y_{reg}$, we conclude.
\end{proof}

\subsection{Asymptotic stability of the invariant measure}
\label{asym_sec}
We prove here some asymptotic stability results for the unique invariant measure of system \eqref{AC_NS} relying on Theorems \ref{KS_thm_weak} and \ref{KS_thm_strong}.
In order to
exploit these results, the idea is again to work with the auxiliary nudged system \eqref{AC_NS_nud} and the key ingredient is again given by the Foias-Prodi estimates of Theorem \ref{FP_thm_complete}. 

Our main result is as follows.
\begin{thm}
\label{weak_conv_ACNS}
Let Assumptions \ref{hyp:structural}-\ref{hyp:additionalAC} hold with $Z=\bHs$, $\alpha_d = 1$, $\alpha_n = 0$,
and let $\overline\beta, \overline N$ satisfy Condition~\ref{conditionFP}.
Let $N\geq\overline N$, $\beta\geq\overline\beta$, 
let
Assumption~\ref{hyp:nondeg} hold for some $M \ge \overline N$,
and let $\vartheta$ be the unique invariant measure of system \eqref{AC_NS}.
Then,
denoting by $\mathsf{d}$ the metric on $\b X \times Y$ induced by the norm $\|\cdot\|^2_{\bHs} +\beta \|\cdot\|^2_{V_1}$,
the following holds:
\begin{itemize}
\item[(i)]
$\vartheta$ is asymptotically $\mathcal{W}_{\mathsf{d}}$-stable in probability, that is, for all $\varepsilon>0$,
\[
\vartheta\left\{(\bx, y)\ : \ \mathcal{W}_{\mathsf{d}}\left(P^*_n\delta_{(\b x, y)}, \vartheta\right)>\varepsilon\right\}\rightarrow 0, \qquad n \rightarrow \infty;
\]
\item[(ii)] $\vartheta$ is asymptotically weakly stable
almost surely, that is 
\[
P_n^*\delta_{(\bx, y)} \overset{\mathsf{d}}{\rightharpoonup}  \vartheta \quad   \forall\,  (\bx, y) \in \b X \times Y_{reg},
\qquad n \rightarrow \infty.
\]
\end{itemize}
\end{thm}
\begin{proof} \textit{Claim (i).} The proof of the first claim relies on Theorem \ref{KS_thm_weak}.
As in the proof of Theorem \ref{asinvmeas},
we choose the Polish space $(S, \mathcal{S})=(\b X \times Y, \mathscr{B}(\b X\times Y))$ with the metric $\rho=\mathsf{d}$
and $Q= \b X \times Y_{s_0}$,
so that from Theorem \ref{asinvmeas} we know that the unique invariant measure $\vartheta$ satisfies $\vartheta(\b X \times Y_{s_0})=1$. Moreover, we observe that the semigroup $P$ is $\mathsf{d}$-Feller in view of Theorem \ref{th:Feller}.
As in the proof of Theorem \ref{asinvmeas}, for $(\bx, y), (\widetilde \bx, \widetilde y)$ in $\b X \times Y_{s_0}$ we consider the measure $\xi:=\xi_{(\bx, y), (\widetilde \bx, \widetilde y)}$ on $(\b X \times Y)^{\mathbb{N}} \times (\b X \times Y)^{\mathbb{N}}$ given by the law of the associated random vector $\left\{\left(\bu(n), \varphi(n) \right),\left( \widetilde {\bu}(n),\widetilde{\varphi}(n)\right)\right\}_{n\in \mathbb{N}}$, where $(\bu, \varphi)$ and $(\widetilde \bu, \widetilde \varphi)$ solve systems \eqref{AC_NS} and \eqref{AC_NS_nud}, respectively, with corresponding initial data $(\bx, y), (\widetilde \bx, \widetilde y)$. In view of Proposition \ref{Girsanov_prop},  the joint law $\xi$ of the pair $(\bu, \varphi)$, $(\widetilde \bu, \widetilde \varphi)$ is a generalised coupling from the class $\hat{\mathcal{C}}\left( \mathbb{P}_{(\bx, y)}, \mathbb{P}_{(\widetilde{\bx}, \widetilde y)}\right)$ which satisfies the additional condition $\pi_1(\xi)=\mathbb{P}_{(\bx, y)}$.
To conclude, it thus remains to prove condition \eqref{coupling_cond_weak} in Theorem \ref{KS_thm_weak}.  To this end, let us fix $\varepsilon>0$.
From Theorem \ref{FP_thm_complete} and  the Markov inequality we infer
\[
\mathbb{P}\left(\|\bu (n)- \widetilde{\bu}(n)\|^2_{\bHs}+ \beta \|\varphi(n)- \widetilde{\varphi}(n)\|^2_{V_1} > \varepsilon\right) 
\le \frac{1}{\varepsilon}\mathbb{E}\left[\|\bu (n)- \widetilde{\bu}(n)\|^2_{\bHs}+ 
\beta \|\varphi(n)- \widetilde{\varphi}(n)\|^2_{V_1} \right]
\le \frac{C}{\varepsilon n^p},
\]
for a $p >1$ and $C$ a positive constant depending on the structural parameters and the initial data. We thus infer 
\[
\lim_{n \rightarrow \infty}\mathbb{P} \left( \|\bu(n)- \widetilde{\bu}(n)\|^2_{\bHs}+ \beta \|\varphi(n)- \widetilde{\varphi}(n)\|^2_{V_1}  \le \varepsilon\right)=1,
\]
which implies \eqref{coupling_cond_weak}.
This concludes the proof of the first claim.
\paragraph{\textit{Claim (ii).}} The proof of the second claim relies on
Theorem \ref{KS_thm_strong} instead. Here, however, we make a 
different choice of the space $(S,\mathcal S)$, by setting
$S=\b X \times Y_{reg}$ and $\mathcal S=\mathscr B(\b X\times Y)\cap (\b X\times Y_{reg})$. 
As far as the metric is concerned, we set $d:=1\wedge\mathsf d_{|\b X\times Y_{reg}}$ as
the restriction to $\b X\times Y_{reg}$ 
of the metric $\mathsf d$ on $\b X\times Y$.
Let us note that $(S,d)$ is not complete now,
and that ${\overline S}^{d}=\b X\times Y$. 
In order to make $S$ a Polish space, 
we need to carefully define the stronger metric $\rho$ on $S$
complying with all the requirements of Section~\ref{KS_sec}.
To this end, for every $(\b x_1, y_1), (\b x_2, y_2)\in \b X\times Y_{reg}=S$ we define
\[
\rho((\b x_1, y_1), (\b x_2, y_2))
:=\|\b x_1-\b x_2\|_{\bHs}
+\|y_1-y_2\|_{V_2}
+\|\tanh^{-1}(y_1)-\tanh^{-1}(y_2)\|_{L^\infty(\OO)}.
\]
Note that $\rho$ is well-defined since $y_1,y_2\in Y_{reg}$, 
hence $|y_1|, |y_2|<1$ almost everywhere in $\OO$. Moreover, 
it is not difficult to show that $\rho$ is indeed complete on $S$
and that $(S,\rho)$ is separable (because $V_2\embed L^\infty(\OO)$), so that $(S,\rho)$ is a Polish space, and that $d$ is continuous with respect to $\rho$. Furthermore, for every
$(\b x,y)\in S$ and
$n\in\enne$ one can set 
$\rho_n^{(\b x,y)}:\overline S^{d}=\b X\times Y\to[0,+\infty)$ as
\[
  \rho_n^{(\b x,y)}(\b z, w):=\rho((\b x,y), (\b z,w_n))=
  \|\b x-\b z\|_{\bHs}+\|y-w_n\|_{V_2} + 
  \|\tanh^{-1}(y)-\tanh^{-1}(w_n)\|_{L^\infty(\OO)},
\]
where $w_n\in V_2$ is the unique solution to the singular perturbation problem  
\[
  \begin{cases}
  w_n-\frac1n\Delta w_n + \frac1n\Psi_{s_0}'(w_n) = w \quad&\text{in } \OO,\\
  w_n=0 \quad&\text{in } \partial \OO.
  \end{cases}
\]
Note that $\rho_n^{(\b x,y)}$ is well-defined since $w_n \in Y_{reg}$.
Indeed, 
for $w\in Y$, the elliptic problem above has a unique solution
$w_n\in V_2$ with $\Psi_{s_0}'(w_n)\in H$, so that 
in particular $w_n \in V_2\cap \mathcal A_{2(s_0+1)}$,
where $V_2\cap \mathcal A_{2(s_0+1)} = Y_{reg}$ since
$2(s_0+1)>2$.
Moreover, if $(\b z,w)\in S=\b X\times Y_{reg}$, 
it is immediate to see by the maximum principle that 
$\sup_n\|w_n\|_{L^\infty(\OO)}\leq\|w\|_{L^\infty(\OO)}<1$,
from which it easily follows that 
$w_n\to w$ in $V_2$ as $n\to\infty$, 
hence also in $L^\infty(\OO)$,
so that $\rho_n^{(\b x,y)}(\b z,w)\to\rho((\b x,y),(\b z,w))$.
Otherwise, if $(\b z,w)\in \overline S^{d}\setminus S 
=\b X\times (Y\setminus Y_{reg})$
instead, then it holds that either $w\in V_1\setminus V_2$
or $\|w\|_{L^\infty(\OO)}=1$.
In the former case one has that $\|w_n\|_{V_2}\to+\infty$
as $n\to\infty$, whereas in the latter case on has that 
$\sup_{n}\|w_n\|_{L^\infty(\OO)}=1$ hence 
$\|\tanh^{-1}(w_n)\|_{L^\infty(\OO)}\to+\infty$.
In both cases, it follows that $\rho_n^{(\b x,y)}(\b z,w)\to+\infty$.
This shows that the sequence $(\rho_n^{(\b x,y)})_n$
satisfies the requirements of Section~\ref{KS_sec}.
We are then ready to apply 
Theorem~\ref{KS_thm_strong}.
To this end, note that $\{(\bu(n), \varphi(n))\}_{n\in\enne}$
is still a Markov chain also on $(S,\mathcal S)$, 
and is $d$-Feller with the choice of the weaker metric $d$.
Moreover, condition \eqref{coupling_cond_strong}
has been already checked in the proof of (i).
Hence, Theorem~\ref{KS_thm_strong} applies and yields
also the statement (ii). This concludes the proof.
\end{proof} 

\begin{remark}
 We notice that the thesis (ii) of Theorem \ref{weak_conv_ACNS} is not stronger than the thesis (i). In fact, in general the weak convergence is weaker that the convergence in the 1-Wasserstein distance (see e.g. \cite[Lemma 5.1.7 and Proposition 7.1.5]{AGSbook}). 
\end{remark}

\section*{Acknowledgments}
The present research has been supported by MUR, 
grant Dipartimento di Eccellenza 2023-2027.
The authors are members of Gruppo Nazionale per l’Analisi Matematica, la Probabilità e le loro Applicazioni (GNAMPA), Istituto Nazionale di Alta Matematica (INdAM), and gratefully acknowledge financial support through the project CUP-E53C22001930001.
The authors gratefully  acknowledge the financial support of the project  ``Prot. P2022TX4FE\_02 -  Stochastic particle-based anomalous reaction-diffusion models with heterogeneous interaction for radiation therapy'' financed by the European Union - Next Generation EU, Missione 4-Componente 1-CUP: D53D23018970001.


\appendix  
\section{Tightness criteria} \label{sec:appA}
In this first appendix, we collect a compactness and a tightness result, starting with the former.
\begin{lem}
\label{tight_cor}
Let $p>2$, let $\gamma \in(\frac 1p,\frac12)$ and 
$\delta\in(0,\frac12)$. For any $T > 0$, the embeddings
\begin{align*}
L^\infty(0,T;V_1) \cap W^{\gamma,p}(0,T;V_1^*)  &\embed 
\C([0,T];H)\cap \C_{\text{w}}([0,T];V_1),\\
L^2(0,T;V_2)\cap W^{\gamma,p}(0,T;V_1^*) &\embed L^2(0,T;V_1),\\
L^\infty(0,T;\bHs) \cap W^{\gamma,p}(0,T;\bVsd)  &\embed \C([0,T];D(\b{A}^{-\delta}))
\cap \C_{\text w}([0,T];\bHs),\\
L^2(0,T;\bVs) \cap \cap W^{\gamma,p}(0,T;\bVsd)  &\embed 
L^2(0,T;\bHs)
\end{align*}
are compact.
\end{lem}
\begin{proof}
Let us only prove the claim for the first two embeddings, as the proof for the last two is analogous. 
The compactness of the embeddings
\[
L^\infty(0,T;V_1)\cap W^{\gamma,p}(0,T;V_1^*)  \embed \C([0,T];H), \qquad
L^2(0,T;V_2)\cap W^{\gamma,p}(0,T;V_1^*)  \embed  L^2(0,T;V_1)
\]
follows from \cite[Corollary 5]{simon}. In order to show that the embedding $L^\infty(0,T;V_1)\cap W^{\gamma,p}(0,T;V_1^*)  \embed  \C_{\text{w}}([0,T];V_1)$ is compact, let $K \subset L^\infty(0,T;V_1)\cap W^{\gamma,p}(0,T;V_1^*)$ be bounded. For any $\alpha \in(0,\gamma-\frac1p)$ one has
the continuous embedding $W^{\gamma,p}(0,T;V_1^*)\embed 
\C^{0, \alpha}(0,T;V_1^*)$, and therefore it 
holds that 
\[M:=\sup_{f\in K}\|f\|_{\C^{0, \alpha}(0,T;V_1^*)}<+\infty,\]
hence
\[
\lim_{\varepsilon \to 0} \sup_{f\in K} 
\sup_{\vert t-s\vert\le \varepsilon} \Vert f(t)-f(s)\Vert_{V_1^*}
\leq M\lim_{\varepsilon \to 0}\varepsilon^\alpha=0.
\] 
This shows that $K$ is uniformly equicontinuous in $\C([0,T];V_1^*)$.
Hence, since $K$ is also bounded in $L^{\infty}(0,T;V_1)$,
we can conclude by applying \cite[Proposition 4.2]{BHW-2019} (see also \cite[Appendix B]{FZ}) and infer that $K$ is relatively compact in $\C_{\text w}([0,T];V_1)$.
\end{proof}
\begin{prop}
\label{prop_tight}
Let $p>2$, let $\gamma\in(\frac 1p,\frac12)$, and $T>0$. Let 
$\{\bu_n\}_{n \in \mathbb{N}}$ and
$\{\varphi_n\}_{n \in \mathbb{N}}$ be sequences of adapted stochastic processes 
with values in $\trajU$ and $\trajF$, respectively, 
such that 
\begin{align*}
&\sup_{n \in \mathbb{N}}
\|\bu_{n}\|_{L^p(\Omega;L^\infty(0,T;\bHs)\cap L^2(0,T;\bVs)\cap W^{\gamma,p}(0,T; \b V_\sigma^*))} <+\infty,\\
&\sup_{n \in \mathbb{N}}
\|\varphi_{n}\|_{L^p(\Omega;L^\infty(0,T;V_1)\cap L^{2}(0,T;V_2)\cap
W^{\gamma,p}(0,T;V_1^*))} <+\infty.
\end{align*}
Then, the family of laws of $\{\bu_n\}_{n \in \mathbb{N}}$ 
is tight on $\trajU$, and the 
family of laws of $\{\varphi_n\}_{n \in \mathbb{N}}$ 
is tight on $\trajF$.
\end{prop}
\begin{proof}
Let us prove the second claim only, as the first one can be proven analogously, replacing the spaces accordingly and exploiting the corresponding estimates. 
For any $R>0$, let $B_R$ be the closed ball of radius $R$ in the space $L^\infty(0,T;V_1) \cap L^2(0,T;V_2)\cap W^{\gamma,p}(0,T;V_1^*)$. By Lemma \ref{tight_cor}, $B_R$ is compact in $\trajF$. The Markov inequality yields, for any $n \in \mathbb{N}$,
\begin{align*}
\mathbb{P}\left(\varphi_n \in B_R^C \right)&=\mathbb{P} \left( \|\varphi_n\|_{L^\infty(0,T;V_1) \cap L^2(0,T;V_2)\cap W^{\gamma,p}(0,T;V_1^*)}>R\right)
\\
&\le \frac{1}{R^p}\mathbb{E}\left[ \|\varphi_n\|^p_{L^\infty(0,T;V_1) \cap L^2(0,T;V_2)\cap W^{\gamma,p}(0,T;V_1^*)}\right]
\le \frac{C}{R^p},
\end{align*}
where
\[
C = \sup_{n \in \mathbb{N}}
\|\varphi_{n}\|_{L^p(\Omega;L^\infty(0,T;V_1)\cap L^{2}(0,T;V_2)\cap
W^{\gamma,p}(0,T;V_1^*))}.
\]
Thus we conclude that, for any $\varepsilon>0$, there exists a compact set $B_R$ of $\trajF$, with $R=R(\varepsilon)>0$, such that 
\begin{equation*}
\sup_{n \in \mathbb{N}} \Law_\P(\varphi_n)(B_R^C) \le \varepsilon.
\end{equation*}
This concludes the proof by definition of tightness.
\end{proof}

    
\section{A priori estimates} \label{sec:appB}
In this second appendix, we provide key a priori estimates on system \eqref{AC_NS}. We always work either under 
Assumptions \ref{hyp:structural}--\ref{hyp:diffusionAC} with $Z=\bVsd$ or with Assumptions \ref{hyp:structural}--\ref{hyp:additionalAC} with $Z=\bHs$, so that the system 
\eqref{AC_NS} admits a unique probabilistically-strong solution 
by Theorem~\ref{th:probstrongsol} and Proposition~\ref{prop:uniqueness}.

\subsection{Moments estimates} Let us start by providing estimates on the moments of the solution of \eqref{AC_NS}.
\begin{prop}
    \label{a_priori_est}
Let either Assumptions \ref{hyp:structural}--\ref{hyp:diffusionAC} hold with $Z=\bVsd$ or Assumptions \ref{hyp:structural}--\ref{hyp:additionalAC} hold with $Z=\bHs$.
Let $(\bx,y) \in \b X \times Y$ and let 
$(\bu,\varphi)$ be the unique strong solution
to the system \eqref{AC_NS} starting from the initial state
$(\b x,y)$. Then, for all $p>2$,
\begin{enumerate}[(i)]
    \item there exists a constant $C_1 = C_1(p,\,L_{G_2},\,L_F,\,|\OO|) > 0$, 
    independent of $\nu$ and $\beta$, such that, 
    for all $t > 0$,
    \[
    \|\varphi\|_{L^p_w(\Omega; L^\infty(0,t;H))} + 
    \|\sqrt \beta \nabla \varphi\|_{L^p_\cP(\Omega; L^2(0,t;\bH))} \leq C_1\left(\sqrt t + \|y\|_H\right);
    \]
    \item there exists a constant $C_2 = C_2(p,\,L_{G_1},\,L_{G_2},\,L_F,\,|\OO|,\,\|F\|_{\C([-1,1]}) > 0$, independent of $\nu$ and $\beta$, such that, for all $t > 0$,
    \begin{multline*}
        \|\bu\|_{L^p_w(\Omega; L^\infty(0,t;\bHs))} + \|\sqrt\nu\bu\|_{L^p_\cP(\Omega; L^2(0,t;\bVs))}+\|\sqrt \beta \nabla \varphi\|_{L^p_w(\Omega; L^\infty(0,t;\bH))}\\  + \|\beta\Delta\varphi\|_{L^p_\cP(\Omega;L^2(0,t;H))}  + \|F'(\varphi)\|_{L^p_\cP(\Omega;L^2(0,t;H)} \leq C_2\left(1 + \sqrt t + \|\bx\|_{\bHs} + \|\sqrt \beta \nabla y\|_\bH + \|y\|_H\right);
    \end{multline*}
    \item for any $T > 0$, there exist constants $C_3 = C_3(T,\,p,\,L_{G_1},\,L_{G_2},\,L_F,\,|\OO|,\,\|F\|_{\C([-1,1])},\,\nu,\,\beta) > 0$ and $\gamma \in \left(\frac 1p, \frac 12\right)$ such that
    \begin{align*}
    \|\bu\|_{L^{\frac p2}(\Omega;W^{\gamma,p}(0,T;\bVsd)} +
    \|\varphi\|_{L^p(\Omega;W^{\gamma,p}(0,T;V_1^*))}
    &\le C_3\left(1 + \|\bx\|_{\bHs}^2 + \|y\|_\bH^2 + \|y\|_{V_1}^2\right).
\end{align*}
\end{enumerate}
\end{prop}
\begin{proof} For the sake of clarity, we divide the argument into steps.
\paragraph{\textit{First estimate}} The It\^o formula for the square $H$-norm of $\varphi$ yields, for any $t \ge 0$, $\mathbb{P}$-a.s.
\begin{multline}
\frac{1}{2}\|\varphi(t)\|_{H}^2 +
\int_0^t\int_\OO\left[\beta|\nabla \varphi(s)|^2 + 
\bu(s)\cdot\nabla\varphi(s)\varphi(s)+
F'(\varphi(s)) \varphi(s) \right] \, \mathrm{d}s 
\\ = \dfrac{1}{2}\|y\|_{H}^2 + \int_0^t \left(\varphi(s), G_{2} (\varphi(s))\,\mathrm{d}W_2(s)\right)_H + \frac{1}{2}\int_0^t \|G_{2}(\varphi(s))\|^2_{\LL_{HS}(U_2, H)} \: \mathrm{d}s.
\label{stima1}
\end{multline} 
Since $\div\bu=0$, integration by parts yields
\[
  \int_\OO\bu(s)\cdot\nabla\varphi(s)\varphi(s)
  =\frac12\int_\OO\bu(s)\cdot\nabla(|\varphi(s)|^2)
  =-\frac12\int_\OO\div\bu(s)|\varphi(s)|^2=0.
\]
Moreover, by means of \eqref{F'_prop} we estimate
\begin{align*}
\int_0^t (F'(\varphi(s)),\varphi(s))_H\, {\rm d}s \ge L_F \int_0^t \|\varphi(s)\|^2_H\, {\rm d}s - 2L_F|\OO|t,
\end{align*}
while from \eqref{HS_norm} we immediately get
\begin{align*}
\frac 12 \int_0^t \|G_2(u(s))\|^2_{\mathcal L_{HS}(U_2,H)}\, {\rm d}s \le \frac{L_{G_2}|\OO|t}{2}.
\end{align*}
On account of the previous observations, raising \eqref{stima1} to the power $\frac p2$, taking the supremum in time and $\P$-expectations we have 
\begin{multline*}
\frac 1{2^{\frac p2}} \E\sup_{r \in [0,t]}\|\varphi(r)\|^p_H
+\E\left(
\int_0^t\left(\beta\|\nabla \varphi(s)\|^2_H+L_F\|\varphi(s)\|_H^2\right)\, {\rm d}s
\right)^{\frac p2}
\\ \leq \dfrac {3^{\frac p2-1}}{2^{\frac p2}} \|y\|^p_H + 3^{\frac p2-1}\left( \frac{L_{G_2}}{2}+2L_F\right)^{\frac{p}2} |\OO|^\frac p2 t^\frac p2
+3^{\frac p2-1}\E\sup_{r \in [0,t]}\left|
\int_0^r (\varphi(s), G_2(\varphi(s))\: \d W_2(s))_H
\right|^{\frac{p}2}.
\end{multline*}
By means of the Burkholder--Davis--Gundy and Young inequality and
again \eqref{HS_norm},  we infer then
\begin{equation}
\label{stima_fin_2}
\E\sup_{r \in [0,t]}\|\varphi(r)\|^p_H+
\E\left(\int_0^t\beta\|\nabla\varphi(s)\|^2_{\b H}\, {\rm d}s \right)^{\frac p2}
\leq C_1\left(t^{\frac p2}+\|y\|^p_H\right),
\end{equation}
where the constant $C_1 > 0$ only depends on $p,\,L_{G_2},\,L_F$ and $|\OO|$. This shows the first estimate.

\paragraph{\textit{Second estimate.}}
We apply the It\^o formula for twice Fréchet-differentiable functionals to free energy functional 
\[
\mathcal{E}(\b{u}, \varphi)
:= \dfrac{1}{2}\|\b{u}\|_{\bHs}^2 + 
\dfrac{\beta}{2}\|\nabla \varphi\|_\bH^2 + \int_\OO F(\varphi).
\]
By exploiting classical Yosida-type approximations on the nonlinearity $F$,
as done in \cite{Bertacco21} in detail (see also \cite[Section 3.7]{DPGS}), we obtain for every $t \ge 0$, and $\mathbb{P}$-almost surely, the energy equality 
\begin{multline}
\dfrac{1}{2}\|\b{u}(t)\|^2_{\bHs} 
+\frac{\beta}2 \|\nabla \varphi(t)\|^2_H
+\|F(\varphi(t))\|_{L^1(\OO)}\\+ \int_0^t\left[
\nu\|\nabla\b{u}(s)\|^2_{\bHs}
+\beta^2\|\Delta \varphi(s)\|^2_H 
+\|F'(\varphi(s))\|_H^2
\right]\, \mathrm{d}s
+2\beta\int_0^t\int_{\OO} 
F^{\prime \prime}(\varphi(s))|\nabla \varphi(s)|^2 \, {\rm d}s\\
=  \dfrac{1}{2}\|\bx\|^2_{\bHs} +
\frac {\beta}2 \|\nabla y\|^2_H +\|F(y)\|_{L^1(\OO)}
+ \dfrac{1}{2}\int_0^t \|G_{1}(\b{u}(s))\|^2_{\LL_{HS}(U_1, \bHs)} \: \mathrm{d}s\\
\frac 12 \int_0^t \sum_{k\in\enne}\int_{\OO}
\left[\beta|g_k'(\varphi(s))\nabla \varphi(s)|^2 + 
F''(\varphi(s))|g_k(\varphi(s))|^2\right]\,\d s \\
+\int_0^t\left(\bu(s), G_{1}( \bu(s))\,\d  W_1(s)\right)_{\bHs}
+ \int_0^t \left( -\beta\Delta \varphi(s) + F^\prime (\varphi(s))
, G_2(\varphi(s))\,{\rm d}W_2(s)\right)_H. 
\label{stima6}
\end{multline}
Let us estimate the various terms appearing in \eqref{stima6}.
By Assumption \ref{hyp:potential} we have that
\begin{align*}
2\beta\int_0^t\int_\OO F''(\varphi(s))|\nabla \varphi(s)|^2 \, {\rm d}s \ge -2\beta L_F\int_0^t\|\nabla \varphi(s)\|_\bH^2\,\d s.
\end{align*}
Moreover, thanks to Assumption \ref{hyp:diffusionNS} we get 
\begin{equation*}
 \dfrac{1}{2}\int_0^t \|G_{1}(\b{u}(s))\|^2_{\LL_{HS}(U_1, \bHs)} \: \mathrm{d}s
 \le \frac{L_{G_1}t}{2}.
\end{equation*}
while from Assumption \ref{hyp:diffusionAC} we infer 
\[
\frac 12 \int_0^t \sum_{k\in\enne}\int_{\OO}
\left[\beta|g_k'(\varphi(s))\nabla \varphi(s)|^2 + 
F''(\varphi(s))|g_k(\varphi(s))|^2\right]\,\d s
 \le \frac{L_{G_2}}{2}\int_0^t 
 \beta\|\nabla \varphi(s)\|^2_H\, {\rm d}s
+\frac{L_{G_2}|\OO|t}{2}.
\]
Assumption \ref{hyp:diffusionNS}, jointly with the Burkholder-Davis-Gundy and the Young inequalities yields
\begin{align*}
\E\sup_{r \in [0,t]}
\left|
\int_0^r\left(\bu(s), G_{1}( \bu(s))\,\d  W_1(s)\right)_{\bHs}
\right|^{\frac p2}
&\le \frac 1{2^{{\frac p2}+1}} \mathbb{E} \sup_{r \in [0,t]} \|\bu(r)\|^p_\bHs +C_2
L_{G_1}^{\frac p2}|\OO|^{\frac p2}t^{\frac p2}
\end{align*}
and, on account of \eqref{HS_norm},
\begin{multline*}
\E\sup_{r \in [0,t]}
\left|\int_0^r \left( 
-\beta\Delta \varphi(s)+
F^\prime (\varphi(s)), G_2(\varphi(s))\,{\rm d}W_2(s)\right)_H 
\right|^{\frac p2}
\\ \le \frac 12 \mathbb{E} 
\left(\int_0^t 
\left(\beta^2\|\Delta \varphi(s)\|_H^2
+\|F^\prime(\varphi(s))\|^2_{H}\right)\,{\rm d}s\right)^{\frac p2} +C_2L_{G_2}^\frac p2|\OO|^\frac p2,
\end{multline*}
where the constant $C_2 > 0$ depends on $p,\,L_{G_1},\,L_{G_2},\,L_F,\,|\OO|$ and $\|F\|_{\C([-1,1])}$, while is independent of $\beta$.
Hence, by collecting the above results in \eqref{stima6}, raising it to the power $\frac p2$, taking the supremum in time and $\P$-expectations, rearranging the terms we obtain
\begin{equation*}
\begin{split}
& \E\sup_{r\in[0,t]}\left(\|\bu(r)\|_{\bHs}^p
     +\|\sqrt \beta \nabla\varphi(r)\|_{\b H}^p\right)
    +\E\left(\int_0^t  \left( 
    \|\sqrt \nu\bu(s)\|_{\bVs}^2
    +\|\beta\Delta \varphi(s)\|_{H}^2+
    \|F'(\varphi(s))\|_H^2\right)\, \d s\right)^{\frac p2}\\ & \hspace{2cm}\leq C_2\left( 1+t^\frac p2+  \|\bx\|^p_{\bHs} 
+\|\sqrt \beta\nabla y\|^p_\bH+
\mathbb{E}\left(\int_0^t 
\|\sqrt \beta\nabla \varphi(s)\|^2_H\, {\rm d}s\right)^{\frac p2}\right) \\
& \hspace{2cm} \leq C_2\left( 1+t^\frac p2+  \|\bx\|^p_{\bHs} 
+\|\sqrt \beta\nabla y\|^p_\bH+\|y\|^p_H\right),
\label{stima_fin_3}
\end{split}
\end{equation*}
where we also used the first estimate and possibly updated the constant $C_2$. This shows the second estimate.

\paragraph{\textit{Third estimate}.}
In what follows, we let $T > 0$ to be a fixed final time, and the symbol $C$ denotes a constant possibly depending on $p,\,L_{G_1},\,L_{G_2},\,L_F,\,|\OO|,\|F\|_{\C([-1,1])},\,\nu$ and $\beta$, whose value may be updated from line to line.
First, note that a comparison argument and the first claim imply that 
\begin{equation}
\label{es5}
\|w\|_{L^p_\cP(\Omega;L^{2}(0,T;H))}
\le C\left(1 + \|\bx\|_{\bHs} + \|\nabla y\|_\bH + \|y\|_H\right).
\end{equation}
Secondly, note that Assumptions \ref{hyp:diffusionNS} and \ref{hyp:diffusionAC}, jointly again the first claim yield	\begin{align*}
			\|G_{1}(\b{u})\|_{L^\infty_\cP(\Omega\times(0,T);\cL^2(U_1,\bHs)))}+
			\|G_{2}(\varphi)\|_{L^\infty_\cP(\Omega \times (0,T);\cL^2(U_2,H))} 
   &\leq C,
		\end{align*}
 so that as a consequence of \cite[Lemma 2.1]{fland-gat},
	 the estimate on the It\^{o} integrals
		\begin{align} 
			\left\| \int_0^\cdot G_{1}(\b{u}(s))\,\mathrm{d}W_1(s)
			\right\|_{L^q_\cP(\Omega;W^{k,q}(0,T;\bHs))}+
			\label{eq:estimate10}
			\left\| \int_0^\cdot G_{2}(\varphi(s))\,\mathrm{d}W_2(s)
			\right\|_{L^q_\cP(\Omega; W^{k,q}(0,T; H))}
			&\leq C
		\end{align}
		 holds for every $k \in (0,\frac{1}{2})$ and $q\geq 2$, provided that $C$ may depend on $k$ and $q$ as well.
In order to prove the third estimate, let us interpret the equation for $\varphi$ as an equality in $V^*_1$, i.e., for every $t\in[0,T]$, $\P$-a.s.		
    \begin{align*}
		\langle \varphi(t),\psi\rangle_{V_1^*, V_1} &= -
		\int_0^t\!\int_\OO\left[ \bu(s)\cdot\nabla \varphi(s)
		+  w(s)\right]\psi\,\d s + \langle y,\psi\rangle_{V_1^*, V_1} +
		\left(\int_0^t G_{2}( \varphi(s))\,\d  W_2(s), \psi\right)_{H}
    \end{align*}
for all $\psi \in V_1$ such that $\|\psi\|_{V_1} = 1$. The H\"{o}lder inequality yields
		\begin{equation*}
			\left| \int_\OO \psi\bu(s)\cdot\nabla \varphi(s)	
        +  \psi w(s) \right| \leq  \|\bu(s)\|_\bHs
			\|\nabla \varphi(s)\|_{\b{L}^4(\OO)}
			+ \|w(s)\|_H,
		\end{equation*}
which implies thanks once again to the first claim, interpolation, and \eqref{es5} that
		\[
			\left\| \int_0^\cdot\!\int_\OO\left[ \bu(s)\cdot\nabla \varphi(s)
			+  w(s)\right]\d s\right\|_{L^p_\cP(\Omega; H^1(0,T;V_1^*))} 
   \leq {C}\left(1 + \|\bx\|_{\bHs}^2 + \| \nabla y\|_\bH^2 + \|y\|_H^2\right).
		\]
		Then, recalling that $H^1(0,T; V_1^*)\embed W^{\gamma,p}(0,T; V_1^*)$
  for all $\gamma\in(\frac1p, \min\{\frac2p, \frac12\})$,
  by \cite[Lemma 2.1]{DPGS}, a
  comparison argument in the equation for $\varphi$ yields directly
  by using \eqref{eq:estimate10} that
		\begin{equation} \label{eq:estimate11}
\|\varphi\|_{L^p_\cP(\Omega;W^{\gamma,p}(0,T;V_1^*))} \le 
   C_3\left(1 + \|\bx\|_{\bHs}^2 + \|\nabla y\|_\bH^2 + \|y\|_H^2\right),
		\end{equation}
        where the constant $C_3 > 0$ depends on $T,\, p,\,L_{G_1},\,L_{G_2},\,L_F,\,|\OO|,\,\|F\|_{\C([-1,1])}$ and $\beta$.
Similarly, we consider the weak formulation of the Navier--Stokes equations, i.e., for any $t \in [0,T]$, $\mathbb{P}$-a.s.,
\begin{align*}
		 \langle \bu(t),\bv\rangle_{\bVsd,\bVs} &= -
			\int_0^t\left[\nu\ip{\b A \bu(s)}{\bv}_{\b V_\sigma^*, \b V_\sigma}
			+\ip{\b B( \bu(s),  \bu(s))}{\bv}_{\b V_\sigma^*, \b V_\sigma}
			\right]\,\d s \\
			&\quad-\int_0^t\int_\OO w(s)\nabla \varphi(s)\cdot\bv\,\d s
			+ (\bx,\bv)_{\bHs} +
			\left(\int_0^t G_{1}( \bu(s))\,\d  W_1(s), \bv\right)_{\bHs},
		\end{align*}
		for all $\b{v} \in \bVs$ such that $\|\b{v}\|_{\bVs} = 1$.
From the first claim, the continuity of $\b{A}$,
as well as the well-known inequality
		\[
			\|\b{B}( \bu,  \bu)\|_{\bVsd}
			\leq C\|\bu\|_\bHs\|\bu\|_\bVs,
		\]
we infer that
		\begin{multline*}
		    \left\|\int_0^\cdot
  \left[\nu\ip{\b A \bu(s)}{\bv}_{\b V_\sigma^*, \b V_\sigma}
			+\ip{\b B( \bu(s),  \bu(s))}{\bv}_{\b V_\sigma^*, \b V_\sigma}
			\right]\,\d s
  \right\|_{L^p_\cP(\Omega; H^1(0,T;\bVsd))}
		\\ \leq C\left(1 + \|\bx\|^2_{\bHs} + \|\nabla y\|^2_\bH + \|y\|_H^2\right).
		\end{multline*}
		Furthermore, since by the H\"{o}lder and Gagliardo--Nirenberg 
		inequalities, we have
		\[
				\left|\int_\OO w(s)\nabla\varphi(s)\cdot \b{v}
				\right|  \leq \|w(s)\|_H
				\|\nabla\varphi(s)\|_{\b{L}^3(\OO)}\|\bv\|_{\b{L}^6(\OO)}
				\leq		C\|w(s)\|_H\|\varphi\|_{L^6(\OO)}^\frac{1}{2}
				\|\varphi(s)\|_{V_2}^\frac{1}{2},
		\]
thanks to the first claim and \eqref{es5} we get
		\[
		\left\|\int_0^\cdot w(s)\nabla \varphi(s)\,\d s
		\right\|_{L^\frac p2(\Omega; W^{1,\frac{4}{3}}(0,T;\bVsd))} 
  \leq C\left(1 + \|\bx\|_{\bHs} + \|\nabla y\|_\bH + \|y\|_H\right).
		\]
  Recalling that $W^{1,\frac43}(0,T; \b V_\sigma^*)
  \embed W^{\gamma,p}(0,T; \b V_\sigma^*)$
  for all $\gamma\in(\frac1p, \min\{\frac4{3p}, \frac12\})$,
  by \cite[Lemma 2.1]{DPGS},
	summing up and owing to \eqref{eq:estimate10} we arrive at
		\begin{equation*} \label{eq:estimate12}
			\|\bu\|_{L^{\frac p2}(\Omega;
			W^{\gamma,p}(0,T;\bVsd))} \leq C_3 \left(1 + \|\bx\|_{\bHs} + \|\nabla y\|_\bH + \|y\|_H\right),
		\end{equation*}
    where the constant $C_3 > 0$ depends on $T,\,p,\,L_{G_1},\,L_{G_2},\,L_F,\,|\OO|,\,\|F\|_{\C([-1,1])},\,\nu$ and $\beta$.
    The proof is complete.
\end{proof}

\begin{prop}
\label{a_priori_est2}
Let either Assumptions \ref{hyp:structural}--\ref{hyp:diffusionAC} hold with $Z=\bVsd$ or Assumptions \ref{hyp:structural}--\ref{hyp:additionalAC} hold with $Z=\bHs$.
Let $(\bx,y) \in \b X \times Y$ and let 
$(\bu,\varphi)$ be the unique strong solution
to the system \eqref{AC_NS} starting from the initial state
$(\b x,y)$. Then, for every $q\in\mathbb N_+$, 
there exists a
constant $C_4=C_4(q,\lambda_1, L_{G_1},L_{G_2},L_F,|\OO|,\|F\|_{\C([-1,1])})>0$,
independent of $\nu$ and $\beta$, such that
\begin{align*}
&\sup_{t\geq0}\E\|\bu(t)\|_\bHs^{2q}+
\sup_{t\geq0}\E\|\varphi(t)\|_H^{2q}+
\sup_{t\geq0}\E\|\sqrt\beta\nabla\varphi(t)\|_{\bH}^{2q}\leq C_4\left(1+\nu^{-q} +
\|\bx\|_{\bHs}^{2q}+
\|y\|_H^{2q}+\|\sqrt \beta\nabla y\|_{\b H}^{2q}\right).
\end{align*}
\end{prop}
\begin{proof}
    Let us show the case $q=1$ first. Going back to \eqref{stima6} and dropping the positive terms involving the potential $F$ to the left hand side, proceeding as in the proof of Proposition \ref{a_priori_est} yield
    \begin{multline*}
\dfrac{1}{2}\|\b{u}(t)\|^2_{\bHs} 
+\frac\beta2 \|\nabla \varphi(t)\|^2_{\b H}
+ \int_0^t\left[
\nu\|\nabla\b{u}(s)\|^2_{\bHs}
+\beta^2\|\Delta \varphi(s)\|^2_H 
\right]\, \mathrm{d}s\\
\leq  \dfrac{1}{2}\|\bx\|^2_{\bHs} 
+\frac {\beta}2 \|\nabla y\|^2_H +\|F(y)\|_{L^1(\OO)}
+ \dfrac 12 \left( L_{G_1} + L_{G_2}|\OO|\right)t + \frac{L_{G_2}}{2}\int_0^t \beta\|\nabla \varphi(s)\|^2_H\, {\rm d}s\\
+\int_0^t\left(\bu(s), G_{1}( \bu(s))\,\d  W_1(s)\right)_{\bHs}
+ \int_0^t \left( w(s)
, G_2(\varphi(s))\,{\rm d}W_2(s)\right)_H. 
\end{multline*}
Next, from the estimate \eqref{stima1}, using all assumptions yields 
\begin{multline*} 
\frac{1}{2}\|\varphi(t)\|_{H}^2 +
\int_0^t\beta\|\nabla\varphi(s)\|_\bH^2\,\d s+
L_F\int_0^t\|\varphi(s)\|_H^2\,\d s\\ 
\leq \dfrac{1}{2}\|y\|_{H}^2 +
(2L_F + L_{G_2})|\OO|t
+\int_0^t \left(\varphi(s), G_{2} (\varphi(s))\,\mathrm{d}W_2(s)\right)_H.
\end{multline*}
Hence, by multiplying the second inequality by $L_{G_2}$
and by summing up, one infers that 
    \begin{multline*}
\dfrac{1}{2}\|\b{u}(t)\|^2_{\bHs} 
+\frac{\beta}{2} \|\nabla \varphi(t)\|^2_H
+\frac{L_{G_2}}{2}\| \varphi(t)\|^2_H
\\ 
+ \int_0^t\left[
\nu\|\nabla\b{u}(s)\|^2_{\bHs}
+\beta^2\|\Delta \varphi(s)\|^2_H 
+L_{G_2}\beta\|\nabla\varphi(s)\|_H^2
+L_FL_{G_2}\|\varphi(s)\|_H^2
\right]\, \mathrm{d}s
\\ 
\leq  \dfrac{1}{2}\|\bx\|^2_{\bHs} 
+\frac {\beta}{2} \|\nabla y\|^2_H 
+\|F(y)\|_{L^1(\OO)}
+ \frac{L_{G_2}}{2}\|y\|_H^2
+\dfrac12\left[
(L_{G_1} + L_{G_2}|\OO|)
+{L_{G_2}}(2L_F + L_{G_2})|\OO|\right]t
\\ 
+\int_0^t\left(\bu(s), G_{1}( \bu(s))\,\d  W_1(s)\right)_{\bHs}
+ \int_0^t \left( L_{G_2}\varphi(s)+w(s)
, G_2(\varphi(s))\,{\rm d}W_2(s)\right)_H.
\end{multline*}
We deduce that there exists a constant $C > 0$, depending only 
on $L_{G_1},\,L_{G_2},\,L_F,\,|\OO|,\,\|F\|_{\C([-1,1])}$, such that 
for all $t\geq0$
\begin{multline}
    \label{star_phi}
\|\b{u}(t)\|^2_{\bHs} 
+\|\varphi(t)\|^2_{H} + \beta \|\nabla \varphi(t)\|^2_{\b H}
+\int_0^t\left[
\nu\|\nabla\b{u}(s)\|^2_{\bHs}
+\beta^2\|\Delta \varphi (s)\|^2_H
+\|\varphi(s)\|^2_{H} + \beta \|\nabla \varphi(s)\|^2_{\b H}
\right]\, \mathrm{d}s
\\
\leq  \|\bx\|^2_{\bHs} + \|y\|_H^2 + \beta \|\nabla y\|_{\b H}^2+ C\left[  1
+t
+\int_0^t\left(\bu(s), G_{1}( \bu(s))\,\d  W_1(s)\right)_{\bHs} \right. 
\\
\left. + \int_0^t \left(\varphi(s)+w(s)
, G_2(\varphi(s))\,{\rm d}W_2(s)\right)_H \right].
\end{multline}
Taking expectations and differentiating \eqref{star_phi} in time, using the Poincaré inequality as well, we infer that 
\begin{multline*}
  \frac{\d}{\d t}\left(
  \E\|\bu(t)\|_{\bHs}^2
  +\E\|\varphi(t)\|^2_{H} + \beta\E
  \|\nabla \varphi(t)\|^2_{\b H}\right)
   + \left(\lambda_1\nu\E\|\bu(t)\|_{\bHs}^2 +
\E\|\varphi(t)\|^2_{H} + \beta\E \|\nabla \varphi(t)\|^2_{\b H} 
\right)
  \leq 
   C
\end{multline*}
for all $t > 0$, so that the Gronwall lemma implies
\begin{align*}
  \E\|\bu(t)\|_{\bHs}^2
  +\E\|\varphi(t)\|^2_{H} + \beta\E  \|\nabla \varphi(t)\|^2_{\b H}
     &\leq \left[\|\bx\|_{\bHs}^2+\|y\|_{H}^2 + 
   \beta\|\nabla y\|^2_{\b H}\right]e^{-((\lambda_1\nu)\wedge1) t}
  +C\frac{(1-e^{-((\lambda_1\nu)\wedge1) t})}{(\lambda_1\nu)\wedge1}
  \\ &\leq \left[\|\bx\|_{\bHs}^2+\|y\|_{H}^2 + 
   \beta\|\nabla y\|^2_{\b H}\right]
   +C(1+(\lambda_1\nu)^{-1}),
\end{align*}
for all $t \geq 0$ and therefore the claim is shown with $q=1$. 
Let us show now the claim for a general $q+1$, with $q\in\enne_+$
is fixed.
To this end, consider the real-valued process
\[
t\mapsto \|\bu(t)\|_{\bHs}^2+\|\varphi(t)\|_{H}^2 + \beta \| \nabla \varphi(t)\|^2_{\b H}.
\]
Applying the It\^o formula for real-valued processes
to its $(q+1)$th-power formally yields, by performing the analogous computations as above, that
    \begin{multline}
    \label{stima_star}
\left(\|\bu(t)\|_{\bHs}^2+\|\varphi(t)\|_{H}^2 +
\beta \| \nabla \varphi(t)\|^2_{\b H}\right)^{q+1}
\\ +(q+1)\int_0^t\left(\|\bu(s)\|_{\bHs}^2+\|\varphi(s)\|_{H}^2 + \beta \| \nabla \varphi(s)\|^2_{\b H}\right)^{q}
\left[
 \nu\|\nabla\b{u}(s)\|^2_{\bHs}
+\|\varphi(s)\|_{H}^2 + \beta \| \nabla \varphi(s)\|^2_{\b H}
\right]\, \mathrm{d}s
\\ +(q+1)\int_0^t\left(\|\bu(s)\|_{\bHs}^2+\|\varphi(s)\|_{H}^2 + \beta \| \nabla \varphi(s)\|^2_{\b H}\right)^{q}
\beta^2\| \Delta \varphi(s)\|^2_{H}\, \mathrm{d}s
 \\
\leq \left(\|\bx\|^2_{\bHs} + \|y\|_{H}^2 + \beta\|\nabla y\|^2_{\b H}\right)^{q+1}
+C\left[(q+1)\int_0^t\left(\|\bu(s)\|_{\bHs}^2+\|\varphi(s)\|_{H}^2 + \beta \| \nabla \varphi(s)\|^2_{\b H}\right)^{q}\,\d s\right.
 \\
+(q+1)\int_0^t\left(\|\bu(s)\|_{\bHs}^2+\|\varphi(s)\|_{H}^2 + \beta \| \nabla \varphi(s)\|^2_{\b H}\right)^q
\left(\bu(s), G_{1}( \bu(s))\,\d  W_1(s)\right)_{\bHs}
 \\ 
+ (q+1)\int_0^t \left(\|\bu(s)\|_{\bHs}^2+\|\varphi(s)\|_{H}^2 + \beta \| \nabla \varphi(s)\|^2_{\b H}\right)^{q}
\left(\varphi(s)-\beta \Delta \varphi(s)
, G_2(\varphi(s))\,{\rm d}W_2(s)\right)_H
\\
+\frac{1}2(q+1)q\int_0^t
\left(\|\bu(s)\|_{\bHs}^2+\|\varphi(s)\|_{H}^2 + 
\beta \| \nabla \varphi(s)\|^2_{\b H}\right)^{q-1}
\|(\bu(s), G_1(\bu(s)))_\bHs\|_{\mathcal L_{HS}(U_1,\mathbb R)}^2\,\d s
 \\
\left.+\frac{1}2(q+1)q\int_0^t
\left(\|\bu(s)\|_{\bHs}^2+\|\varphi(s)\|_{H}^2 + 
\beta \| \nabla \varphi(s)\|^2_{\b H}\right)^{q-1}
\|(\varphi(s)-\beta \Delta \varphi(s), G_2(\varphi(s)))_H\|_{\mathcal L_{HS}(U_2,\mathbb R)}^2\,\d s\right].
\end{multline}
By the Young inequality we have 
\begin{multline*}
  C(q+1)\E\int_0^t\left(\|\b{u}(s)\|^2_{\bHs} 
+\|\varphi(s)\|^2_{H}
+\beta\|\nabla\varphi(s)\|^2_{\b H}\right)^{q}\,\d s \\
\leq 
\frac{(q+1)((\lambda_1\nu)\wedge1)}4
\E\int_0^t\left(\|\b{u}(s)\|^2_{\bHs} 
+\|\varphi(s)\|^2_{H}
+\beta\|\nabla\varphi(s)\|^2_{\b H}\right)^{q+1}\,\d s
+\frac{4^q(q+1)C^{q+1}}{((\lambda_1\nu)\wedge1)^{q}}t.
\end{multline*}
Thanks to Assumption \ref{hyp:diffusionNS}, we have
\begin{align*}
&\frac C2(q+1)q\int_0^t
\left(\|\bu(s)\|_{\bHs}^2+\|\varphi(s)\|_{H}^2 
+ \beta \| \nabla \varphi(s)\|^2_{\b H}\right)^{q- 1}
\|(\bu(s), G_1(\bu(s)))_\bHs\|_{\mathcal L_{HS}(U_1,\mathbb R)}^2\,\d s\\
&\qquad \leq \frac C2(q+1)q L_{G_1}
\int_0^t
\left(\|\bu(s)\|_{\bHs}^2+\|\varphi(s)\|_{H}^2 + \beta
\| \nabla \varphi(s)\|^2_{\b H}\right)^{q- 1}
\|\bu(s)\|_{\bHs}^2\,\d s\\
&\qquad\leq \frac C2(q+1)q L_{G_1}
\int_0^t
\left(\|\bu(s)\|_{\bHs}^2+\|\varphi(s)\|_{H}^2 + \beta \| \nabla \varphi(s)\|^2_{\b H}\right)^{q}\,\d s\\
&\qquad\leq\frac{(q+1)((\lambda_1\nu)\wedge1)}4
\E\int_0^t\left(\|\b{u}(s)\|^2_{\bHs} 
+\|\varphi(s)\|^2_{H}
+\beta\|\nabla\varphi(s)\|^2_{\b H}\right)^{q+1}
+\frac{(2q)^q(q+1)C^{q+1}}{2((\lambda_1\nu)\wedge1)^{q}}t,
\end{align*}
while thanks to Assumption \ref{hyp:diffusionAC} one has, 
possibly updating the value of the constant $C$,
for all $s\in[0,t]$,
\begin{align}
\label{stima_traccia}
&
\|(\varphi(s)-\beta \Delta \varphi(s), G_2(\varphi(s)))_H\|_{\mathcal L_{HS}(U_2,\mathbb R)}^2
\notag \\
&\qquad \leq  2
\sum_{k\in\enne}
\left[
\left(\int_\OO
\varphi(s)g_k(\varphi(s))
\right)^2+
\left( \beta \int_\OO
g'_k(\varphi(s))|\nabla\varphi(s)|^2
\right)^2
\right]
\notag \\
&\qquad \leq  2 L_{G_2}
\left(|\OO|^2
+ \beta^2 \|\nabla\varphi(s)\|_H^4
\right)
\leq C(1+\beta^2\|\Delta\varphi(s)\|_H^2),
\end{align}
so that the Young inequality yields, by possibly updating again the value of the constant $C$ also depending on $q$,
\begin{multline*}
\frac C2(q+1)q\int_0^t
\left(\|\bu(s)\|_{\bHs}^2+\|\varphi(s)\|_{H}^2 + \beta \| \nabla \varphi(s)\|^2_{\b H}\right)^{q-1}
\|(\varphi(s)-\beta\Delta\varphi(s), G_2(\varphi(s)))_H\|_{\mathcal L_{HS}(U_2,\mathbb R)}^2\,\d s\\
\leq 
\frac{q+1}2\int_0^t\left(\|\bu(s)\|_{\bHs}^2+\|\varphi(s)\|_{H}^2 + \beta \| \nabla \varphi(s)\|^2_{\b H}\right)^{q}
\beta^2\| \Delta \varphi(s)\|^2_{H}\, \mathrm{d}s +
C^2\int_0^t
\beta^2\| \Delta \varphi(s)\|^2_{H}\, \mathrm{d}s.
\end{multline*}
This implies, by taking \eqref{star_phi} into account, 
by taking expectations,
and by possibly updating the value of $C$, that
    \begin{align*}
&\left(\E\|\b{u}(t)\|^{2}_{\bHs} 
+\E\|\varphi(t)\|^{2}_{H}+
\beta\E\|\nabla \varphi(t)\|^{2}_{\b H}\right)^{q+1}\\
&\qquad+\frac{q+1}2((\lambda_1\nu)\wedge1)\int_0^t
\left(\E\|\b{u}(s)\|^{2}_{\bHs} 
+\E\|\varphi(s)\|^{2}_{H}
+\beta\E\|\nabla\varphi(s)\|^{2}_{\b H}\right)^{q+1}\, \mathrm{d}s\\
&\leq  C\left(1+\|\bx\|^{2}_{\bHs} + 
\|y\|_{H}^{2}+
\beta\|\nabla y\|_{\b H}^{2}\right)^{q+1}
+\frac{C}{(\lambda_1\nu)^q\wedge1}t.
\end{align*}
It follows for every $t\geq0$ that 
\begin{align*}
&\left(\E\|\b{u}(t)\|^{2}_{\bHs} 
+\E\|\varphi(t)\|^{2}_{H}+
\beta\E\|\nabla \varphi(t)\|^{2}_{\b H}\right)^{q+1}\\
&\leq C\left(1+\|\bx\|^{2}_{\bHs} + 
\|y\|_{H}^{2}+
\beta\|\nabla y\|_{\b H}^{2}\right)^{q+1}
e^{-\frac{q+1}2((\lambda_1\nu)\wedge 1)t}
+ \frac {2C}{(q+1)((\lambda_1\nu)^{q+1}\wedge 1)}
(1-e^{-\frac{q+1}2((\lambda_1\nu)\wedge 1)t})\\
&\leq C\left(1+\|\bx\|^{2}_{\bHs} + 
\|y\|_{H}^{2}+
\beta\|\nabla y\|_{\b H}^{2}\right)^{q+1}
+ \frac {2C}{q+1}(1+(\lambda_1\nu)^{-(q+1)}),
\end{align*}
and this concludes the proof.
\end{proof}

\begin{prop}
    \label{a_priori_est2_bis}
Let either Assumptions \ref{hyp:structural}--\ref{hyp:diffusionAC} hold with $Z=\bVsd$ or Assumptions \ref{hyp:structural}--\ref{hyp:additionalAC} hold with $Z=\bHs$.
Let $(\bx,y), (\widetilde{\b x}, \widetilde y) \in \b X \times Y$,
let $N\in\enne_+$, $\eta=\frac{\lambda_N\nu}2$, and let 
$(\widetilde \bu,\widetilde \varphi)$ be the unique strong solution
to the system \eqref{AC_NS_nud} starting from the initial state
$(\widetilde{\b x},\widetilde y)$.
Then, for every $q\in\mathbb N_+$, 
there exists a
constant $C_4=C_4(N,
q, \lambda_1, L_{G_1},L_{G_2},L_F,|\OO|,\|F\|_{\C([-1,1])})>0$,
independent of $\nu$ and $\beta$, such that
\begin{align*}
\sup_{t\geq0}\E\|\widetilde\bu(t)\|_\bHs^{2q}&+
\sup_{t\geq0}\E\|\widetilde\varphi(t)\|_H^{2q}+
\sup_{t\geq0}\E\|\sqrt\beta\nabla\widetilde\varphi(t)\|_{\bH}^{2q}\\
&\leq C_4\left(1+\nu^{-q} +
\|\widetilde \bx\|_{\bHs}^{2q}
+\|\widetilde y\|_H^{2q}+\|\sqrt \beta\nabla\widetilde y\|_{\b H}^{2q}
+\|\bx\|_{\bHs}^{2q}+
\|y\|_H^{2q}+\|\sqrt \beta\nabla y\|_{\b H}^{2q}\right).
\end{align*}
\end{prop}
\begin{proof}
    The proof follows by arguing as in Proposition~\ref{a_priori_est2}, by noting that the nudged term on
    the right-hand yields 
    \[
    \frac{\lambda_N\nu}2\int_0^t(P_N(\bu(s)-\widetilde\bu(s),
    \widetilde\bu(s))_\bH\,\d s
    \leq \frac{\lambda_N\nu}4\int_0^t\|\bu(s)\|_\bH^2\,\d s
    -\frac{\lambda_N\nu}4\int_0^t
    \|P_N\widetilde\bu(s)\|_\bH^2\,\d s,
    \]
    and by exploiting \eqref{star_phi} as well as
    the estimate of Proposition~\ref{a_priori_est2} 
    (see also \cite[Lemma A.3]{FZ23}).
\end{proof}

\subsection{Estimates in probability}
We start with the following four preliminary lemmata. The first three provide some probability estimates for suitable functionals of the solution process. The last one is a general result that holds true for any martingale process.
\begin{lem}
\label{lemma1}
Let either Assumptions \ref{hyp:structural}--\ref{hyp:diffusionAC} hold with $Z=\bVsd$ or Assumptions \ref{hyp:structural}--\ref{hyp:additionalAC} hold with $Z=\bHs$.
Let $(\bx,y) \in \b X \times Y$ and let 
$(\bu,\varphi)$ be the unique strong solution
to the system \eqref{AC_NS} starting from the initial state
$(\b x,y)$. Then, 
there exists a constant 
$\mathfrak{C}_1=\mathfrak{C}_1(\lambda_1,L_{G_1},L_{G_2},L_F,|\OO|,\|F\|_{\C([-1,1])})>0$,
independent of $\nu$ and $\beta$,
such that 
\begin{align}
\label{a_priori_estimate3}
    &\nu\int_0^t\|\bu(\tau)\|_{\bVs}^2\|\bu(\tau)\|_\bHs^2\,\d \tau
    \leq \mathfrak{C}_1\left(1 + (1+\nu^{-2}) t + \|\bx\|_{\bHs}^4 
    + \|y\|_{H}^4 +\beta^2\|\nabla y\|_{\b H}^4 
 + \mathfrak{M}_1(t)
    +\mathfrak{M_2}(t)\right)
\end{align}
for all $t\geq0$, $\P$-almost surely, 
where the martingales 
$\mathfrak{M}_1$ and $\mathfrak{M}_2$ are given by
\begin{align*}
    \mathfrak{M}_1(t)&:=
    \int_0^t\left({1+}\|\b{u}(s)\|^2_{\bHs} 
+\|\varphi(s)\|_{H}^2
+\beta\|\nabla \varphi(s)\|^2_{\b H}\right)
\left(\bu(s), G_{1}( \bu(s))\,\d  W_1(s)\right)_{\bHs}\,,\quad t\geq0\,,\\
\mathfrak{M}_2(t)&:=\int_0^t
\left({1+}\|\b{u}(s)\|^2_{\bHs} 
+\|\varphi(s)\|_{H}^2
+\beta\|\nabla \varphi(s)\|^2_{\b H}\right)
\left( w(s)+\varphi(s), 
G_2(\varphi(s))\,{\rm d}W_2(s)\right)_H\,, \quad t\geq0\,.
\end{align*}
\end{lem}
\begin{proof}
    By proceeding as in the proof of Proposition~\ref{a_priori_est2}, from \eqref{stima_star}
    with $q+1=2$, we infer that there exists a constant $C>0$, depending on $L_{G_1},L_{G_2},L_F$, such that 
        \begin{align*}
&\left(\|\b{u}(t)\|^2_{\bHs} 
+\|\varphi(t)\|_H^2
+\beta\|\nabla \varphi(t)\|^2_H\right)^2\\
&\qquad+ \int_0^t
\left(\|\b{u}(s)\|^2_{\bHs} 
+\|\varphi(s)\|_{H}^2
+\beta\|\nabla \varphi(s)\|^2_{\b H}\right)
\left(\nu\|\nabla\b{u}(s)\|^2_{\bHs}
+\|\varphi(s)\|_{H}^2
+\beta\|\nabla\varphi(s)\|_{\b H}^2
\right)\, \mathrm{d}s\\
&\leq C\left(1+\|\bx\|^4_{\bHs} 
+\|y\|_{H}^4+ \beta^2\|\nabla y\|_{\b H}^4\right)
+C\int_0^t\left(\|\b{u}(t)\|^2_{\bHs} 
+\|\varphi(s)\|_{H}^2
+\beta\|\nabla \varphi(s)\|^2_{\b H}\right)\,\d s\\
&\qquad+C\int_0^t\left(\|\b{u}(t)\|^2_{\bHs} 
+\|\varphi(s)\|_{H}^2
+\beta\|\nabla \varphi(s)\|^2_{\b H}\right)
\left(\bu(s), G_{1}( \bu(s))\,\d  W_1(s)\right)_{\bHs}\\
&\qquad+ C\int_0^t
\left(\|\b{u}(t)\|^2_{\bHs} 
+\|\varphi(s)\|_{H}^2
+\beta\|\nabla \varphi(s)\|^2_{\b H}\right)
\left( w(s)+\varphi(s), 
G_2(\varphi(s))\,{\rm d}W_2(s)\right)_H\,
\\
& \qquad {
+\int_0^t
\|(\bu(s), G_1(\bu(s)))_\bHs\|_{\mathcal L_{HS}(U_1,\mathbb R)}^2\,\d s
+\int_0^t\|(\varphi(s)-\beta \Delta \varphi(s), G_2(\varphi(s)))_H\|_{\mathcal L_{HS}(U_2,\mathbb R)}^2\,\d s.
}
\end{align*}

{
From Assumption \ref{hyp:diffusionNS}, estimates \eqref{stima_traccia} and \eqref{star_phi} and the Young inequality we deduce
\begin{multline*}
\int_0^t \|(\bu(s), G_1(\bu(s)))_\bHs\|_{\mathcal L_{HS}(U_1,\mathbb R)}^2\,\d s
+ \int_0^t\|(\varphi(s)-\beta \Delta \varphi(s), G_2(\varphi(s)))_H\|_{\mathcal L_{HS}(U_2,\mathbb R)}^2\,\d s
\\
\le C(1+t+\|\bx\|^4_{\bHs} 
+\|y\|_{H}^4+ \beta^2\|\nabla y\|_{\b H}^4)
+ L_{G_1}\int_0^t\|\bu(s)\|^2_{\bHs}\, {\rm d}s
\\
+C\int_0^t
\left(\bu(s), G_{1}( \bu(s))\,\d  W_1(s)\right)_{\bHs}
+ C\int_0^t
\left( w(s)+\varphi(s), 
G_2(\varphi(s))\,{\rm d}W_2(s)\right)_H.
\end{multline*}
Thus we obtain the estimate 
 \begin{align*}
&\left(\|\b{u}(t)\|^2_{\bHs} 
+\|\varphi(t)\|_H^2
+\beta\|\nabla \varphi(t)\|^2_H\right)^2\\
&\qquad+ \int_0^t
\left(\|\b{u}(s)\|^2_{\bHs} 
+\|\varphi(s)\|_{H}^2
+\beta\|\nabla \varphi(s)\|^2_{\b H}\right)
\left(\nu\|\nabla\b{u}(s)\|^2_{\bHs}
+\|\varphi(s)\|_{H}^2
+\beta\|\nabla\varphi(s)\|_{\b H}^2
\right)\, \mathrm{d}s\\
&\leq C\left(1+t+\|\bx\|^4_{\bHs} 
+\|y\|_{H}^4+ \beta^2\|\nabla y\|_{\b H}^4\right)
+C\int_0^t\left(\|\b{u}(t)\|^2_{\bHs} 
+\|\varphi(s)\|_{H}^2
+\beta\|\nabla \varphi(s)\|^2_{\b H}\right)\,\d s\\
&\qquad+C\int_0^t\left(1+\|\b{u}(t)\|^2_{\bHs} 
+\|\varphi(s)\|_{H}^2
+\beta\|\nabla \varphi(s)\|^2_{\b H}\right)
\left(\bu(s), G_{1}( \bu(s))\,\d  W_1(s)\right)_{\bHs}\\
&\qquad+ C\int_0^t
\left(1+\|\b{u}(t)\|^2_{\bHs} 
+\|\varphi(s)\|_{H}^2
+\beta\|\nabla \varphi(s)\|^2_{\b H}\right)
\left( w(s)+\varphi(s), 
G_2(\varphi(s))\,{\rm d}W_2(s)\right)_H\,.
\end{align*}
}
By the Poincar\'e and Young inequalities, 
\begin{align*}
&C\int_0^t\left(\|\b{u}(t)\|^2_{\bHs} 
+\|\varphi(s)\|_{H}^2
+\beta\|\nabla \varphi(s)\|^2_{\b H}\right)\,\d s\\
&\leq \frac C{{(1 \wedge\nu)}^2}t
+
\frac14\int_0^t
\left(\|\b{u}(s)\|^2_{\bHs} 
+\|\varphi(s)\|_{H}^2
+\beta\|\nabla \varphi(s)\|^2_{\b H}\right)
\left(\nu\|\nabla\b{u}(s)\|^2_{\bHs}
+\|\varphi(s)\|_{H}^2
+\beta\|\nabla\varphi(s)\|_{\b H}^2
\right)\, \mathrm{d}s
\end{align*}
Hence, rearranging the terms above the thesis is proven.
\end{proof}

\begin{lem}
\label{lemma2}
Let either Assumptions \ref{hyp:structural}--\ref{hyp:diffusionAC} hold with $Z=\bVsd$ or Assumptions \ref{hyp:structural}--\ref{hyp:additionalAC} hold with $Z=\bHs$.
Let $(\bx,y) \in \b X \times Y$ and let 
$(\bu,\varphi)$ be the unique strong solution
to the system \eqref{AC_NS} starting from the initial state
$(\b x,y)$. Then, if $\alpha_d=1$ and $\alpha_n=0$,
there exists a constant 
$\mathfrak{C}_2=\mathfrak{C}_2(L_{G_1},L_{G_2},L_F,|\OO|,\|F\|_{\C([-1,1])})>0$,
independent of $\nu$ and $\beta$,
such that  
\begin{align}
    \label{a_priori_estimate4}
    &\beta^2\int_0^t\|\varphi(s)\|^2_{V_2}\,\d s
    \leq \mathfrak{C}_2\left(1 + 
    \|\bx\|_{\bHs}^2 + \|y\|_{H}^2 + \beta\|\nabla y\|_{\b H}^2
    +  t + \mathfrak{M_3}(t)
    +\mathfrak{M_4}(t)\right),
\end{align}
for all $t\geq0$, $\P$-almost surely,
where the martingales $\mathfrak{M}_3$ and $\mathfrak{M}_4$ are given by
\begin{align*}
    \mathfrak{M}_3(t)&:=
    \int_0^t \left(\bu(s), G_{1}( \bu(s))\,\d  W_1(s)\right)_{\bHs}\,,\quad t\geq0,\\
\mathfrak{M}_4(t)&:=\int_0^t
\left( w(s)+\varphi(s), 
G_2(\varphi(s))\,{\rm d}W_2(s)\right)_H\,, \quad t\geq0\,.
\end{align*}
\end{lem}
\begin{proof}
The thesis follows reasoning as in the proof of Proposition \ref{a_priori_est2} (see estimate \eqref{star_phi}). 
\end{proof}

\begin{lem}
\label{lemma3}
Let Assumptions \ref{hyp:structural}--\ref{hyp:additionalAC} hold with $Z=\bHs$.
Let $(\bx,y), (\widetilde{\b x}, \widetilde y) \in \b X \times Y$,
let $N\in\enne_+$, $\eta=\frac{\lambda_N\nu}2$, 
let 
$(\bu, \varphi)$ be the unique strong solution
to the system \eqref{AC_NS} starting from the initial state
$(\b x, y)$, and let 
$(\widetilde \bu,\widetilde \varphi)$ be the unique strong solution
to the system \eqref{AC_NS_nud} starting from the initial state
$(\widetilde{\b x},\widetilde y)$.
Then, setting $\gamma:=\frac{s_0+1}{s_F}>2$,
there exists a
constant $\mathfrak{C}_3=\mathfrak{C}_3(s_0, s_F, \lambda_1, L_{G_1},L_{G_2},L_F,|\OO|,\|F\|_{\C([-1,1])})>0$,
independent of $N$, $\nu$ and $\beta$, such that
\begin{align}
\label{stimaF''}
\int_0^t \|F''(\varphi(s))\|_{L^{\gamma}(\OO)}^2 \: \d s 
&\le \mathfrak{C}_3\left( 1 + \int_\OO \Psi_{s_0}(y) + t +\mathfrak{M}_5(t)\right), \\
\label{stimaF''nud}
\int_0^t \|F''(\widetilde{\varphi}(s))\|_{L^{\gamma}(\OO)}^2 \: \d s &\le \mathfrak{C}_3\left( 1 + \int_\OO \Psi_{s_0}(\widetilde y) + t +\widetilde{\mathfrak{M}}_5(t)\right),
\end{align}
for all $t\geq0$, $\P$-almost surely, 
where the martingales $\mathfrak{M}_5$ and 
$\widetilde{\mathfrak{M}}_5$ are given by 
\begin{align*}
\mathfrak{M}_5(t)&:= \int_0^t \left(\Psi'_{s_0}(\varphi(s)), G_2 (\varphi(s)) \d W_2(s) \right)_H, \qquad t \ge 0,\\
\widetilde{\mathfrak{M}}_5(t)&:= \int_0^t \left(\Psi'_{s_0}(\widetilde{\varphi}(s)), G_2 (\widetilde{\varphi}(s)) \d W_2(s) \right)_H, \qquad t \ge 0.
\end{align*}
\end{lem}
\begin{proof}
By Assumption~\ref{hyp:additionalF} and Lemma~\ref{lem:F''} we get,
by possibly updating the constant $C$,
    \begin{align*}
\int_0^t \|F''(\varphi(s)\|_{L^{\gamma}(\OO)}^\gamma \: \d s 
& \leq L_F^\gamma \int_0^t \int_\OO \left( 1 + \Psi_{s_F}(\varphi(s)) \right)^\gamma \: \d s 
\\
& \leq C\left( t +  \int_0^t \int_\OO  \Psi_{s_F}^{\gamma}(\varphi(s))  \: \d s\right) 
 = C\left(t + \int_0^t \int_\OO  \Psi_{s_0+1}(\varphi(s))  \: \d s\right) \\
& \leq C\left(1+ t + \int_\OO \Psi_{s_0}(y) + \int_0^t \left(\Psi'_{s_0}(\varphi(s)), G_2 (\varphi(s)) \,\d W_2(s) \right)_H\right).
    \end{align*}
Since $\gamma>2$ by Assumption~\ref{hyp:additionalAC}, the Young inequality yields
\[
\int_0^t \|F''(\varphi(s)\|_{L^{\gamma}(\OO)}^2 \: \d s 
\le C  \int_0^t\left( 1+ \|F''(\varphi(s)\|_{L^{\gamma}(\OO)}^\gamma \right)\: \d s. 
\]
Putting together the two estimates above the first claim follows. 
The second claim is proven analogously, by noting that the convective term in the Allen–Cahn equation (and thus the nudged term in the Navier-Stokes equations) does not have an impact on the estimate.  
\end{proof}

\begin{lem}
\label{lemma4}
Let $M$ be a real square-integrable martingale and let 
$[M]$ be its quadratic variation process. Then, for every $q>2$, 
there exists a constant $c_q>0$ such that,
for every $R,T>0$,
\[
\mathbb{P}\left( \sup_{t\ge T}\left(M(t)-t-2\right) > R  \right)\le c_q \sum_{m \ge \lfloor T\rfloor}\frac{\mathbb{E}\left[ [M](m+1)^{\frac q2}\right]
}{(R+m+2)^q}.
\]
\end{lem}

\begin{proof}
We observe that
\[
\left\{\sup_{t \ge T}\left(  M(t)-t-2\right) > R \right\} \subset \bigcup_{m \ge \lfloor T\rfloor}\left\{ \sup_{t \in [m, m+1)}\left( M(t)-t-2\right) \ge R\right\}.
\]
 On the other hand, notice that for every $m \ge 0$,
\[
\left\{ \sup_{t \in [m, m+1)}\left( M(t)-t-2\right) \ge R\right\} \subset \left\{ M^*(m+1) \ge R+m+2\right\},
\] 
where
\[
M^*(t):= \sup_{s \in [0,t] }| M(s)|.
\]
Thus, 
from the Markov inequality we infer that
\begin{align*}
\mathbb{P}\left(\sup_{t \ge T}\left( M(t)-t-2\right) > R \right) &\le \sum_{m \ge  \lfloor T\rfloor}\mathbb{P}\Big( M^*(m+1) \ge R+m+2\Big)
\le \sum_{m \ge  \lfloor T\rfloor}\frac{\mathbb{E}\left[M^*(m+1)^q\right]} {(R+m+2)^q},
\end{align*}
so that from the Burkholder-Davis-Gundy inequality
\[
\mathbb{E}\left[M^*(t)^q\right] \le c_q  \mathbb{E}\left[ \left[M \right](t)^{\frac q2}\right], \qquad t \ge 0.
\]
The proof follows.
\end{proof}

The following three propositions provide some estimates in probability for suitable functionals of the solution process. Their proof is based on the above lemmata.

\begin{prop}
    \label{prop1}
    Let either Assumptions \ref{hyp:structural}--\ref{hyp:diffusionAC} hold with $Z=\bVsd$ or Assumptions \ref{hyp:structural}--\ref{hyp:additionalAC} hold with $Z=\bHs$.
Let $(\bx,y) \in \b X \times Y$ and let 
$(\bu,\varphi)$ be the unique strong solution
to the system \eqref{AC_NS} starting from the initial state
$(\b x,y)$ and let $\mathfrak C_1$ be the constant 
given by Lemma~\ref{lemma1}.
Then, for every $q>2$, 
there exists a constant $C=C(q,\lambda_1, L_{G_1},L_{G_2},L_F,|\OO|,\|F\|_{\C([-1,1])})>0$,
independent of $\nu$ and $\beta$, such that, 
for every $R,T>0$
it holds 
\begin{multline*}
\mathbb{P}\left(\sup_{t \ge T}\left( \nu\int_0^t\|\bu(s)\|_{\bVs}^2\|\bu(s)\|_\bHs^2\,\d s - (\mathfrak{C}_1+2)\left(1 + (1+\nu^{-2}) t
+ \|\bx\|^4_{\bHs}+ \|y\|^4_{H}
+\beta^2\|\nabla y\|_\bH^4\right)\right) > R\right) \\
\le \frac{C}{(T+R)^{\frac q2-1}}\left(1+\nu^{-\frac32}+\|\bx\|_{\bHs}^{3q}+ 
\|y\|_{H}^{3q} + \beta^{\frac32q}\|\nabla y\|_{\bH}^{3q}\right).
\end{multline*}
\end{prop}
\begin{proof}
By subtracting $(t+2)$ to both sides of inequality \eqref{a_priori_estimate3} we obtain the estimate 
\[
\nu\int_0^t\|\bu(s)\|_{\bVs}^2\|\bu(s)\|_\bHs^2\,\d s - ({\mathfrak{C}_1}+2)\left(1 + (1+\nu^{-2}) t
+ 
\|\bx\|_{\bHs}^4 + \|y\|_{V_1}^4 \right) \le \mathcal{M}(t) -t -2,
    \quad\forall\,t\geq0\,,\quad\P\text{-a.s.},
\]
where $\mathcal M:=\mathfrak{C}_1( \mathfrak{M}_1  +\mathfrak{M_2})$
and $ \mathfrak{M}_1$, $\mathfrak{M_2}$ are the martingales of Lemma~\ref{lemma1}.
Thus, thanks to Lemma~\ref{lemma4} we infer, for any $T \ge 0$,
\begin{multline}
\label{P1proof}
\mathbb{P}\left(\sup_{t \ge T}\left( \nu\int_0^t\|\bu(s)\|_{\bVs}^2\|\bu(s)\|_\bHs^2\,\d s - (\mathfrak{C}_1+2)\left(1 +(1+\nu^{-2}) t + \|\bx\|^4_{\bHs}+ \|y\|^4_{H}
+\beta^2\|\nabla y\|_\bH^4\right)\right) > R\right)
     \\
    \le 
     \mathbb{P}\left(\sup_{t \ge T}
     \left(\left(\mathcal{M}(t)\right) - t-2\right) \ge R\right)
\le c_q \sum_{m \ge \lfloor T\rfloor}\frac{\mathbb{E}\left[ [\mathcal{M}](m+1)^{\frac q2}\right]}{(R+m+2)^q}.
\end{multline}
We estimate the $\frac q2$-moment of the quadratic variation process $\{[\mathcal{M}](t)\}_{t \ge 0}$ as follows.
By the Young inequality we estimate, for any $t \ge 0$, we have
\begin{multline*}
[\mathfrak{M}_1](t) \le \int_0^t \left( {1+}\|\bu(s)\|^2_{\bHs}+ \|\varphi(s)\|^2_{H}
+\beta\|\nabla\varphi(s)\|_\bH^2\right)^2
\|(\bu(s), G_1(\bu(s)))_{\bHs}\|^2_{\mathcal{L}_{HS}(U_1, \mathbb{R})}\, {\rm d}s
\\
\le L_{G_1}\int_0^t \left({1+} \|\bu(s)\|^2_{\bHs}+ \|\varphi(s)\|^2_{H}
+\beta\|\nabla\varphi(s)\|_\bH^2
\right)^2\|\bu(s)\|_{\bHs}^2\, {\rm d}s\\
\le C\int_0^t \left({1+} \|\bu(s)\|^6_{\bHs}+
\|\varphi(s)\|^6_{H}
+\beta^3\|\nabla\varphi(s)\|_\bH^6\right)\, {\rm d}s,
\end{multline*}
while by reasoning as in the proof of Proposition \ref{a_priori_est2} and by exploiting the Young inequality, for any $t \ge 0$ we have 
\begin{multline*}
[\mathfrak{M}_2](t)
\le \int_0^t \left({1+}\|\bu(s)\|^2_{\bHs}+
\|\varphi(s)\|^2_{H}
+\beta\|\nabla\varphi(s)\|_\bH^2\right) 
\|(\varphi(s)+w(s), G_2(\varphi(s)))_H\|_{\mathcal L_{HS}(U_2,\mathbb R)}^2\, {\rm d} s 
\\
\le C \int_0^t \left( {1+}\|\bu(s)\|^2_{\bHs}+
\|\varphi(s)\|^2_{H}
+\beta\|\nabla\varphi(s)\|_\bH^2\right)(1 +\beta^2\|\nabla\varphi(s)\|_{\bH}^4)\, {\rm d} s\\
\le C \int_0^t \left({1+}\|\bu(s)\|^6_{\bHs}+\|\varphi(s)\|^6_{H}
+\beta^3\|\nabla\varphi(s)\|_\bH^6\right)\, {\rm d} s.
\end{multline*}
Hence, noting that
\[
[\mathcal{M}](t)\leq 2\mathfrak{C}_1^2 \left( [\mathfrak{M}_1](t)+[\mathfrak{M}_2](t)\right), \qquad t \ge 0,
\]
by the H\"older inequality and Proposition~\ref{a_priori_est2}
we infer that
\begin{multline*}
\mathbb{E} \left[[\mathcal{M}](t)^{\frac q2}\right]\le 
C \left( \mathbb{E} \left[[\mathfrak{M}_1](t)^{\frac q2}\right]+
\mathbb{E} \left[[\mathfrak{M}_2](t)^{\frac q2}\right]
\right)\\
\leq C t^{\frac{q}2-1}
 \mathbb{E}\left[ \int_0^t\left({1+}\|\bu(s)\|^{3q}_{\bHs}+ \|\varphi(s)\|^{3q}_{H}
+\beta^{\frac32q}\|\nabla\varphi(s)\|_\bH^{3q}\right)\, {\rm d}s\right]\\
 \le C(t+1)^{\frac{q}{2}}
 \left(1+\nu^{-\frac32q} + \|\bx\|^{3q}_{\bHs}+
 \|y\|^{3q}_{H}
+\beta^{\frac32q}\|\nabla y\|_\bH^{3q}\right),
\end{multline*}
where we have updated the constant $C$ form line to line, as usual.
From the above estimate and \eqref{P1proof} we infer
\begin{multline*}
\mathbb{P}\left(\sup_{t \ge T}\left( \nu\int_0^t\|\bu(s)\|_{\bVs}^2\|\bu(s)\|_\bHs^2\,\d s - (\mathfrak{C}_1+2)\left(1 +(1+\nu^{-2})t  + \|\bx\|^4_{\bHs}+ \|y\|^4_{H}
+\beta^2\|\nabla y\|_\bH^4\right)\right) > R\right)\\
\le C\left(1+\nu^{-\frac32q} + \|\bx\|^{3q}_{\bHs}+
 \|y\|^{3q}_{H}
+\beta^{\frac32q}\|\nabla y\|_\bH^{3q}\right)
\sum_{m \ge \lfloor T\rfloor}\frac{(m+2)^{\frac q2}}{(R+m+2)^q}
\\
\le C\left(1+\nu^{-\frac32q} + \|\bx\|^{3q}_{\bHs}+
 \|y\|^{3q}_{H}
+\beta^{\frac32q}\|\nabla y\|_\bH^{3q}\right)
\sum_{m \ge \lfloor T\rfloor}\frac{1}{(R+m+2)^{\frac q2}}\\
\le C \frac{\left(1+\nu^{-\frac32q} + \|\bx\|^{3q}_{\bHs}+
 \|y\|^{3q}_{H}
+\beta^{\frac32q}\|\nabla y\|_\bH^{3q}\right)}{(T+R)^{\frac q2 -1}},
\end{multline*}
and this concludes the proof.
\end{proof}

\begin{prop}
    \label{prop2}
    Let either Assumptions \ref{hyp:structural}--\ref{hyp:diffusionAC} hold with $Z=\bVsd$ or Assumptions \ref{hyp:structural}--\ref{hyp:additionalAC} hold with $Z=\bHs$.
Let $(\bx,y) \in \b X \times Y$ and let 
$(\bu,\varphi)$ be the unique strong solution
to the system \eqref{AC_NS} starting from the initial state
$(\b x,y)$ and let $\mathfrak C_2$ be the constant 
given by Lemma~\ref{lemma2}.
Then, if $\alpha_d=1$ and $\alpha_n=0$, for every $q>2$, 
there exists a constant $C=C(q,\lambda_1, L_{G_1},L_{G_2},L_F,|\OO|,\|F\|_{\C([-1,1])})>0$,
independent of $\nu$ and $\beta$, such that, 
for every $R,T>0$
it holds 
\begin{multline*}
\mathbb{P}\left(\sup_{t \ge T}\left( \beta^2 \int_0^t\|\varphi(s)\|_{V_2}^2\,\d s 
- (\mathfrak{C}_2+2)\left(1 +\|\bx\|^2_{\bHs}+ 
\|y\|_{H}^2 + \beta\|\nabla y\|_{\b H}^2
+t \right) \right) > R)\right)\\
\le \frac{C}{(T+R)^{\frac q2-1}}\left(1+\nu^{-q} + \|\bx\|_{\bHs}^{2q}+
\|y\|_{H}^{2q}
+\beta^q\|\nabla y\|_\bH^{2q}\right).
\end{multline*}
\end{prop}

\begin{proof}
We proceed as in the proof of Proposition \ref{prop1}.
By subtracting $(t+2)$ to both sides of inequality \eqref{a_priori_estimate4} we obtain the estimate 
\[
\beta^2\int_0^t\|\varphi(s)\|_{V_2}^2\,\d s - 
({\mathfrak{C}_2}+2)\left(1 + \|\bx\|_{\bHs}^2 + 
\|y\|_{H}^2 + \beta\|\nabla y\|_\bH^2 
+  t \right) \le \mathcal{M}(t) -t -2,
    \quad\forall\,t\geq0\,,\quad\P\text{-a.s.},
\]
where $\mathcal M:=\mathfrak{C}_2( \mathfrak{M}_3  +\mathfrak{M}_4)$
and $ \mathfrak{M}_3$, $\mathfrak{M}_4$ are the martingales of Lemma~\ref{lemma2}.
Thus, thanks to Lemma~\ref{lemma4} we infer, for any $T \ge 0$,
\begin{multline}
\label{P2proof}
    \mathbb{P}\left(\sup_{t \ge T}\left( 
    \beta^2\int_0^t\|\varphi(s)\|_{V_2}^2\,\d s - 
({\mathfrak{C}_2}+2)\left(1 + \|\bx\|_{\bHs}^2 + 
\|y\|_{H}^2 + \beta\|\nabla y\|_\bH^2 
+  t \right)\right) > R)\right)
     \\
    \le 
     \mathbb{P}\left(\sup_{t \ge T}
     \left(\left(\mathcal{M}(t)\right) - t-2\right) > R\right)
\le C_q \sum_{m \ge \lfloor T\rfloor}\frac{\mathbb{E}\left[ [\mathcal{M}](m+1)^{\frac q2}\right]}{(R+m+2)^q}.
\end{multline}
By reasoning as in the proof of Proposition \ref{prop1}, 
it is immediate to see that
\begin{align*}
 \mathbb{E} \left[[\mathfrak{M}_3](t)^{\frac q2}\right]  
 &\le C(t+1)^{\frac{q}{2}}\left(1+\nu^{-\frac q2} + \|\bx\|^{q}_{\bHs}+\|y\|^{q}_{H}
 +\beta^{\frac q2}\|\nabla y\|_\bH^q\right),\\
\mathbb{E} \left[[\mathfrak{M}_4](t)^{\frac q2}\right]  
 &\le C(t+1)^{\frac{q}{2}}\left(1+\nu^{-q} + \|\bx\|^{2q}_{\bHs}+\|y\|^{2q}_{H}
 +\beta^{q}\|\nabla y\|_\bH^{2q}\right),
\end{align*}
hence
\[
\mathbb{E} \left[[\mathcal{M}](t)^{\frac q2}\right]  
  \le C(t+1)^{\frac{q}{2}}\left(1+\nu^{-q} + \|\bx\|^{2q}_{\bHs}+\|y\|^{2q}_{H}
 +\beta^{q}\|\nabla y\|_\bH^{2q}\right).
\]
From the above estimate and \eqref{P2proof} we infer,
by proceeding as in the proof of Proposition~\ref{prop1}, that
\begin{multline*}
\mathbb{P}\left(\sup_{t \ge T}\left( 
    \beta^2\int_0^t\|\varphi(s)\|_{V_2}^2\,\d s - 
({\mathfrak{C}_2}+2)\left(1 + \|\bx\|_{\bHs}^2 + 
\|y\|_{H}^2 + \beta\|\nabla y\|_\bH^2 
+  t \right)\right) > R)\right)\\
\le C\left(1+\nu^{-q} + \|\bx\|^{2q}_{\bHs}+\|y\|^{2q}_{H}
 +\beta^{q}\|\nabla y\|_\bH^{2q}\right)
 \sum_{m \ge \lfloor T\rfloor}\frac{1}{(R+m+2)^{\frac q2}}\\
\le C\frac{\left(1+\nu^{-q} + \|\bx\|^{2q}_{\bHs}+\|y\|^{2q}_{H}
 +\beta^{q}\|\nabla y\|_\bH^{2q}\right)}{(T+R)^{\frac q2 -1}},
\end{multline*}
and this concludes the proof.
\end{proof}

\begin{prop}
    \label{prop3}
    Let Assumptions \ref{hyp:structural}--\ref{hyp:additionalAC} hold with $Z=\bHs$.
Let $(\bx,y), (\widetilde{\b x}, \widetilde y) \in \b X \times Y$,
let $N\in\enne_+$, $\eta=\frac{\lambda_N\nu}2$, 
let $(\bu,\varphi)$ be the unique strong solution
to the system \eqref{AC_NS} starting from the initial state
$(\b x,y)$,
let 
$(\widetilde \bu,\widetilde \varphi)$ be the unique strong solution
to the system \eqref{AC_NS_nud} starting from the initial state
$(\widetilde{\b x},\widetilde y)$,
and let $\mathfrak C_3$ be the constant 
given by Lemma~\ref{lemma3}.
Then, setting $\gamma:=\frac{s_0+1}{s_F}>2$,
for every $q>2$
there exists a constant $C=C(q,s_0, s_F, \lambda_1, L_{G_1},L_{G_2},L_F,|\OO|,\|F\|_{\C([-1,1])})>0$,
independent of $N$, $\nu$ and $\beta$, such that, 
for every $R,T>0$
it holds
\begin{multline*}
\mathbb{P}\left(\sup_{t \ge T}\left(  \int_0^t  
\left(\|F''({\varphi}(s))\|_{L^{\ell}(\OO)}^2+ \|F''(\widetilde{\varphi}(s))\|_{L^{\ell}(\OO)}^2\right) \: \d s
\right.\right.\\
\left.\left.- 2(\mathfrak{C}_3+1)\left( 1 + \int_\OO \Psi_{s_0}(y) + \int_\OO \Psi_{s_0}(\widetilde y)+ t\right)\right) > R)\right)
\le \frac{C}{(T+R)^{\frac q2-1}}.
\end{multline*}
\end{prop}
\begin{proof}
We proceed as in the proof of Proposition \ref{prop1}.
We subtract $\left(\frac t2 +1\right)$ to both sides of inequalities \eqref{stimaF''}, \eqref{stimaF''nud} and sum them up. We obtain the estimate 
\begin{multline*}
\int_0^t \left( \|F''(\varphi(s))\|_{L^{\ell}(\OO)}^2+ \|F''(\widetilde{\varphi}(s))\|_{L^{\ell}(\OO)}^2 \right) \,\d s -\\ 2({\mathfrak{C}_3}+1)\left(1 + \int_\OO \Psi_{s_0}(y)  +  \int_\OO \Psi_{s_0}(\widetilde y) + t \right) \le  \mathcal{M}(t) -t -2,
    \quad\forall\,t\geq0\,,\quad\P\text{-a.s.},
\end{multline*}
where $\mathcal M:=\mathfrak{C}_3( \mathfrak{M}_5  +\mathfrak{M_6})$
and $ \mathfrak{M}_5$, $\mathfrak{M_6}$ are the martingales of Lemma~\ref{lemma3}.
From Lemma~\ref{lemma4} we infer, for any $T \ge 0$,
\begin{multline}
\label{P3proof}
\mathbb{P}\left(\sup_{t \ge T}\left(  \int_0^t  
\left(\|F''({\varphi}(s))\|_{L^{\ell}(\OO)}^2+ \|F''(\widetilde{\varphi}(s))\|_{L^{\ell}(\OO)}^2\right) \: \d s
\right.\right.\\
\left.\left.- 2(\mathfrak{C}_3+1)\left( 1 + \int_\OO \Psi_{s_0}(y) + \int_\OO \Psi_{s_0}(\widetilde y)+ t\right)\right) > R)\right)\\
    \le 
     \mathbb{P}\left(\sup_{t \ge T}\left[\left(\mathcal{M}(t)\right) - t-2\right] > R\right)
\le C_q \sum_{m \ge \lfloor T\rfloor}\frac{\mathbb{E}\left[ [\mathcal{M}](m+1)^{\frac q2}\right]}{(R+m+2)^q}.
\end{multline}
Noting that by Assumption~\ref{hyp:additionalAC}
one has
\[
[\mathfrak M_5](t) + [\mathfrak M_6](t)
\leq 2L_{G_2}|\OO|t,
\]
we get
\[
\mathbb{E}\left[[\mathcal{M}](t)^{\frac q2}\right]
 \le C(t+1)^{\frac{q}{2}}, \qquad t \ge 0.
\]
Therefore, together with \eqref{P3proof} we obtain
\begin{multline*}
\mathbb{P}\left(\sup_{t \ge T}\left(  \int_0^t  
\left(\|F''({\varphi}(s))\|_{L^{\ell}(\OO)}^2+ \|F''(\widetilde{\varphi}(s))\|_{L^{\ell}(\OO)}^2\right) \: \d s
\right.\right.\\
\left.\left.- 2(\mathfrak{C}_3+1)\left( 1 + \int_\OO \Psi_{s_0}(y) + \int_\OO \Psi_{s_0}(\widetilde y)+ t\right)\right) > R)\right)\\
\le C \sum_{m \ge \lfloor T\rfloor}\frac{(m+2)^{\frac q2}}{(R+m+2)^q}
\le C \sum_{m \ge \lfloor T\rfloor}\frac{1}{(R+m+2)^{\frac q2}}
\le \frac{C}{(T+R)^{\frac q2 -1}},
\end{multline*}
and this concludes the proof.
\end{proof}

\section*{Data availability statement}
No new data were created or analysed in this study. Data sharing is not applicable to this article.

\printbibliography

\end{document}